# On the Euler+Prandtl expansion for the Navier-Stokes equations

Igor Kukavica, Trinh T. Nguyen, Vlad Vicol, and Fei Wang

ABSTRACT. We establish the validity of the Euler+Prandtl approximation for solutions of the Navier-Stokes equations in the half plane with Dirichlet boundary conditions, in the vanishing viscosity limit, for initial data which are analytic only near the boundary, and Sobolev smooth away from the boundary. Our proof does not require higher order correctors, and works directly by estimating an $L^1$-type norm for the vorticity of the error term in the expansion Navier-Stokes−(Euler+Prandtl). An important ingredient in the proof is the propagation of local analyticity for the Euler equation, a result of independent interest.

CONTENTS

## 1. Introduction

In this paper, we consider the Navier-Stokes system

$$\partial_t \boldsymbol{u}^{\rm NS} - \epsilon^2 \Delta \boldsymbol{u}^{\rm NS} + \boldsymbol{u}^{\rm NS} \cdot \nabla \boldsymbol{u}^{\rm NS} + \nabla p^{\rm NS} = 0 \tag{1.1}$$

$$\operatorname{div} \boldsymbol{u}^{\rm NS} = 0 \tag{1.2}$$

on the domain $\mathbb{H} = \mathbb{T} \times \mathbb{R}_+ = \{(x,y) \in \mathbb{T} \times \mathbb{R} \colon y \geq 0\}$, where $\mathbb{T} = [-\pi, \pi]$, with the no-slip boundary condition

$$\boldsymbol{u}^{\rm NS}|_{y=0} = 0 \tag{1.3}$$

and with an incompressible initial datum

$$\boldsymbol{u}^{\rm NS}|_{t=0} = u_0^{\rm NS}\,. \tag{1.4}$$

Throughout the paper, we denote the kinematic viscosity by $\epsilon^2$. Our goal is to establish, with a concise proof, the *Euler+Prandtl approximation* for solutions of (1.1)–(1.4) in the vanishing viscosity limit $\epsilon \to 0$, for initial data that are analytic only near the boundary of the domain, and are Sobolev smooth away from the boundary.

**1.1. Previous results.** One of the fundamental problems in mathematical fluid dynamics is to determine whether the solutions of the Navier-Stokes equations (1.1)–(1.2) converge to the solution of the Euler equations

$$\partial_t \boldsymbol{u}^{\mathrm{E}} + \boldsymbol{u}^{\mathrm{E}} \cdot \nabla \boldsymbol{u}^{\mathrm{E}} + \nabla p^{\mathrm{E}} = 0 \tag{1.5}$$

$$\operatorname{div} \boldsymbol{u}^{\mathrm{E}} = 0 \tag{1.6}$$

in the inviscid limit $\epsilon \to 0$. In [30], Kato showed that the inviscid limit holds in the energy norm $L^\infty(0, T, L^2(\mathbb{H}))$ if and only if the energy dissipation in a thin layer of size $\epsilon^2$ vanishes as $\epsilon \to 0$, i.e.,

$$\lim_{\epsilon \to 0} \epsilon^2 \int_0^T \int_{\{y \lesssim \epsilon^2\}} |\nabla \boldsymbol{u}^{\mathrm{NS}}|^2 dx dy dt = 0. \tag{1.7}$$

We refer the reader to [3, 6, 7, 31, 33, 49, 57, 58, 61] for refinements and extensions based on Kato's original argument; cf. also the recent review [48]. These results assume explicit properties that the sequence of Navier-Stokes solutions must obey on $[0, T]$ as $\epsilon \to 0$ in order for them to have a strong $L^\infty_t L^2_x$ Euler limit. On the other hand, verifying these conditions based on the knowledge of the initial datum only is in general an outstanding open problem. We emphasize that to date, even the question of whether the weak $L^2_t L^2_x$ inviscid limit holds (against test functions compactly supported in the interior of the domain), remains open. Conditional results have been established recently in terms of interior structure functions [9, 11], or in terms of interior vorticity concentration measures [8].

In his seminal 1908 paper, Prandtl postulated that the solution of the Navier-Stokes equations can be written as

$$\boldsymbol{u}^{\mathrm{NS}}(x, y, t) = \boldsymbol{u}^{\mathrm{E}}(x, y, t) + \left(\widetilde{u}^{\mathrm{P}}\left(x, \frac{y}{\epsilon}, t\right), \epsilon \bar{v}^{\mathrm{P}}\left(x, \frac{y}{\epsilon}, t\right)\right) + \mathcal{O}(\epsilon), \tag{1.8}$$

where $\boldsymbol{u}^{\mathrm{E}}$ denotes the solution of the Euler equations and $\widetilde{u}^{\mathrm{P}}$, $\bar{v}^{\mathrm{P}}$ are components of the solution of the Prandtl boundary layer equations (see (2.10) below). While the well-posedness [1, 10, 18, 28, 35, 37, 40, 42, 50, 54, 55] and the ill-posedness [14, 26, 19, 41] regimes for the Prandtl equations are by now well-understood, establishing the validity of the expansion (1.8) presents a number of outstanding challenges.

In [55, 56], Sammartino-Caflisch establish the validity of the Prandtl expansion and hence the strong inviscid limit in the energy norm, for well-prepared and analytic initial data $u_0$, in the sense that $u_0$ satisfies the Prandtl ansatz (1.8) at time $t = 0$, and are analytic in both the $x$ and $y$ variables on the entire half space. They construct solutions of Euler and Prandtl in suitable analytic spaces in $x$ and $y$, carefully analyze the error terms in the expansion (1.8), and show that they remain $\mathcal{O}(\epsilon)$ for an $\mathcal{O}(1)$ time interval by an abstract Cauchy-Kowalevski theorem. This strategy has been proven successful for treating the case of a channel [43, 34] and the exterior of a disk [5].

In [47], Maekawa established the validity of the expansion (1.8) for Sobolev smooth initial vorticity that is compactly supported away from the boundary, by using the vorticity boundary condition in [2, 46] and controlling the weak interaction between the Prandtl solutions near the boundary and the Euler solution far away from the boundary.

We refer the reader to [59] for an energy based proof of the Caflisch-Sammartino result, and [12, 13] for a proof of Maekawa's result in 2D and 3D respectively, which relies solely on energy methods.

Recently, in [53], Nguyen and the second author establish the strong inviscid limit in $L^\infty(0, T, L^2(\mathbb{H}))$ for analytic initial data, and for the first time, avoid completely the use of Prandtl boundary layer correctors (1.8). Instead, they appeal to the vorticity formulation, give precise pointwise bounds for the associated Green's function, and work in a suitable analytic boundary-layer function spaces that control the pointwise behavior of solutions of the Navier-Stokes equations. In this paper we use the pointwise estimates for the Green function of the Stokes problem from [53]; cf. Lemma 7.2 below.

In [38, 39], the first and the last two authors established the strong inviscid limit in the energy norm, for initial data that is only analytic close to the boundary of the domain, and has finite Sobolev regularity in the complement (see also [60] in the 3D case). These works thus close the gap between the Sammartino-Caflisch [55, 56], which assumes the analyticity on the entire half-plane, and the Maekawa [47] results, which assumes that the initial vorticity vanishes identically near the boundary. Up to now, the class of initial data considered in [38, 39] appears to be the largest class of initial data that the strong inviscid limit is known to hold, in the absence of structural or symmetry assumptions. Note that neither [53] nor [38, 39] establish the validity of the expansion (1.8), which is the main result of this paper.

Recently in [16, 17], Gerard-Varet, Maekawa, and Masmoudi improved the classical results of Sammartino-Caflisch to Gevrey perturbations in $x$ and Sobolev perturbation in $y$ for shear flows of the Prandtl type, when the Prandtl shear flow is both monotonic and concave. Lastly, we mention that the vanishing viscosity limit is also known to hold in the presence of certain symmetry assumptions on the initial data, which is maintained by the flow; see e.g. [4, 27, 32, 44, 45, 48, 51, 52, 20] and references therein. Also, the very recent works [15, 25, 24, 29] establish the vanishing viscosity limit and the validity of the Prandtl expansion for the stationary Navier-Stokes equation, in certain regimes.

It is worth noting that in all the above cases the Prandtl expansion (1.8) is valid, and thus the Kato criterion (1.7) holds. However, in general there is a large discrepancy between the question of the vanishing viscosity limit in the energy norm, and the problem of the validity of the Prandtl expansion.

In the negative direction of the Prandtl asymptotic expansion, we refer the reader to the works [**21, 22, 23**] of Grenier and Nguyen, which show that the Prandtl expansion (1.8) is in general false at the level of Sobolev regularity.

### 1.2. The present paper.

The main purpose of this paper is two-fold.

First, we establish the Prandtl asymptotic expansion (1.8) for initial data that is only analytic near the boundary, and is Sobolev regular in the complement. When compared to [**38, 39**], the main difficulty here is that the Euler equation is not a priori well-suited for propagating regularity which is analytic near the boundary of the domain, and only Sobolev away from the domain. The main reasons are that the pressure is nonlocal and the equation is not parabolic. This essential fact is established in Theorem 5.1 below. The proof consists of three steps. First, we obtain the analyticity of the Euler solution with respect to the operators $y\partial_y$ and $\partial_x$ (i.e., in an analytic wedge), by approximating the Euler solution via the Navier-Stokes solutions as in [**38, 39**]. Since the Euler data is uniformly analytic up to the boundary, it belongs to the initial space required by [**38**]. In the second step, we use Montel's theorem for normal families, to obtain that the family of the Navier-Stokes solutions, which are analytic in a wedge, have a subsequence which converges to the solution of the Euler equation, which is then analytic in a wedge. In the third step we bootstrap the analyticity to uniform by using the following strategy. The solution of the Euler equations is analytic uniformly on any line which is at a positive distance from the boundary. This provides analyticity of $v^{\mathrm{E}}$ on every such line. Note that, in addition, $v^{\mathrm{E}} = 0$ on the boundary $\partial\mathbb{H}$. Therefore, we may perform a localized analytic energy proof, which takes advantage of the boundary condition on the lower boundary and the uniform interior analyticity strictly inside the domain to propagate the local analyticity forward in time.

Secondly, we note that in the previous works where the Prandtl expansion was justified, a further asymptotic expansion of the error term was used, by considering correctors given by the linearization of Navier-Stokes about the Euler and Prandtl solutions, with suitable boundary conditions. Our main improvement is to obtain the convergence directly, without resorting to further expansions, by using the $L^1$ based analytic spaces from [**53, 38, 39**]. As a consequence of this simpler approximation procedure, our main result requires fewer compatibility conditions between the Euler, Navier-Stokes, and Prandtl initial data, when compared to [**55, 56**].

The paper is structured as follows. In Section 2, we introduce the Euler+Prandtl approximation of Navier-Stokes, at the level of the vorticity. The main theorem concerning the expansion (1.8) is stated in Section 3, along with a corollary, which states that the high order $\mathcal{O}(\epsilon)$ estimate on the error also holds in the uniform norm. In Section 4, we recall the equation for the error and introduce the necessary norms, along with some preliminary results. Sections 5 and 6 contain the necessary analytic bounds for the Euler and Prandtl equations, respectively. Sections 7.2 and 7.3 contain the analytic and Sobolev estimates needed in the proof of the main result. The proof of the main theorem is then provided in Section 9, while the proof of the corollary are given in Section 10.

## 2. The Euler+Prandtl approximation in the vorticity form

In order to describe the Euler+Prandtl approximation of solutions to the Navier-Stokes equation, it is convenient to work with the vorticity formulations of the Navier-Stokes, Euler, and Prandtl equations. We describe these next.

**The Navier-Stokes vorticity.** We denote the components of the Navier-Stokes velocity as $\boldsymbol{u}^{\mathrm{NS}} = (u^{\mathrm{NS}}, v^{\mathrm{NS}})$ and let the associated vorticity be given by

$$\omega^{\mathrm{NS}} = \nabla^{\perp} \cdot \boldsymbol{u}^{\mathrm{NS}} = \partial_x v^{\mathrm{NS}} - \partial_y u^{\mathrm{NS}} .$$

The Navier-Stokes vorticity satisfies

$$\partial_t \omega^{\mathrm{NS}} - \epsilon^2 \Delta \omega^{\mathrm{NS}} = -(u^{\mathrm{NS}} \partial_x + v^{\mathrm{NS}} \partial_y) \omega^{\mathrm{NS}} .$$

in $\mathbb{H}$, with the boundary condition given by (cf. [**2, 46, 47**])

$$\epsilon^2 (\partial_y + |\partial_x|) \omega^{\mathrm{NS}} = \partial_y \Delta^{-1} \left( (u^{\mathrm{NS}} \partial_x + v^{\mathrm{NS}} \partial_y) \omega^{\mathrm{NS}} \right) |_{y=0} .$$

**The Euler vorticity**. Away from the boundary $\{y = 0\}$, that is for $y \gtrsim \epsilon$, the Navier-Stokes vorticity shall be shown to be well-approximated by the Euler vorticity, which we denote as

$$\omega^{\mathrm{E}} = \nabla^{\perp} \cdot \boldsymbol{u}^{\mathrm{E}} = \partial_x v^{\mathrm{E}} - \partial_y u^{\mathrm{E}} .$$

Here, $\boldsymbol{u}^{\mathrm{E}} = (u^{\mathrm{E}}, v^{\mathrm{E}})$ is the smooth solution of the Euler equations in $\mathbb{H}$, i.e., (1.1)–(1.2) with $\epsilon = 0$, with the initial condition

$$(u^{\mathrm{E}}, v^{\mathrm{E}})|_{t=0} = (u_0^{\mathrm{E}}, v_0^{\mathrm{E}}) \tag{2.1}$$

specified below and the boundary condition

$$v^{\mathrm{E}}|_{y=0} = 0\,. \tag{2.2}$$

It is convenient to denote by $U^{\mathrm{E}}$ and $P^{\mathrm{E}}$ the trace of the Euler tangential flow and pressure on $\partial\mathbb{H}$, i.e.,

$$U^{\mathrm{E}}(t,x) = u^{\mathrm{E}}(t,x,0) \qquad \text{and} \qquad P^{\mathrm{E}}(t,x) = p^{\mathrm{E}}(t,x,0). \tag{2.3}$$

**The Prandtl vorticity.** Close to the boundary $\{y = 0\}$, that is for $y \lesssim \epsilon$, the Navier-Stokes vorticity is shown below to be well-approximated by the *total boundary layer vorticity* defined in (2.11). We recall that the Prandtl equations for the velocity field $(u^{\mathrm{P}}, \epsilon v^{\mathrm{P}})$, which are functions of $t$, $x$, and the fast normal variable[1]

$$Y = \frac{y}{\epsilon}\,,$$

read as

$$(\partial_t - \partial_{YY})u^{\mathrm{P}} + u^{\mathrm{P}}\partial_x u^{\mathrm{P}} + v^{\mathrm{P}}\partial_Y u^{\mathrm{P}} = -\partial_x P^{\mathrm{E}}\,, \tag{2.4}$$

$$v^{\mathrm{P}} = -\int_0^Y \partial_x u^{\mathrm{P}}\, dY'\,, \tag{2.5}$$

for $(x, Y) \in \mathbb{H}$. The boundary conditions for $u^{\mathrm{P}}$ are

$$u^{\mathrm{P}}|_{Y=0} = 0 \quad \text{and} \quad u^{\mathrm{P}}|_{Y\to\infty} = U^{\mathrm{E}}\,. \tag{2.6}$$

The classical Prandtl vorticity, defined as

$$\Omega^{\mathrm{P}} = \partial_Y u^{\mathrm{P}}\,, \tag{2.7}$$

satisfies the equation

$$\partial_t\Omega^{\mathrm{P}} - \partial_Y^2\Omega^{\mathrm{P}} + u^{\mathrm{P}}\partial_x\Omega^{\mathrm{P}} + v^{\mathrm{P}}\partial_Y\Omega^{\mathrm{P}} = 0 \tag{2.8}$$

in $\mathbb{H}$, with the boundary conditions

$$\partial_y\Omega^{\mathrm{P}}|_{Y=0} = \partial_x P^{\mathrm{E}} \qquad \text{and} \qquad \Omega^{\mathrm{P}}|_{Y\to\infty} = 0\,. \tag{2.9}$$

The Prandtl velocity component $u^{\mathrm{P}}$ may then be computed from the vorticity as $u^{\mathrm{P}}(x,Y) = \int_0^Y \Omega^{\mathrm{P}}(x,Y')dY'$. The boundary layer velocity vector is then given by $(\widetilde{u}^{\mathrm{P}}, \epsilon\bar{v}^{\mathrm{P}})$, where

$$\widetilde{u}^{\mathrm{P}} = u^{\mathrm{P}} - U^{\mathrm{E}} \qquad \text{and} \qquad \bar{v}^{\mathrm{P}} = \int_Y^\infty \partial_x\widetilde{u}^{\mathrm{P}} dY'\,. \tag{2.10}$$

We introduce the total boundary layer vorticity by

$$\omega^{\mathrm{P}} = (-\partial_y, \partial_x)\cdot(\widetilde{u}^{\mathrm{P}}, \epsilon\bar{v}^{\mathrm{P}}) = \epsilon\partial_x\bar{v}^{\mathrm{P}} - \frac{1}{\epsilon}\partial_Y\widetilde{u}^{\mathrm{P}} = \epsilon\partial_x\bar{v}^{\mathrm{P}} - \frac{1}{\epsilon}\Omega^{\mathrm{P}}\,. \tag{2.11}$$

**The Euler+Prandtl expansion.** In terms of the vorticity, the Euler+Prandtl expansion of the Navier-Stokes solution is

$$\omega^{\mathrm{NS}} = \omega^{\mathrm{E}} + \omega^{\mathrm{P}} + \epsilon\omega_{\mathrm{e}}\,, \tag{2.12}$$

where $\omega_{\mathrm{e}}$ is the *error vorticity*. To prove the validity of the Euler+Prandtl expansion amounts to showing that the error vorticity $\omega_{\mathrm{e}}$ is $\mathcal{O}(1)$ with respect to $\epsilon$ uniformly in time, in a suitable norm in space. We achieve this in Theorem 3.1 below. Since we prove the validity of the expansion uniformly in time which is $\epsilon$-independent, the initial data for the Navier-Stokes equation has to be compatible with (2.12), which we explain next.

**Compatible initial data.** By compatible initial data $\boldsymbol{u}_0^{\mathrm{NS}} = (u_0^{\mathrm{NS}}, v_0^{\mathrm{NS}})$ and $\boldsymbol{u}_0^{\mathrm{E}} = (u_0^{\mathrm{E}}, v_0^{\mathrm{E}})$ we mean that

$$u_0^{\mathrm{NS}}(x,y) = u_0^{\mathrm{E}}(x,y) + \widetilde{u}_0^{\mathrm{P}}(x,Y) + \epsilon u_{\mathrm{e}0}(x,y)\,, \tag{2.13}$$

$$v_0^{\mathrm{NS}}(x,y) = v_0^{\mathrm{E}}(x,y) + \epsilon\bar{v}_0^{\mathrm{P}}(x,Y) + \epsilon v_{\mathrm{e}0}(x,y)\,, \tag{2.14}$$

where $(\widetilde{u}_0^{\mathrm{P}}, \epsilon\bar{v}_0^{\mathrm{P}})$ are defined from the Prandtl initial datum $u_0^{\mathrm{P}}$ via (2.10), and the error velocity $(u_{\mathrm{e}0}, v_{\mathrm{e}0})$ is incompressible and satisfies boundary conditions which ensure that $u_0^{\mathrm{NS}}(x,0) = v_0^{\mathrm{NS}}(x,0) = 0$, namely $u_{\mathrm{e}0}(x,0) = 0$ and

[1] Throughout the paper, we use the vertical spatial variable $Y$ for the Prandtl variables, and $y$ for all others.

$v_{e0}(x,0) = -\int_0^\infty \partial_x \widetilde{u}_0^{P}(x,Y)dY$. In addition to (2.13)–(2.14), we require that $\omega_{e0}$ is $\mathcal{O}(\epsilon)$ in a suitable norm which is $L^\infty$ based in $x$ and $L^1$ based in $y$ (cf. (7.1) below).

A concrete example for compatible initial data is as follows.[2] The initial data for the modified Prandtl velocity components $\widetilde{u}^{P}$ and $\bar{v}^{P}$ (cf. (2.10)) may be taken as

$$\widetilde{u}_0^{P} = U_0^{E}(x)\varphi'(Y) \qquad \text{and} \qquad \bar{v}_0^{P} = -\partial_x U_0^{E}(x)\varphi(Y)\,, \tag{2.15}$$

where $\varphi$ is a uniformly analytic function which decays sufficiently fast as $Y \to \infty$, along with its derivatives, and satisfies $\varphi'(0) = -1$. The precise assumption is given in (6.5) below. For the initial error velocity components $u_e$ and $v_e$ appearing in (2.13)–(2.14), we may consider

$$u_{e0}(x,y) = -U_0^{E}(x)\psi'(y) \qquad \text{and} \qquad v_{e0}(x,y) = \partial_x U_0^{E}(x)\psi(y)\,, \tag{2.16}$$

where $\psi$ is a uniformly analytic function with $\psi(0) = \varphi(0)$ and $\psi'(0) = 0$, with a sufficient decay as $y \to \infty$. The precise assumption is given in (7.2) below. From (2.16) it follows that the error vorticity $\omega_e = -\partial_y u_e + \partial_x v_e$ at the initial time equals

$$\omega_{e0}(x,y) = \partial_x^2 U_0^{E}(x)\psi(y) + U_0^{E}(x)\psi''(y)\,, \tag{2.17}$$

which is shown in (7.1) to be $\mathcal{O}(1)$. Using (2.13)–(2.16), the properties of $\varphi$ and $\psi$ stated above, and the fact that the Euler data are incompressible and satisfy $v_0^{E} = 0$, it follows that the Navier-Stokes datum is incompressible, and satisfies the correct boundary conditions, namely $u_0^{NS} = v_0^{NS} = 0$ on $\{y = 0\}$.

## 3. Main results

Our main result provides an $\mathcal{O}(\epsilon)$ estimate on the error for the vorticity in the Euler+Prandtl expansion (2.12).

THEOREM 3.1. *Assume that the Navier-Stokes initial datum $\boldsymbol{u}_0^{NS}$ and the Euler initial datum $\boldsymbol{u}_0^{E}$ are* compatible, *as described in* (2.13)–(2.14)*, with the Euler datum that satisfies* (5.1)*, the Prandtl initial vorticity $\Omega_0^{P}$ which satisfies* (6.1)*, and with the initial error vorticity $\omega_{e0}$ that satisfies* (7.1)*. Then, there exists $T_* > 0$, independent of $\epsilon$, such that*

$$\sup_{t\in[0,T_*]} \left|\!\left|\!\left| (\omega^{NS} - \omega^{E} - \omega^{P})(\cdot,t) \right|\!\right|\!\right|_t \le C\epsilon\,, \tag{3.1}$$

*where $C > 0$ is a constant. The norm $|||\cdot|||_t$ is defined in* (4.10)*; it represents a norm which encodes $L^1$-based analyticity near the boundary, and Sobolev regularity away from the boundary.*

REMARK 3.1. An example of compatible initial conditions which satisfies the assumptions of Theorem 3.1 is given by the Prandtl and error of the form (2.15) and (2.16), with functions $\varphi$ and $\psi$ which satisfy certain regularity assumptions, cf. (6.5) and (7.2) respectively.

As a direct consequence of Theorem 3.1, we obtain that at the level of the *velocity*, the Euler+Prandtl approximation of the Navier-Stokes solution is $\mathcal{O}(\epsilon)$ in the uniform norm, with respect to both the tangential and the normal variables. Moreover, at any fixed distance away from the boundary, the same convergence rate holds as $\epsilon \to 0$, even without an additional help of the Prandtl corrector.

COROLLARY 3.2. *Under the assumptions of Theorem 3.1, we have*

$$\sup_{t\in[0,T_*]} \|(u^{NS} - u^{E} - \widetilde{u}^{P}, v^{NS} - v^{E} - \epsilon\bar{v}^{P})(\cdot,t)\|_{L^\infty(\mathbb{H})} \le C\epsilon\,. \tag{3.2}$$

*Also, for any set $K \subset \mathbb{H}$ such that* $\operatorname{dist}(K,\partial\mathbb{H}) > 0$*, we have*

$$\sup_{t\in[0,T_0]} \|(\boldsymbol{u}^{NS} - \boldsymbol{u}^{E})(\cdot,t)\|_{L^\infty(K)} \le C\epsilon\,. \tag{3.3}$$

The proofs of Theorem 3.1 and Corollary 3.2 are given in Sections 9 and 10 respectively. The main idea in the proof of Theorem 3.1 is to estimate the error term in the vorticity equation for Navier-Stokes − Euler − Prandtl, cf. (3.17)–(3.18) below. The remainder of this section is dedicated to deriving this error equation, while in the rest of the paper we perform estimates on it.

[2]Compare with the initial datum compatibility assumption in [**56**, Assumption (2.26)]; the fact that we do not need to include higher order correctors in the expansion (2.12), means that we require a less restrictive set of initial conditions for $u_e$ and $v_e$.

**3.1. The evolution for the error velocity and vorticity.** At the velocity level, the Euler+Prandtl expansion of the Navier-Stokes solution is given by

$$u^{\rm NS} = u^{\rm E} + \widetilde{u}^{\rm P} + \epsilon u_{\rm e} \tag{3.4}$$

$$v^{\rm NS} = v^{\rm E} + \epsilon \bar{v}^{\rm P} + \epsilon v_{\rm e}\,, \tag{3.5}$$

where $\widetilde{u}^{\rm P}$ and $\bar{v}^{\rm P}$ are introduced in (2.10) and where $(u_{\rm e}, v_{\rm e})$ stands for the error velocity. At the initial time $t = 0$, the expressions (3.4)–(3.5) correspond to the definition of compatible initial datum, cf. (2.13)–(2.14). The vorticity for the error $(u_{\rm e}, v_{\rm e})$ is denoted by

$$\omega_{\rm e} = -\partial_y u_{\rm e} + \partial_x v_{\rm e} \tag{3.6}$$

and corresponds to the expansion (2.12).

It is also convenient to introduce the approximate velocity

$$u_{\rm a} = u^{\rm E} + \widetilde{u}^{\rm P} \qquad \text{and} \qquad v_{\rm a} = v^{\rm E} + \epsilon \bar{v}^{\rm P} \tag{3.7}$$

and the approximate vorticity

$$\omega_{\rm a} = -\partial_y u_{\rm a} + \partial_x v_{\rm a} = \omega^{\rm E} - \frac{1}{\epsilon}\Omega^{\rm P} + \epsilon \partial_x \bar{v}^{\rm P}\,. \tag{3.8}$$

The evolution equation for $(u_{\rm e}, v_{\rm e})$ is given by (see [**56**, equations (2.32)–(2.39)])

$$(\partial_t - \epsilon^2\Delta)u_{\rm e} + (u_{\rm e}\partial_x + v_{\rm e}\partial_y)u_{\rm a} + (u_{\rm a}\partial_x + v_{\rm a}\partial_y)u_{\rm e} + \epsilon(u_{\rm e}\partial_x + v_{\rm e}\partial_y)u_{\rm e} + \partial_x p_e = f_1 + \frac{1}{\epsilon} g\partial_y \widetilde{u}^{\rm P} \tag{3.9}$$

$$(\partial_t - \epsilon^2\Delta)v_{\rm e} + (u_{\rm e}\partial_x + v_{\rm e}\partial_y)v_{\rm a} + (u_{\rm a}\partial_x + v_{\rm a}\partial_y)v_{\rm e} + \epsilon(u_{\rm e}\partial_x + v_{\rm e}\partial_y)v_{\rm e} + \partial_y p_e = f_2 \tag{3.10}$$

$$\partial_x u_{\rm e} + \partial_y v_{\rm e} = 0 \tag{3.11}$$

$$u_{\rm e}|_{y=0} = 0 \tag{3.12}$$

$$v_{\rm e}|_{y=0} = g\,, \tag{3.13}$$

where $\Delta = \partial_{xx} + \partial_{yy}$. The function $g$ in (3.9) and (3.13) is defined by

$$g = g(t,x) = -\int_0^\infty \partial_x \widetilde{u}^{\rm P} dY = -\bar{v}^{\rm P}|_{Y=0}\,, \tag{3.14}$$

and at the initial time, we have

$$g(t,x)|_{t=0} = -\bar{v}^{\rm P}(x,0,t)|_{t=0}\,.$$

The forcing terms in (3.9)–(3.10) read

$$\begin{aligned} f_1 &= -\frac{1}{\epsilon}\left(\widetilde{u}^{\rm P}\partial_x(u^{\rm E} - U^{\rm E}) + \partial_x\widetilde{u}^{\rm P}(u^{\rm E} - U^{\rm E}) + \partial_y \widetilde{u}^{\rm P}(v^{\rm E} + y\partial_x U^{\rm E})\right) - \bar{v}^{\rm P}\partial_y u^{\rm E} + \epsilon\Delta u^{\rm E} + \epsilon\partial_x^2 \widetilde{u}^{\rm P} \\ &= -Y\left(\widetilde{u}^{\rm P}\frac{\partial_x(u^{\rm E} - U^{\rm E})}{y} + \partial_x \widetilde{u}^{\rm P}\frac{u^{\rm E} - U^{\rm E}}{y} + Y\Omega^{\rm P}\frac{v^{\rm E} + y\partial_x U^{\rm E}}{y^2}\right) - \bar{v}^{\rm P}\partial_y u^{\rm E} + \epsilon\Delta u^{\rm E} + \epsilon\partial_x^2\widetilde{u}^{\rm P} \end{aligned} \tag{3.15}$$

and

$$\begin{aligned} f_2 &= -\left(\partial_t \bar{v}^{\rm P} + u_{\rm a}\partial_x \bar{v}^{\rm P} + v_{\rm a}\partial_y \bar{v}^{\rm P} + \bar{v}^{\rm P}\partial_y v^{\rm E}\right) - \frac{1}{\epsilon}\widetilde{u}^{\rm P}\partial_x v^{\rm E} + \epsilon\Delta v_{\rm a} \\ &= -\left(\partial_t \bar{v}^{\rm P} + u_{\rm a}\partial_x \bar{v}^{\rm P} + Y\frac{v_{\rm a}}{y}\partial_Y \bar{v}^{\rm P} + \bar{v}^{\rm P}\partial_y v^{\rm E}\right) - Y\widetilde{u}^{\rm P}\frac{\partial_x v^{\rm E}}{y} + \epsilon\Delta v_{\rm a}\,. \end{aligned} \tag{3.16}$$

From (3.9)–(3.13), we obtain that $\omega_{\rm e}$ obeys the boundary value problem

$$(\partial_t - \epsilon^2\Delta)\omega_{\rm e} = F \qquad \text{in } \mathbb{H} \tag{3.17}$$

$$\epsilon^2(\partial_y + |\partial_x|)\omega_{\rm e} = \left(\partial_y(-\Delta_D)^{-1}F\right)|_{y=0} + |\partial_x|\int_0^\infty \partial_t \widetilde{u}^{\rm P} dY \qquad \text{on } \partial\mathbb{H}\,, \tag{3.18}$$

where

$$F = -(u_{\rm e}\partial_x + v_{\rm e}\partial_y)\omega_{\rm a} - \frac{1}{\epsilon^2} g\partial_Y \Omega^{\rm P} - (u_{\rm a}\partial_x + v_{\rm a}\partial_y)\omega_{\rm e} - \epsilon(u_{\rm e}\partial_x + v_{\rm e}\partial_y)\omega_{\rm e} + (\partial_x f_2 - \partial_y f_1)\,. \tag{3.19}$$

The boundary condition (3.18) may be derived proceeding similarly to [47], by combining (3.12) and (3.19). Observe that the second boundary term in (3.18) may be written as $|\partial_x| \int_0^\infty \partial_t \widetilde{u}^{\rm P} dY = \frac{\partial_x}{|\partial_x|}\partial_t g$. Recall that the evolution equation for $\widetilde{u}^{\rm P}$ reads

$$(\partial_t - \partial_{YY})\widetilde{u}^{\rm P} + \widetilde{u}^{\rm P}\partial_x U^{\rm E} + U^{\rm E}\partial_x \widetilde{u}^{\rm P} + \widetilde{u}^{\rm P}\partial_x \widetilde{u}^{\rm P} + (\bar{v}^{\rm P} - YU^{\rm E})\partial_Y \widetilde{u}^{\rm P} = 0\,, \tag{3.20}$$

where

$$\bar{v}^{\rm P}(Y) = \int_Y^\infty \partial_x \widetilde{u}^{\rm P}\, dY' \tag{3.21}$$

(cf. [56, equation (2.20)]). Lastly, observe that using the definition $g = -\partial_x(\int_0^\infty \widetilde{u}^{\rm P} dY)$, we rewrite the integral in the last term on the right side of (3.18) as

$$\begin{aligned}\int_0^\infty \partial_t \widetilde{u}^{\rm P} dY &= \int_0^\infty \big(\partial_{YY}\widetilde{u}^{\rm P} - \partial_x(\widetilde{u}^{\rm P}U^{\rm E}) - \widetilde{u}^{\rm P}\partial_x\widetilde{u}^{\rm P} - v^{\rm P}\partial_y \widetilde{u}^{\rm P}\big)\, dY\\ &= \int_0^\infty \big(\partial_Y \Omega^{\rm P} - \partial_x(\widetilde{u}^{\rm P}U^{\rm E}) - \widetilde{u}^{\rm P}\partial_x\widetilde{u}^{\rm P} - \partial_x(\widetilde{u}^{\rm P} + U^{\rm E})\widetilde{u}^{\rm P}\big)\, dY\\ &= -\Omega^{\rm P}|_{Y=0} + U^{\rm E} g - 2\partial_x U^{\rm E}\int_0^\infty \widetilde{u}^{\rm P} dY - \partial_x \int_0^\infty (\widetilde{u}^{\rm P})^2 dY\,,\end{aligned} \tag{3.22}$$

where we used (3.20) in the first equality and thus the boundary condition in (3.18) reads

$$\begin{aligned}\epsilon^2(\partial_y + |\partial_x|)\omega_{\rm e} &= \big(\partial_y(-\Delta_D)^{-1}F\big)|_{y=0} - |\partial_x|\,\Omega^{\rm P}|_{Y=0} + |\partial_x|\,U^{\rm E} g\\ &\quad - 2\,|\partial_x|\,\partial_x U^{\rm E}\int_0^\infty \widetilde{u}^{\rm P} dY - |\partial_x|\,\partial_x \int_0^\infty (\widetilde{u}^{\rm P})^2 dY \text{ on } \partial\mathbb{H}\,.\end{aligned} \tag{3.23}$$

Since the error vorticity equation (3.17) has a forcing term which depends on the Euler and Prandtl solutions, it is natural that we first perform suitable analytic and Sobolev estimates for these Euler (cf. Section 5) and Prandtl (cf. Section 6) solutions, with the initial conditions given by (2.13)–(2.14). Prior to this, in the following section we introduce the functional framework in which these estimates are performed.

## 4. The functional framework

**4.1. The base analytic norms.** For $\mu \in (0,1]$ we define the complex domains

$$\Omega_\mu = \{z \in \mathbb{C} : 0 \le \mathbb{R}{\rm e}\, z \le 1, |\mathbb{I}{\rm m}\, z| \le \mu \mathbb{R}{\rm e}\, z\} \cup \{z \in \mathbb{C} : 1 \le \mathbb{R}{\rm e}\, z \le 1+\mu, |\mathbb{I}{\rm m}\, z| \le 1 + \mu - \mathbb{R}{\rm e}\, z\} \tag{4.1}$$

and

$$\widetilde{\Omega}_\mu = \big\{Z \in \mathbb{C} : 0 \le \mathbb{R}{\rm e}\, Z, |\mathbb{I}{\rm m}\, Z| \le \mu \mathbb{R}{\rm e}\, Z\big\}\,. \tag{4.2}$$

We note that the domain $\widetilde{\Omega}_\mu$ is much larger than the domain $\Omega_\mu$, and allows $\mathbb{R}{\rm e}\, Z$ to be arbitrarily large, while the domain $\Omega_\mu$ is located near the boundary $0 \le \mathbb{R}{\rm e}\, y \le 1 + \mu$. We use $f_\xi(y) \in \mathbb{C}$ to denote the Fourier transform of $f(x,y)$ with respect to the $x$ variable at frequency $\xi \in \mathbb{Z}$, i.e., $f(x,y) = \sum_{\xi\in\mathbb{Z}} f_\xi(y) e^{ix\xi}$.

We define three types of analytic norms, $Y_{\lambda,\mu}$, $Y_{\lambda,\mu,\infty}$, and $P_{\lambda,\mu,\infty}$. The principal purpose of the $Y_{\lambda,\mu}$ norm is to control the remainder of the Prandtl expansion, the main role of the $Y_{\lambda,\mu,\infty}$ norm is to estimate the Euler solution in analytic spaces, while the $P_{\lambda,\mu,\infty}$ norm bounds the Prandtl solution in the domain $\widetilde{\Omega}_\mu$. Let $\lambda, \mu \in (0,1]$.

- For a complex function $f(y)$ defined on $\Omega_\mu$ , let

$$\|f\|_{\mathcal{L}^1_\mu} = \sup_{0 \le \theta < \mu} \|f\|_{L^1(\partial\Omega_\theta)}\,, \tag{4.3}$$

  and for a complex function $f(x,y)$ defined on the domain $\mathbb{T}\times\Omega_\mu$, we introduce the $L^1_y$-based analytic norm

$$\|f\|_{Y_{\lambda,\mu}} = \sum_{\xi\in\mathbb{Z}} \|e^{\lambda(1+\mu-y)|\xi|} f_\xi\|_{\mathcal{L}^1_\mu}\,. \tag{4.4}$$

- For a complex valued function $f(x,y)$ defined on $\mathbb{T}\times\Omega_\mu$, we define the $L^\infty_y$-based analytic norm

$$\|f\|_{Y_{\lambda,\mu,\infty}} = \sum_{\xi\in\mathbb{Z}} \|e^{\lambda(1+\mu-y)|\xi|} f_\xi\|_{L^\infty(\Omega_\mu)}\,. \tag{4.5}$$

  If $f = f(x)$ is independent of $y$ and only depends on $x \in \mathbb{T}$, we replace the norm $\|f_\xi\|_{L^\infty(\Omega_\mu)}$ simply by $|f_\xi|$, and still use the notation in (4.5).

- For a function $f(x,Y)$ defined on the domain $\mathbb{T}\times\widetilde{\Omega}_\mu$, we define the $L^\infty_Y$-based analytic norm
$$\|f\|_{P_{\lambda,\mu,\infty}} = \sum_{\xi\in\mathbb{Z}} e^{\lambda(1+\mu)|\xi|}\|f_\xi\|_{L^\infty(\widetilde{\Omega}_\mu)}. \tag{4.6}$$
If $f=f(x)$ is independent of $Y$ and only depends on $x\in\mathbb{T}$, e.g. trace terms at $Y=0$ or terms which are integrated in $Y$, we replace the norm $\|f_\xi\|_{L^\infty(\widetilde{\Omega}_\mu)}$ simply by $|f_\xi|$, and still use the notation in (4.6).

Note that both the $Y_{\lambda,\mu}$ and $Y_{\lambda,\mu,\infty}$ norms only require the corresponding function to be analytic in $x$ near the boundary $\{y=0\}$, whereas the $P_{\lambda,\mu,\infty}$ norm requires also analyticity at $Y$-large. Moreover, unlike in [**38**], the $Y_{\lambda,\mu,\infty}$ norm is not weighted in the $y$ variable.

**4.2. The Sobolev norms.** To control the Sobolev part of a function $f$ away from the boundary, for $\mu>0$ we introduce
$$\|f\|_{S_\mu} = \sum_{\xi}\|yf_\xi\|_{L^2(y\geq 1+\mu)}\,. \tag{4.7}$$
Note that the $S_\mu$ norm is $\ell^1_\xi$, so that in view of the Hausdorff-Young inequality, we have $\|yf\|_{L^\infty_x L^2_y(y\geq 1+\mu)}\leq\|f\|_{S_\mu}$.

Using (4.4) and (4.7) we also define
$$\|f\|_{Y_{\lambda,\mu}\cap S_\mu} = \|f\|_{Y_{\lambda,\mu}} + \|f\|_{S_\mu}\,.$$
Note that the norm $Y_{\lambda,\mu}\cap S_\mu$ controls the $L^1$ norm in the analytic region $0\leq\mathbb{R}\mathrm{e}\,y\leq 1+\mu$, and a weighted $L^2$ norm in the Sobolev region $y\geq 1+\mu$.

As a genuine $L^2_{x,y}$-based Sobolev norm we choose
$$\|f\|_S^2 = \|yf\|^2_{L^2(y\geq 1/2)} = \sum_{\xi\in\mathbb{Z}}\|yf_\xi\|^2_{L^2(y\geq 1/2)}$$
and denote the higher derivative version by
$$\|f\|_Z = \sum_{0\leq i+j\leq 3}\|\partial_x^i\partial_y^j f\|_S = \sum_{0\leq i+j\leq 3}\|y\partial_x^i\partial_y^j f\|_{L^2(y\geq 1/2)}\,. \tag{4.8}$$
Note that since $(1+|\xi|)^{-1}\in\ell^2_\xi$, we have the lossy estimates $\|f\|_{S_\mu}\leq\|f\|_S+\|\partial_x f\|_S\leq\|f\|_{S_\mu}+\|\partial_x f\|_{S_\mu}$.

**4.3. The cumulative error norm.** Finally, we define the norm $|\!|\!|\cdot|\!|\!|_t$ which appears in Theorem 3.1.

Before doing so, we *fix* two sufficiently small parameters $\lambda_*,\mu_*\in(0,1]$, which are independent of $\epsilon$, and only depend on the parameter $\lambda_0$ which appears in the assumptions on the Euler datum (cf. (5.1)) and the Prandtl datum (cf. (6.1)), and the parameters $\mu_2,\lambda_2$ which appear in the assumption on the initial error vorticity (cf. (7.1)). The precise values of $\lambda_*,\mu_*$ are given in (7.3) below; at this point we only emphasize that these parameters are determined in terms of the datum, and that they are independent of $\epsilon$. Lastly, let $\gamma_*\geq 2$ be a sufficiently large parameter representing the rate of decay of the analyticity radius. This parameter is also independent of $\epsilon$, and its value shall be determined in the proof (see the line above (9.2)). Throughout the paper, the time parameter is chosen to satisfy $0\leq t\leq\min\{1,\mu_*/(2\gamma_*)\}$, so that $t\leq 1$ and $\mu_*-\gamma_* t\geq\mu_*/2>0$; in fact, we let $t\in[0,T_*]$, where $T_*\in(0,1]$ is independent of $\epsilon$, is given explicitly in (7.3).

To treat the loss of a derivative in the nonlinear terms, in terms of the parameters $\mu_*$ and $\gamma_*$ discussed above, we use (4.4) to define the cumulative $L^1_y$-based analytic norm
$$\|f(t)\|_{Y(t)} = \sup_{0<\mu<\mu_*-\gamma_* t}\left(\sum_{i+j\leq 1}\|\partial_x^i(y\partial_y)^j f\|_{Y_{\lambda_*,\mu}} + (\mu_*-\mu-\gamma_* t)^{1/3}\sum_{i+j=2}\|\partial_x^i(y\partial_y)^j f\|_{Y_{\lambda_*,\mu}}\right), \tag{4.9}$$
for all $0\leq t\leq T_*$. Lastly, for the same range of $t$, using (4.8) we denote by
$$|\!|\!|\omega|\!|\!|_t = \|\omega(t)\|_{Y(t)} + \|\omega(t)\|_Z \tag{4.10}$$
the cumulative error vorticity norm.

REMARK 4.1 (**Implicit constants**). We emphasize that throughout the paper the implicit constants in the symbols $\lesssim$ are *never* allowed to depend on the large parameters $\epsilon^{-1}$, $\gamma_*$, and $T_*^{-1}$. These implicit constants are however allowed to depend on parameters independent of $\epsilon$ and $\gamma_*$, such as $\lambda_0,\lambda_1,\lambda_2,\lambda_*,\mu_0,\mu_1,\mu_2,\mu_*$, or $\kappa$.

**4.4. Functional inequalities.** We recall several useful properties of the norms introduced in (4.4)–(4.6). First, from the Cauchy integral formula, we deduce the following inequality (cf. also [**53**, Lemma 2.2]).

LEMMA 4.1 (**Analytic recovery**). *For $0 \leq \mu < \widetilde{\mu} < \mu_* - \gamma_* s$, we have*

$$\sum_{i+j=1} \|\partial_x^i (y\partial_y)^j f\|_{Y_{\lambda,\mu}} \lesssim \frac{1}{\widetilde{\mu}-\mu}\|f\|_{Y_{\lambda,\widetilde{\mu}}} ,$$

*where the implicit constant is universal.*

We omit the proof of Lemma 4.1 and refer the reader to [**53**]. In the next lemma, we record a number of useful product estimates concerning the analytic norms. Similar bounds to the ones stated in (4.11) below were previously established in [**53**] and [**38**].

LEMMA 4.2 (**Product estimates**). *For $\lambda, \mu \in (0,1]$, we have the inequalities*

$$\|f(x,Y)g(x,y)\|_{Y_{\lambda,\mu}} \lesssim \|f(x,Y)\|_{P_{\lambda,\mu,\infty}}\|g(x,y)\|_{Y_{\lambda,\mu}} \tag{4.11a}$$

$$\|f(x,Y)g(x,y)\|_{Y_{\lambda,\mu}} \lesssim \epsilon\|(1+Y)^{3/2}f(x,Y)\|_{P_{\lambda,\mu,\infty}}\|g(x,y)\|_{Y_{\lambda,\mu,\infty}} \tag{4.11b}$$

$$\|f(x,y)g(x,y)\|_{Y_{\lambda,\mu}} \lesssim \|f(x,y)\|_{Y_{\lambda,\mu,\infty}}\|g(x,y)\|_{Y_{\lambda,\mu}} \tag{4.11c}$$

$$\|f(x,Y)g(x,y)\|_{S_\mu} \lesssim \epsilon^\theta\|Y^\theta f(x,Y)\|_{P_{\lambda,\mu,\infty}}\|g(x,y)\|_{S_\mu} \tag{4.11d}$$

$$\|f(x,Y)g(x,y)\|_{S_\mu} \lesssim \epsilon^\theta\|Y^\theta f(x,Y)\|_{P_{\lambda,\mu,\infty}}\left(\|g(x,y)\|_{L^2(y\geq 1/2)} + \|\partial_x g(x,y)\|_{L^2(y\geq 1/2)}\right) \tag{4.11e}$$

$$\|f(x,Y)g(x,y)\|_{S_\mu} \lesssim \epsilon^\theta\|Y^\theta f(x,Y)\|_{P_{\lambda,\mu,\infty}}\left(\|g(x,y)\|_{L^2_xL^\infty_y(y\geq 1/2)} + \|\partial_x g(x,y)\|_{L^2_xL^\infty_y(y\geq 1/2)}\right) \tag{4.11f}$$

$$\|f(x,y)g(x,y)\|_{S_\mu} \lesssim \|f\|_{S_\mu}\left(\|g\|_{L^2_xL^\infty_y(y\geq 1+\mu)} + \|\partial_x g\|_{L^2_xL^\infty_y(y\geq 1+\mu)}\right) , \tag{4.11g}$$

*for any $\theta \geq 2$, whenever the right sides of the above inequalities are finite. For simplicity of notation, we write $Y$ instead of $\mathbb{R}\mathrm{e}\, Y$ for the weights on the right sides.*

PROOF OF LEMMA 4.2. We first observe that for an analytic function $f(x,Y)$ defined on $\mathbb{T}\times\widetilde{\Omega}_\mu$, with $Y = y/\epsilon$, the function $(x,y) \mapsto f(x,y/\epsilon)$ is analytic in $\Omega_\mu$, since $y \in \Omega_\mu$ implies $Y \in \widetilde{\Omega}_\mu$. This observation is used throughout the proof.

Since the $Y_{\lambda,\mu}$ norm contains an $L^1$ norm with respect to the $y$ variable along the polygonal path $\partial\Omega_\theta$ with $\theta < \mu$, and since we have $dy = \epsilon dY$ and $(1+\mathbb{R}\mathrm{e}\, Y)^{-\frac{3}{2}} \in L^1_Y$, we have a useful bound

$$\begin{aligned}\|f(x,Y)\|_{Y_{\lambda,\mu}} &\leq \|(1+\mathbb{R}\mathrm{e}\, Y)^{-3/2}\|_{\mathcal{L}^1_\mu}\|(1+\mathbb{R}\mathrm{e}\, Y)^{3/2}f(x,Y)\|_{P_{\lambda,\mu,\infty}} \\ &\lesssim \epsilon\|(1+Y)^{3/2}f(x,Y)\|_{P_{\lambda,\mu,\infty}} ,\end{aligned} \tag{4.12}$$

where the implicit constant is universal, and we omitted the real part of the weight appearing on the right side. Next, we note that by the definition of the domain $\Omega_\mu$, we have

$$\|f(x,y)\|_{Y_{\lambda,\mu}} \lesssim \|f(x,y)\|_{Y_{\lambda,\mu,\infty}} , \tag{4.13}$$

where the implicit constant is universal. The above two estimates bound the $L^1$-based analytic norm, in terms of those based on $L^\infty$.

Next, we consider product estimates, and start with (4.11a). Again, using that $y \in \Omega_\mu$ implies $Y = y/\epsilon \in \widetilde{\Omega}_\mu$, from the Hölder inequality we obtain

$$\begin{aligned}\|fg\|_{Y_{\lambda,\mu}} &= \sum_\xi \|e^{\lambda(1+\mu-y)_+|\xi|}\sum_{\xi'} f_{\xi'}(Y)g_{\xi-\xi'}(y)\|_{\mathcal{L}^1_\mu} \\ &\leq \sum_\xi\sum_{\xi'}\|e^{\lambda(1+\mu-y)_+|\xi-\xi'|}g_{\xi-\xi'}(y)\|_{\mathcal{L}^1_\mu}\sup_{Y\in\widetilde{\Omega}_\mu}|f_{\xi'}(Y)e^{\lambda(1+\mu)|\xi'|}| \\ &\leq \|g\|_{Y_{\lambda,\mu}}\|f\|_{P_{\lambda,\mu,\infty}}.\end{aligned}$$

Similarly for (4.11b), we appeal to the above argument and to the proof of (4.12), to obtain

$$\begin{aligned}\|fg\|_{Y_{\lambda,\mu}} &\leq \sum_\xi\sum_{\xi'}\|e^{\lambda(1+\mu)|\xi'|}f_{\xi'}(Y)\|_{\mathcal{L}^1_\mu}\sup_{y\in\Omega_\mu}|e^{\lambda(1+\mu-y)_+|\xi-\xi'|}g_{\xi-\xi'}(y)| \\ &\leq \epsilon\|(1+Y)^{3/2}f\|_{P_{\lambda,\mu,\infty}}\|g\|_{Y_{\lambda,\mu,\infty}} .\end{aligned} \tag{4.14}$$

The inequality (4.11c) is a consequence of the Hölder inequality in $y$ on the domain $\Omega_\mu$.

In order to prove the bound (4.11d), we note that by the definition of the $S_\mu$ norm in (4.7), Hölder's inequality in $y$, and the fact that $y \geq 1+\mu$ implies that $Y = y/\epsilon \geq (1+\mu)/\epsilon \geq 1/\epsilon$, we deduce that

$$\begin{aligned}\|fg\|_{S_\mu} &\leq \sum_\xi \sum_{\xi'} \|f_{\xi'}\|_{L^\infty(y\geq 1+\mu)} \|y g_{\xi-\xi'}\|_{L^2(y\geq 1+\mu)} \\ &\leq \|g\|_{S_\mu} \sum_\xi \epsilon^\theta \left\|Y^\theta f_\xi\right\|_{L^\infty(Y\geq 1/\epsilon)} \\ &\leq \|g\|_{S_\mu} \epsilon^\theta \left\|Y^\theta f_\xi\right\|_{P_{\lambda,\mu,\infty}} ,\end{aligned}$$

for any $\lambda, \theta \geq 0$. In a similar fashion we may establish (4.11e) as

$$\begin{aligned}\|fg\|_{S_\mu} &\leq \sum_\xi \sum_{\xi'} \epsilon \|Y f_{\xi'}\|_{L^\infty(y\geq 1+\mu)} \|g_{\xi-\xi'}\|_{L^2(y\geq 1+\mu)} \\ &\leq \sum_\xi \epsilon^\theta \left\|Y^\theta f_\xi\right\|_{L^\infty(Y\geq 1/\epsilon)} \sum_{\xi''} \|g_{\xi''}\|_{L^2(y\geq 1+\mu)} \\ &\lesssim \epsilon^\theta \left\|Y^\theta f_\xi\right\|_{P_{\lambda,\mu,\infty}} \left(\|g\|_{L^2(y\geq 1/2)} + \|\partial_x g\|_{L^2(y\geq 1/2)}\right) ,\end{aligned}$$

where in the last inequality we have used Plancherel, and the fact $(1+|\xi''|)^{-1} \in \ell^2_{\xi''}$. The proof of (4.11f) is similar as we have

$$\begin{aligned}\|fg\|_{S_\mu} &\leq \sum_\xi \sum_{\xi'} \epsilon^2 \left\|Y^2 f_{\xi'}\right\|_{L^\infty(y\geq 1+\mu)} \left\|y^{-1} g_{\xi-\xi'}\right\|_{L^2(y\geq 1+\mu)} \\ &\leq \sum_\xi \epsilon^\theta \left\|Y^\theta f_\xi\right\|_{L^\infty(Y\geq 1/\epsilon)} \sum_{\xi''} \|g_{\xi''}\|_{L^\infty(y\geq 1+\mu)} \\ &\lesssim \epsilon^\theta \left\|Y^\theta f_\xi\right\|_{P_{\lambda,\mu,\infty}} \left(\|g\|_{L^2_x L^\infty_y(y\geq 1/2)} + \|\partial_x g\|_{L^2_x L^\infty_y(y\geq 1/2)}\right)\end{aligned}$$

since $\|y^{-1}\|_{L^2(y\geq 1+\mu)} \lesssim 1$.

The last inequality, (4.11g) follows directly from the definition (4.7) and Hölder's inequality

$$\|fg\|_{S_\mu} \leq \sum_\xi \sum_{\xi'} \|y f_{\xi'}\|_{L^2(y\geq 1+\mu)} \|g_{\xi-\xi'}\|_{L^\infty(y\geq 1+\mu)} \quad \leq \|f\|_{S_\mu} \left(\|g\|_{L^2_x L^\infty_y(y\geq 1+\mu)} + \|\partial_x g\|_{L^2_x L^\infty_y(y\geq 1+\mu)}\right) ,$$

which concludes the proof. □

Next, we recall the following elliptic estimates for the velocity; for a proof, we refer the reader to [**38**, Lemma 6.3] and [**39**, Lemmas 4.2 and 5.1].

LEMMA 4.3 (**Elliptic estimates**). *Let $(u,v)$ be the velocity obtained from the vorticity $\omega$ via the Biot-Savart law, cf.* (7.35)–(7.36) *with $g=0$. For $\mu \in (0, \mu_* - \gamma_* t)$ and $\lambda \in (0, \lambda_*]$, we have the estimates*

$$\|\partial_x^i (y\partial_y)^j u\|_{Y_{\lambda,\mu,\infty}} \lesssim \|\partial_x^{i+j}\omega\|_{Y_{\lambda,\mu}\cap S_\mu} + j\left(\|\omega\|_{Y_{\lambda,\mu}} + \|y\partial_y \omega\|_{Y_{\lambda,\mu}}\right)$$

*and*

$$\left\|\partial_x^i (y\partial_y)^j \left(\frac{v}{y}\right)\right\|_{Y_{\lambda,\mu,\infty}} \lesssim \|\partial_x^{i+1}\omega\|_{Y_{\lambda,\mu}\cap S_\mu} ,$$

*for all non-negative integers $i,j$ such $i+j \leq 1$. For the Sobolev norm away from the boundary, one has*

$$\sum_{i+j=3} \left(\|\partial_x^i \partial_y^j u\|_{L^2_{x,y}(y\geq 1/4)} + \|\partial_x^i \partial_y^j v\|_{L^2_{x,y}(y\geq 1/4)}\right) \lesssim \|\!|\omega|\!\|_t$$

*and*

$$\sum_{i+j\leq 2} \left(\|\partial_x^i \partial_y^j u\|_{L^\infty_{x,y}(y\geq 1/4)} + \|\partial_x^i \partial_y^j v\|_{L^\infty_{x,y}(y\geq 1/4)}\right) \lesssim \|\!|\omega|\!\|_t , \tag{4.15}$$

*for all $t \in [0, T_*]$.*

## 5. Uniform analyticity of the Euler solution in a strip

In this section, we estimate the solution of the Euler equations (1.5)–(1.6) posed on the half-space $\mathbb{H} = \mathbb{T} \times \mathbb{R}_+$ with the boundary condition (2.2) and the initial condition (2.1). We require the initial data to be uniformly analytic in $x$ and $y$ near the boundary. Away from the boundary, i.e., for $y \geq 2$, we only require Sobolev regularity. These assumptions are stated in terms of the initial vorticity $\omega_0^{\rm E}$. Namely, we assume that $\omega_{0,\xi}^{\rm E}(y)$ is analytic in the domain $\{y \in \mathbb{C} : 0 < \mathbb{R}{\rm e}\, y < 2,\, |\mathbb{I}{\rm m}\, y| \leq 2\}$ with values in the $L^1_\xi$ space with the weight $e^{\lambda_0|\xi|}$, is continuous on the closure, and satisfies

$$\sum_{\xi\in\mathbb{Z}} e^{\lambda_0|\xi|} \sup_{0\leq \mathbb{R}{\rm e}\, y\leq 2, |\mathbb{I}{\rm m}\, y|\leq 2} |\omega_{0,\xi}^{\rm E}(y)| + \sum_{i+j\leq 4} \|y\partial_x^i (y\partial_y)^j \omega_0^{\rm E}\|_{L^2(y\geq 1/2)} \lesssim 1\,, \tag{5.1}$$

for some $\lambda_0 \in (0,1]$. We allow all the constants to depend on $\lambda_0$. Note that $\omega_0^{\rm E}$ satisfies the assumptions on the initial data in [**38**, Theorem 3.1]. Our goal in this section is to establish the bounds stated in Lemma 5.7 below. To this end, we first prove that if the initial Euler data satisfies (5.1), then the solution of the Euler equations remains analytic near the boundary, locally in time.

THEOREM 5.1. *Assume that (5.1) holds, and let $\omega^{\rm E}$ be the vorticity corresponding to the unique solution of the Cauchy problem for the Euler equations (1.5)–(1.6), (2.1), with the initial vorticity $\omega_0^{\rm E}$. Then there exists $T_0 \in (0,1]$ such that*

$$\sum_{i+j\leq 4} \|y\partial_x^i\partial_y^j \omega^{\rm E}(t)\|^2_{L^2(y\geq 1/2)} \lesssim 1 \tag{5.2}$$

*and*

$$\sum_{0\leq i+j\leq 4} \|\partial_x^i\partial_y^j \boldsymbol{u}^{\rm E}(t)\|_{L^2_{x,y}(y\geq 1/2)} + \sum_{0\leq i+j\leq 3} \|\partial_x^i\partial_y^j \boldsymbol{u}^{\rm E}(t)\|_{L^\infty_{x,y}(y\geq 1/2)} \lesssim 1\,, \tag{5.3}$$

*for $t \in [0,T_0]$. Moreover, the vorticity $\omega^{\rm E}$ and the velocity $\mathbf{u}^{\rm E}$ are uniformly real-analytic in $(x,y) \in \mathbb{T}\times[0,1]$ in the sense that there exists a constant $\zeta_0 \in (0,1]$ such that*

$$\sum_{i,j} \frac{\zeta_0^{i+j}}{(i+j)!} \|\partial_x^i\partial_y^j \omega^{\rm E}\|_{L^\infty(\mathbb{T}\times[0,1])} \lesssim 1 \tag{5.4}$$

*and*

$$\sum_{i,j} \frac{\zeta_0^{i+j}}{(i+j)!} \|\partial_x^i\partial_y^j \boldsymbol{u}^{\rm E}\|_{L^\infty(\mathbb{T}\times[0,1])} \lesssim 1\,, \tag{5.5}$$

*for $t \in [0,T_0]$.*

The inequalities (5.4) and (5.5) assert the uniform analyticity up to $y = 0$, instead of only analyticity in a wedge.

We divide the proof of Theorem 5.1 into several steps. First, we obtain the interior analyticity of solutions, which is asserted in the next lemma.

LEMMA 5.2. *Assume that $\omega_0^{\rm E}$ satisfies (5.1). Then there exist constants $T_0, \mu_0 \in (0,1]$ and $C \geq 1$ such that we have (5.2), (5.3), and*

$$\sum_\xi e^{|\xi|/C} |\omega_\xi^{\rm E}(t,y)| \lesssim 1, \qquad y \in \Omega_{\mu_0} \cap \{y : 1/2 < \mathbb{R}{\rm e}\, y < 1+\mu_0\}\,, \tag{5.6}$$

*for all $t \in [0,T_0]$.*

PROOF OF LEMMA 5.2. In order to apply [**38**, Theorem 3.1], note that we have

$$\sum_{i+j\leq 2} \left(\|\partial_x^i (y\partial_y)^j \omega_0\|_{Y_{\epsilon_0,1}} + \|\partial_x^i (y\partial_y)^j \omega_0\|_{Y_{\epsilon_0,1,\infty}}\right) + \sum_{i+j\leq 3} \|y\partial_x^i (y\partial_y)^j \omega_0^{\rm E}\|_{L^2(y\geq 1/2)} \lesssim 1\,,$$

where $\epsilon_0 = \lambda_0/2$, i.e., the condition (3.1) in [**38**, Theorem 3.1] is fulfilled. Therefore, for every $\epsilon \in (0,1]$ sufficiently small there exists a unique solution to the Navier-Stokes equations $\omega^{{\rm NS},\epsilon}$, with the initial data $\omega_0^{\rm E}$ on a uniform in $\epsilon$ time interval $[0,T_0]$, and on this interval the solutions $\omega^{{\rm NS},\epsilon}$ are uniformly bounded and analytic in $\Omega_{\mu_0}$, i.e.,

$$\|\max\{\epsilon, \mathbb{R}{\rm e}\, y\}\omega^{{\rm NS},\epsilon}\|_{Y_{\epsilon_0,\mu_0,\infty}} \lesssim 1\,, \tag{5.7}$$

for some $\mu_0 \in (0,1]$ which is independent of $\epsilon$; additionally, by [38],

$$\omega^{\mathrm{NS},\epsilon} \to \omega^{\mathrm{E}} \text{ in } C([0,T_0], L^2(\mathbb{H})) \text{ as } \epsilon \to 0. \tag{5.8}$$

Note that this solution is different than the one in (2.12) since it starts from a different initial data. Using (5.7), we get a uniform in $\epsilon$ bound

$$\sum_\xi e^{\epsilon_0(1+\mu_0-\mathbb{R}\mathrm{e}\, y)|\xi|} \max\{\epsilon, \mathbb{R}\mathrm{e}\, y\} |\omega_\xi^{\mathrm{NS},\epsilon}(t,y)| \lesssim 1, \qquad y \in \Omega_{\mu_0}, \quad t \in [0,T_0],$$

which implies

$$\sum_\xi e^{\epsilon_0\mu_0|\xi|/2} |\omega_\xi^{\mathrm{NS},\epsilon}(t,y)| \lesssim 1, \qquad \frac{1}{2} < \mathbb{R}\mathrm{e}\, y < 1 + \frac{\mu_0}{2}, \quad y \in \Omega_{\mu_0}, \tag{5.9}$$

for every $t \in [0,T_0]$. We next claim that the Euler solution satisfies

$$\sum_\xi e^{\epsilon_0\mu_0|\xi|/4} |\omega_\xi^{\mathrm{E}}(t,y)| \lesssim 1, \qquad \frac{1}{2} < \mathbb{R}\mathrm{e}\, y < 1 + \frac{\mu_0}{2}, \quad y \in \Omega_{\mu_0}, \tag{5.10}$$

for $t \in [0,T_0]$. To prove (5.10), first observe that we have (5.8). Fix any $t_0 \in [0,T_0]$ and $m_0 \in \mathbb{N}$. Due to the uniform bound (5.9) at time $t_0$, we may apply the vector version of the Montel theorem and deduce that there exists an analytic function $f$ on $\Omega_0 = \{y \in \Omega_{\mu_0} : 1/2 < \mathbb{R}\mathrm{e}\, y < 1 + \mu_0/2\}$ with values in the space of functions $g$ such that

$$\sum_\xi e^{\epsilon_0\mu_0|\xi|/4} |g(t,y)| < \infty, \qquad \frac{1}{2} < \mathbb{R}\mathrm{e}\, y < 1 + \frac{\mu_0}{2}, \quad y \in \Omega_{\mu_0} \tag{5.11}$$

and a sequence $\epsilon_1, \epsilon_2, \ldots \to 0$ such that $\omega^{\mathrm{NS},\epsilon_j}(t_0)$ converges to $f$ uniformly on compact subsets of $\Omega_0$, with values in the space corresponding to (5.11). By the uniform bound

$$\sum_{\xi=-m_0}^{m_0} e^{\epsilon_0\mu_0|\xi|/2} |\omega_\xi^{\mathrm{NS},\epsilon}(t,y)| \lesssim 1, \qquad \frac{1}{2} < \mathbb{R}\mathrm{e}\, y < 1 + \frac{\mu_0}{2}, \quad y \in \Omega_{\mu_0},$$

for every $m_0 \in \mathbb{N}$ (which is a consequence of (5.9) at $t = t_0$) for $\epsilon = \epsilon_1, \epsilon_2, \ldots$, the function $f$ also satisfies the same bound. Finally, note that $f = \{\omega_\xi^{\mathrm{E}}(t_0)\}_{\xi=-m_0}^{m_0}$ by $\omega^{\mathrm{NS},\epsilon} \to \omega^{\mathrm{E}}$ in $C([0,T_0], L^2(\mathbb{H}))$, and we obtain

$$\sum_{\xi=m_0}^{m_0} e^{\epsilon_0\mu_0|\xi|/2} |\omega_\xi^{\mathrm{E}}(t,y)| \lesssim 1, \qquad \frac{1}{2} < \mathbb{R}\mathrm{e}\, y \le 1 + \frac{\mu_0}{2}, \quad y \in \Omega_{\mu_0}$$

at $t = t_0$, and $t_0 \in [0,T_0]$. Since $m_0 \in \mathbb{N}$ is arbitrary, we obtain (5.10) and (5.6) follows by replacing $\mu_0$ with $\mu_0/2$.

Next, we establish (5.2), which is obtained using a weighted Sobolev estimate with a weight $\phi(y) = (y^2+1)^{1/2}$. First, note that

$$\sum_{|\alpha|\le 4} (\|\partial^\alpha \omega\|_{L^2}^2 + \|\partial^\alpha \boldsymbol{u}\|_{L^2}^2) \lesssim 1, \tag{5.12}$$

by the local $H^4$ existence. The weighted energy $\psi = \sum_{|\alpha|\le 4} \int |\partial^\alpha \omega| \phi^2$ satisfies

$$\begin{aligned} \frac{1}{2}\frac{d}{dt}\psi = 2 \sum_{|\alpha|=4} \int |\partial^\alpha \omega|^2 \phi u \cdot \nabla \phi - \sum_{|\alpha|\le 4, 0\le\beta\le\alpha, |\beta|\le 2, (|\alpha|,\beta)\ne(4,0)} \int \partial^\beta \boldsymbol{u} \cdot \nabla \partial^{\alpha-\beta} \omega \partial^\alpha \omega \phi^2 \\ - \sum_{|\alpha|\le 4, 0\le\beta\le\alpha, |\beta|\ge 3} \int \partial^\beta \boldsymbol{u} \cdot \nabla \partial^{\alpha-\beta} \omega \partial^\alpha \omega \phi^2 . \end{aligned}$$

All the terms are estimated in a straight-forward way by (5.12) and using that all the derivatives of $\phi$ are uniformly bounded. For the first term, we estimate the integral by $\|\partial^\alpha \omega \phi\|_{L^2} \|\partial^\alpha \omega\|_{L^2} \|u\|_{L^\infty}$, for the second term, we bound the integral by $\|\partial^\beta \boldsymbol{u}\|_{L^\infty} \|D^{|\alpha|-|\beta|+1} \omega \phi\|_{L^2} \|\partial^\alpha \omega \phi\|_{L^2}$, while the integral in the third term by $\|\partial^\beta \boldsymbol{u}\|_{L^\infty} \|D^{|\alpha|-|\beta|+1} \omega \phi\|_{L^2} \|\partial^\alpha \omega \phi\|_{L^2}$. We omit further details.

Finally, the inequality (5.3) follows by the Biot-Savart law as in the proof of [39, Lemma 5.1]. □

Next, we provide estimates on a solution of the Euler equation in the region $\mathbb{T} \times [0,1]$ in the analytic norm

$$\|\omega\|_{\widetilde{A}_\tau} = \sum_{|\alpha|\ge 3} \frac{\tau^{|\alpha|-3}}{(|\alpha|-3)!} \widetilde{\delta}^{\alpha_1} \bar{\delta}^{\alpha_2} \|\partial^\alpha \omega\|_{L^2(\mathbb{T}\times[0,1])}, \tag{5.13}$$

where $\alpha = (\alpha_1, \alpha_2)$ and $\tau > 0$. In (5.13), the parameters $\widetilde{\delta}, \bar{\delta} \in (0,1]$ are constants such that

$$\widetilde{\delta}, \bar{\delta} \leq \frac{1}{C} \quad \text{and} \quad \bar{\delta} \leq \frac{\widetilde{\delta}}{C}, \tag{5.14}$$

for a sufficiently large constant $C$, determined in the proof of Lemma 5.6 below. Also, denote by

$$\|\omega\|_{\widetilde{B}_\tau} = \sum_{|\alpha|\geq 4} \frac{\tau^{|\alpha|-3}}{(|\alpha|-2)!} \widetilde{\delta}^{\alpha_1} \bar{\delta}^{\alpha_2} \|\partial^\alpha \omega\|_{L^2(\mathbb{T}\times[0,1])}$$

the corresponding dissipative analytic norm.

The following statement provides an estimate for the Euler vorticity in a uniform analytic norm up to the boundary.

LEMMA 5.3. *Assume that $\omega_0^{\mathrm{E}}$ satisfies (5.1), and suppose that $\omega^{\mathrm{E}}$ is a solution of the Euler equations with the initial data $\omega_0^{\mathrm{E}}$ such that (5.2), (5.3), and (5.6) hold for $t \in [0, T_0]$, for some constant $T_0 > 0$. With $\widetilde{\delta}, \bar{\delta}$ as in (5.14), the function $\omega^{\mathrm{E}}$ satisfies*

$$\sup_{0\leq t\leq T_0} \|\omega^{\mathrm{E}}(t)\|_{\widetilde{A}_{1/C}} \lesssim 1, \tag{5.15}$$

*where $C \geq 1$ is a sufficiently large constant.*

Recall that all constants depend only on $\lambda_0$. Note that since $\omega_0 \in H^4$, by the local existence theory for the Euler equations, by potentially reducing the value of the parameter $T_0$ from Lemma 5.2, we have

$$\|\omega^{\mathrm{E}}(t)\|_{H^4}, \|\boldsymbol{u}^{\mathrm{E}}(t)\|_{H^5} \lesssim 1, \qquad t \in [0, T_0]. \tag{5.16}$$

Before the proof of Lemma 5.3, we state two auxiliary results. In the first one, we show that the analytic norm in $\mathbb{T} \times [1/2, 1+\mu_0/2]$ of the Euler vorticity is bounded.

LEMMA 5.4. *Assume that $\omega_0^{\mathrm{E}}$ and $\omega^{\mathrm{E}}$ are is in Lemma 5.3, and let $\mu_0, T_0$ be as in Lemma 5.2 and (5.16). Then we have*

$$\sum_{|\alpha|\geq 3} \frac{\tau^{|\alpha|-3}}{(|\alpha|-3)!} \widetilde{\delta}^{\alpha_1} \bar{\delta}^{\alpha_2} \|\partial^\alpha \omega^{\mathrm{E}}\|_{L^2(\mathbb{T}\times[1/2,1+\mu_0/2])} \lesssim 1,$$

*for $t \in [0, T_0]$, provided $\widetilde{\delta}, \bar{\delta} \leq 1/C$ for a sufficiently large constant $C$.*

PROOF OF LEMMA 5.4. Fix $t \in [0, T_0]$, and denote $\omega = \omega^{\mathrm{E}}$. By (5.6), we have

$$\sum_{\xi\in\mathbb{Z}} e^{\lambda_0|\xi|/C} |\partial_y^i \omega_{0,\xi}(y)| \lesssim C^i i!, \qquad y \in [1/2, 1+\mu_0/2], \quad i \in \mathbb{N}_0, \tag{5.17}$$

omitting indicating the dependence on $t$. Therefore,

$$\begin{aligned}
\frac{1}{(|\alpha|-3)!} \|\partial_x^{\alpha_1} \partial_y^{\alpha_2} \omega\|_{L^2(\mathbb{T}\times[0,1])} &\lesssim \frac{1}{(|\alpha|-3)!} \sum_\xi \|\xi^{\alpha_1} \partial_y^{\alpha_2} \omega\|_{L_y^2(0,1)} \\
&\lesssim \frac{C^{\alpha_1}\alpha_1!}{(|\alpha|-3)!} \sum_\xi e^{\lambda\mu_0|\xi|/C} \|\partial_y^{\alpha_2}\omega\|_{L_y^2(0,1)} \lesssim \frac{C^{\alpha_1}\alpha_1!}{(|\alpha|-3)!} \sum_\xi e^{\lambda\mu_0|\xi|/C} \|\partial_y^{\alpha_2}\omega\|_{L_y^\infty(0,1)} \\
&\lesssim \frac{C^{\alpha_1}C^{\alpha_2}\alpha_1!\alpha_2!}{(|\alpha|-3)!} \lesssim \frac{C^{|\alpha|}|\alpha|!}{(|\alpha|-3)!} \lesssim C^{|\alpha|},
\end{aligned}$$

where we used (5.17) in the fourth inequality. □

In order to bound the analytic norm of the velocity by the vorticity in a strip (cf. Lemma 5.6 below), we first need to control the analyticity of $v^{\mathrm{E}}$ at $y = 1$. For $\tau > 0$, denote by

$$\|g\|_{\bar{A}_\tau} = \sum_{i\geq 2} \|\partial_x^i g\|_{H^{1/2}(\mathbb{T})} \frac{\tau^{i-2}}{(i-2)!}$$

the boundary analytic norm of a function $g$ defined on $\mathbb{T}$.

LEMMA 5.5. *Let $\omega^E$ be as in Lemma 5.3. Then we have*

$$\|v^{\mathrm{E}}|_{y=1}\|_{\bar{A}_{1/C}} \lesssim 1, \qquad t \in [0, T_0], \tag{5.18}$$

*for a sufficiently large constant $C$.*

PROOF OF LEMMA 5.5. As in the proof of Lemma 5.4, we have

$$\sum_{|\alpha|\geq 3} \frac{1}{C^{|\alpha|}} \frac{\tau^{|\alpha|-3}}{(|\alpha|-3)!} \|\partial^\alpha \omega^{\mathrm{E}}\|_{L^2(\mathbb{T}\times[\frac{1}{2},\frac{3}{2}])} \lesssim 1 \, . \tag{5.19}$$

where $C$ is a sufficiently large constant. Now, the component $v^E$ satisfies the elliptic equation

$$\Delta v^{\mathrm{E}} = \partial_x \omega^{\mathrm{E}} \, ,$$

and then the local elliptic analytic regularity, the bound (5.19), and the Sobolev estimate (5.2) imply

$$\sum_{\alpha\in\mathbb{N}_0^2} \|\partial^\alpha v^{\mathrm{E}}\|_{L^2(\mathbb{T}\times[\frac{3}{4},\frac{5}{4}])} \frac{1}{C^{|\alpha|}} \frac{\tau^{(|\alpha|-3)_+}}{(|\alpha|-3)!} \lesssim 1 \, , \tag{5.20}$$

with a possibly larger $C$. The bound (5.20) then gives (5.18) by using the trace inequality, upon enlarging the constant $C$. □

In the proof of Lemma 5.3, we need to estimate the velocity in terms of the vorticity in the analytic norm. It is important that we provide an estimate with the same analyticity radius; thus, simply appealing to the analytic regularity of the div-curl system is not sufficient.

LEMMA 5.6 (**Elliptic estimates in analytic spaces**). *For $t \in [0, T_0]$, denote $\omega = \omega^{\mathrm{E}}(t)$. Assume that*

$$\|\omega\|_{H^3}, \|\omega\|_{\widetilde{A}_\tau}, \|g\|_{\bar{A}_\tau} < \infty \, ,$$

*for some constant $\tau \in (0, 1]$. Then the function $\boldsymbol{u} = (u, v) = \boldsymbol{u}^{\mathrm{E}}$ is the solution of the elliptic system*

$$\begin{aligned} \operatorname{div} \boldsymbol{u} &= 0 \\ \operatorname{curl} \boldsymbol{u} &= \omega \, , \end{aligned}$$

*with the boundary conditions*

$$\begin{aligned} v|_{y=0} &= 0 \\ v|_{y=1} &= g \, , \end{aligned} \tag{5.21}$$

*and we have*

$$\|\boldsymbol{u}\|_{\widetilde{A}_\tau} \lesssim \|\omega\|_{H^3} + \|\omega\|_{\widetilde{A}_\tau} + \|g\|_{\bar{A}_\tau} \, , \tag{5.22}$$

*provided $\widetilde{\delta}$ and $\bar{\delta}$ satisfy (5.14) for a sufficiently large $C$.*

Applying (5.22) to (5.16) and (5.18), we get

$$\|\boldsymbol{u}^{\mathrm{E}}\|_{\widetilde{A}_{1/C}} \lesssim \|\omega^{\mathrm{E}}\|_{\widetilde{A}_{1/C}} + 1, \qquad t \in [0, T_0] \, , \tag{5.23}$$

where $C$ is sufficiently large.

PROOF OF LEMMA 5.6. We start with an estimate for $v$, which satisfies the Laplace equation

$$\Delta v = \partial_x \omega$$

with the boundary conditions (5.21). Denote

$$\phi(v) = \sum_{i+j\geq 3} \frac{\widetilde{\delta}^i \bar{\delta}^j \tau^{i+j-3}}{(i+j-3)!} \|\partial_x^i \partial_y^j v\|_{L^2} \, , \tag{5.24}$$

where, unless otherwise indicated, the norm is understood to be over the set $\mathbb{T} \times [0, 1]$. To treat the sum (5.24), we employ derivative reduction estimates as follows. For large values of $j$, we use

$$\|\partial_x^i \partial_y^j v\|_{L^2} \lesssim \|\partial_x^{i+1} \partial_y^{j-2} \omega\|_{L^2} + \|\partial_x^{i+1} \partial_y^{j-1} v\|_{L^2} + \|\partial_x^{i+1} \partial_y^{j-2} v\|_{L^2} + \|\partial_x^{i} \partial_y^{j-1} v\|_{L^2}, \qquad j \geq 2 \, , \tag{5.25}$$

while for small values,

$$\|\partial_x^i \partial_y v\|_{L^2} \lesssim \|\partial_x^i \omega\|_{L^2} + \|\partial_x^{i-1} g\|_{H^{3/2}(\Gamma)}, \qquad i \geq 2 \, ,$$

where $\Gamma = \{(x, y) : y = 1\}$ and

$$\|\partial_x^i v\|_{L^2} \lesssim \|\partial_x^{i-1} \omega\|_{L^2} + \|\partial_x^{i-2} g\|_{H^{3/2}(\Gamma)}, \qquad i \geq 3 \, ; \tag{5.26}$$

all three reductions (5.25)–(5.26) follow by using the $H^2$ elliptic regularity for the Laplacian. Now, we replace the inequalities (5.25)–(5.26) in the sum (5.24) according to the values of $j$ obtaining,

$$\phi(v) \lesssim \sum_{i+j\geq 3; j\geq 2} c_{ij}\big(\|\partial_x^{i+1}\partial_y^{j-2}\omega\|_{L^2} + \|\partial_x^{i+1}\partial_y^{j-1}v\|_{L^2} + \|\partial_x^{i+1}\partial_y^{j-2}v\|_{L^2} + \|\partial_x^{i}\partial_y^{j-1}v\|_{L^2}\big)$$
$$+ \sum_{i\geq 2} c_{i1}\big(\|\partial_x^{i}\omega\|_{L^2} + \|\partial_x^{i-1}g\|_{H^{3/2}(\Gamma)}\big) + \sum_{i\geq 3} c_{i0}\big(\|\partial_x^{i-1}\omega\|_{L^2} + \|\partial_x^{i-2}g\|_{H^{3/2}(\Gamma)}\big)\,,$$

where we denoted

$$c_{ij} = \frac{\widetilde{\delta}^i\bar{\delta}^j\tau^{i+j-3}}{(i+j-3)!}\,.$$

Next, we re-index the sums. All the terms involving $v$ may be absorbed into the left hand side under the condition (5.14), where $C$ is a sufficiently large constant, except for some lower order terms, which may be controlled by $\|v\|_{H^4}$. Thus we obtain

$$\phi(v) \lesssim \|\omega\|_{H^3} + \tau\|\omega\|_{\widetilde{A}_\tau} + \|g\|_{\bar{A}(\tau)}$$

since $\|v\|_{H^4} \lesssim \|\omega\|_{H^3}$, completing the inequality for $v$.

In order to treat the first component of the velocity, we split the sum $\phi(u)$ into the sums over regions $i \geq 1$ and $i = 0$. For the first sum, we use the divergence-free condition $\partial_x u = -\partial_y v$ and obtain

$$\sum_{i+j\geq 3; i\geq 1} c_{ij}\|\partial_x^i\partial_y^j u\|_{L^2} \lesssim \sum_{i+j\geq 3; i\geq 1} c_{ij}\|\partial_x^{i-1}\partial_y^{j+1}v\|_{L^2} \lesssim \|v\|_{H^3} + \frac{\widetilde{\delta}}{\bar{\delta}}\phi(v)\,, \tag{5.27}$$

while for $i = 0$, we use $\partial_y u = \partial_x v - \omega$ and write

$$\sum_{i+j\geq 3; i=0} c_{ij}\|\partial_x^i\partial_y^j u\|_{L^2} = \sum_{j\geq 3}\frac{\bar{\delta}^j\tau^{j-3}}{(j-3)!}\|\partial_y^j u\|_{L^2} \lesssim \sum_{j\geq 3}\frac{\bar{\delta}^j\tau^{j-3}}{(j-3)!}\|\partial_x\partial_y^{j-1}v\|_{L^2} + \sum_{j\geq 3}\frac{\bar{\delta}^j}{(j-3)!\tau^{j-3}}\|\partial_y^{j-1}\omega\|_{L^2}\,. \tag{5.28}$$

Summing (5.27) and (5.28), we get $\phi(u) \lesssim \|\omega\|_{H^3} + \tau\|\omega\|_{\widetilde{A}_\tau} + \|g\|_{\bar{A}(\tau)}$, and we obtain (5.22) for $u$. ☐

PROOF OF LEMMA 5.3. First, observe that the bound on the first term in (5.1) implies

$$\sum_{|\alpha|\geq 3}\frac{1}{C^{|\alpha|}(|\alpha|-3)!}\|\partial^\alpha\omega_0^{\mathrm{E}}\|_{L^2(\mathbb{T}\times[0,1])} \lesssim 1,$$

from where

$$\|\omega_0^{\mathrm{E}}\|_{\widetilde{A}_{\tau_0}} \lesssim 1, \tag{5.29}$$

regardless of the values of $\widetilde{\delta}, \bar{\delta} \in (0,1]$. Note that the solution $\omega = \omega^{\mathrm{E}}$ satisfies

$$\partial_t\omega + \boldsymbol{u}\cdot\nabla\omega = 0\,, \tag{5.30}$$

where $\boldsymbol{u} = (u,v) = (u^{\mathrm{E}}, v^{\mathrm{E}})$ is the Euler velocity. Let $\tau(t) = \tau_0 - Ct$, where $C \geq 1$ is a sufficiently large constant determined below. By the product rule, we have

$$\frac{d}{dt}\|\omega\|_{\widetilde{A}_\tau} = \tau'(t)\|\omega\|_{\widetilde{B}_\tau} + \sum_{|\alpha|\geq 3}\frac{\tau^{|\alpha|-3}}{(|\alpha|-3)!}\widetilde{\delta}^{\alpha_1}\bar{\delta}^{\alpha_2}\frac{d}{dt}\|\partial^\alpha\omega\|_{L^2(\mathbb{T}\times[0,1])}\,. \tag{5.31}$$

Next, we compute the time derivative of $\|\partial^\alpha\omega\|_{L^2(\mathbb{T}\times[0,1])}$. With $\alpha \in \mathbb{N}_0$ such that $|\alpha| \geq 3$, apply $\partial^\alpha$ to (5.30), multiply by $\partial^\alpha\omega$, and integrate by parts, obtaining

$$\begin{aligned}
\frac{1}{2}\frac{d}{dt}\|\partial^\alpha\omega\|_{L^2}^2 &= -\sum_{0<\beta\leq\alpha}\binom{\alpha}{\beta}\int_{\mathbb{T}\times[0,1]}(\partial^\beta\boldsymbol{u}\cdot\nabla\partial^{\alpha-\beta}\omega)\partial^\alpha\omega\,dxdy - \int_{\mathbb{T}\times[0,1]}\boldsymbol{u}\cdot\nabla\left(\frac{|\partial^\alpha\omega|^2}{2}\right)dxdy\\
&= -\sum_{0<\beta\leq\alpha}\binom{\alpha}{\beta}\int_{\mathbb{T}\times[0,1]}(\partial^\beta\boldsymbol{u}\cdot\nabla\partial^{\alpha-\beta}\omega)\partial^\alpha\omega\,dxdy - \int_{\mathbb{T}}v(t,x,1)\frac{|\partial^\alpha\omega(t,x,1)|^2}{2}dx\\
&\leq \sum_{0<\beta\leq\alpha}\binom{\alpha}{\beta}\|\partial^\beta\boldsymbol{u}\cdot\nabla\partial^{\alpha-\beta}\omega\|_{L^2}\|\partial^\alpha\omega\|_{L^2} + \mathcal{A}_\alpha(t)\,,
\end{aligned} \tag{5.32}$$

where

$$\mathcal{A}_\alpha(t) = -\frac{1}{2}\int_{\mathbb{T}} v(t,x,1)|\partial^\alpha\omega(t,x,1)|^2 dx$$

denotes the boundary term resulting from integration by parts. Since for all $(t,x)$ we have

$$\begin{aligned} -\frac{1}{2}v(t,x,1)|\partial^\alpha\omega(t,x,1)|^2 &= -\frac{1}{2}\int_0^1 \partial_y(v(\partial^\alpha\omega)^2)dy = -\frac{1}{2}\int_0^1 \big(\partial_y v(\partial^\alpha\omega)^2 + 2v\partial_y(\partial^\alpha\omega)\partial^\alpha\omega\big)dy \\ &\lesssim \|\nabla u\|_{L^\infty_y(0,1)}\|\partial^\alpha\omega\|^2_{L^2_y(0,1)} + \|u\|_{L^\infty(0,1)}\|\partial_y\partial^\alpha\omega\|_{L^2(0,1)}\|\partial^\alpha\omega\|_{L^2_y(0,1)}\,, \end{aligned}$$

we obtain by (5.16)

$$\begin{aligned} \mathcal{A}_\alpha(t) &\lesssim \|\nabla v\|_{L^\infty}\|\partial^\alpha\omega\|^2_{L^2} + \|v\|_{L^\infty}\|\partial_y\partial^\alpha\omega\|_{L^2}\|\partial^\alpha\omega\|_{L^2} \\ &\lesssim \|\partial^\alpha\omega\|^2_{L^2} + \|\partial_y\partial^\alpha\omega\|_{L^2}\|\partial^\alpha\omega\|_{L^2}\,. \end{aligned} \tag{5.33}$$

Combining (5.31), (5.32), and (5.33), we get

$$\begin{aligned} &\frac{d}{dt}\|\omega\|_{\widetilde{A}_\tau} - \tau'(t)\|\omega\|_{\widetilde{B}_\tau} \\ &\lesssim \|\omega\|_{\widetilde{A}_\tau} + \sum_{|\alpha|\ge 3}\frac{\tau^{|\alpha|-3}\widetilde{\delta}^{\alpha_1}\bar{\delta}^{\alpha_2}}{(|\alpha|-3)!}\sum_{0<\beta\le\alpha}\binom{\alpha}{\beta}\|\partial^\beta u\cdot\nabla\partial^{\alpha-\beta}\omega\|_{L^2} + \sum_{|\alpha|\ge 3}\|\partial_y\partial^\alpha\omega\|_{L^2}\frac{\tau^{|\alpha|-3}}{(|\alpha|-3)!}\,, \end{aligned} \tag{5.34}$$

on the interval $[0,T_0]\cap[0,\tau_0/C)$. Using the product rules for analytic norms as in [36], we obtain

$$\frac{d}{dt}\|\omega\|_{\widetilde{A}_\tau} - \tau'(t)\|\omega\|_{\widetilde{B}_\tau} \lesssim \|\omega\|_{\widetilde{A}_\tau} + (1+\|u\|_{H^2}+\|u\|_{\widetilde{A}_\tau})(\|\omega\|_{H^3}+\|\omega\|_{\widetilde{B}_\tau})\,.$$

This inequality, together with (5.16) and (5.23), leads to

$$\frac{d}{dt}\|\omega\|_{\widetilde{A}_\tau} - \tau'(t)\|\omega\|_{\widetilde{B}_\tau} \lesssim \|\omega\|_{\widetilde{A}_\tau} + (1+\|\omega\|_{\widetilde{A}_\tau})(1+\|\omega\|_{\widetilde{B}_\tau})\,.$$

Under the assumption

$$\tau'(t) + C\|\omega\|_{\widetilde{A}_\tau} \le 0\,,$$

where $C$ is a sufficiently large constant, we obtain

$$\frac{d}{dt}\|\omega\|_{\widetilde{A}_\tau} \lesssim 1 + \|\omega\|^2_{\widetilde{A}_\tau}\,.$$

Now, noting also that we have a bound (5.29) for $\omega_0^{\mathrm{E}}$, we conclude by a simple application of a Grönwall argument. □

PROOF OF THEOREM 5.1. Since the inequalities (5.2) and (5.3) are established in Lemma 5.2 above, we only need to prove (5.4) and (5.5). For simplicity, denote $\boldsymbol{u} = \boldsymbol{u}^{\mathrm{E}}$ and $\omega = \omega^{\mathrm{E}}$. By (5.15) in Lemma 5.3, there exists a constant $\zeta_0$ such that

$$\sum_{|\alpha|\ge 4}\frac{\zeta_0^{|\alpha|-4}}{(|\alpha|-2)!}\widetilde{\delta}^{\alpha_1}\bar{\delta}^{\alpha_2}\|\partial^\alpha\omega\|_{L^2(\mathbb{T}\times[0,1])} \lesssim 1, \qquad t\in[0,T_1]\,. \tag{5.35}$$

Since $\widetilde{\delta}$ and $\bar{\delta}$ are constants, we may reduce $\zeta_0$ to obtain

$$\sum_{|\alpha|\ge 4}\frac{\zeta_0^{|\alpha|-4}}{(|\alpha|-2)!}\|\partial^\alpha\omega\|_{L^2(\mathbb{T}\times[0,1])} \lesssim 1, \qquad t\in[0,T_1]\,.$$

Finally, we may use Agmon's inequality to bound $\|\omega\|_{L^\infty(\mathbb{T}\times[0,t])}$ in terms of the $L^2$ norms and further decrease $\zeta_0$ to get (5.4) for $t\in[0,T_1]$. Finally, by (5.23) and (5.35), we get

$$\sum_{m\ge 4}\frac{\zeta_0^{m-4}}{(m-2)!}\sum_{|\alpha|=m}\widetilde{\delta}^{\alpha_1}\bar{\delta}^{\alpha_2}\|\partial^\alpha\boldsymbol{u}\|_{L^2(\mathbb{T}\times[0,1])} \lesssim 1, \qquad t\in[0,T_1]\,,$$

from where, using the same arguments as for the vorticity, we obtain (5.5). □

From Theorem 5.1, we obtain the next statement. The bounds (5.36)–(5.38) in the theorem are used when estimating the remainder of the Prandtl asymptotic expansions.

LEMMA 5.7. *Assume that* (5.1) *holds. Then there exist constants* $T_1 \in (0,1]$, $\lambda_1 \in (0, \lambda_0/2]$, *and* $\mu_1 \in (0, \mu_0]$, *such that for all* $\lambda \in [0, \lambda_1]$, $\mu \in [0, \mu_1]$, *and all* $t \in [0, T_1]$, *we have for the Euler vorticity*

$$\|\partial_x^i \partial_y^j \omega^{\mathrm{E}}\|_{Y_{\lambda,\mu,\infty}} + \|\partial_x^i \partial_y^j \omega^{\mathrm{E}}\|_{Y_{\lambda,\mu}} \lesssim 1\,, \tag{5.36}$$

*for the first velocity component there holds*

$$\|\partial_x^i \partial_y^j u^{\mathrm{E}}\|_{Y_{\lambda,\mu,\infty}} + \left\| \frac{\partial_x^i u^{\mathrm{E}} - \partial_x^i U^{\mathrm{E}}}{y} \right\|_{Y_{\lambda,\mu,\infty}} \lesssim 1\,, \tag{5.37}$$

*while for the second velocity component we have*

$$\|\partial_x^i \partial_y^j v^{\mathrm{E}}\|_{Y_{\lambda,\mu,\infty}} + \left\| \frac{1}{y} \partial_x^i v^{\mathrm{E}} \right\|_{Y_{\lambda,\mu,\infty}} + \left\| \frac{\partial_x^i v^{\mathrm{E}} + y \partial_x^{i+1} U^{\mathrm{E}}}{y^2} \right\|_{Y_{\lambda,\mu,\infty}} \lesssim 1\,, \tag{5.38}$$

*for all* $i + j \leq 3$, *where the implicit constants depend on* $i$ *and* $j$. *Moreover, for the Euler trace* $U^{\mathrm{E}}$ *defined in* (2.3), *we have*

$$\sum_{\xi} e^{\lambda|\xi|} |U^{\mathrm{E}}_\xi| \lesssim 1\,, \tag{5.39}$$

*for* $t \in [0, T_1]$ *and* $\lambda \in [0, \lambda_1]$.

PROOF OF LEMMA 5.7. Let $T_1 \in (0,1]$ be the constant $T_0$ in (5.15). For simplicity of presentation, we shall establish the inequality (5.36) for the first term when $i = j = 0$. The general case, as well as the inequalities (5.37)–(5.39), follow from Theorem 5.1 in the same way. Using the definition of the $Y_{\lambda,\mu,\infty}$ norm, we need to prove

$$\sum_{\xi \in \mathbb{Z}} e^{\lambda(1+\mu)|\xi|} |\omega^{\mathrm{E}}_\xi(t,y)| \lesssim 1\,, \tag{5.40}$$

for $y \in \Omega_\mu \cap \{\mathbb{R}\mathrm{e}\, y \leq 1\}$ and $t \in [0, T_1]$, where $\lambda$ and $\mu$ are sufficiently small constants. Fix $t \in [0, T_1]$. For $j \in \mathbb{N}_0$, denote

$$a_{j,\xi} = \sup_{0 \leq y \leq 1+\mu} |\partial_y^j \omega^{\mathrm{E}}_\xi(t,y)|\,.$$

(Note that the supremum is taken among the real values of $y$.) Using Agmon's inequality in the variable $y$ and (5.4), we have

$$|\partial_y^j \omega^{\mathrm{E}}_\xi(t,y)| \lesssim \|\partial_y^j \omega^{\mathrm{E}}_\xi\|_{L^2_y(0,1+\mu)} + \|\partial_y^{j+1} \omega^{\mathrm{E}}_\xi\|_{L^2_y(0,1+\mu)} \lesssim \frac{j!}{\zeta_0^j} + \frac{(j+1)!}{\zeta_0^{j+1}} \lesssim \frac{(j+1)!}{\zeta_0^j} \lesssim \frac{2^j j!}{\zeta_0^j}\,,$$

for $j \in \mathbb{N}_0$. Therefore, for $R_0 \leq \zeta_0/2$, we obtain the bound

$$\sum_{j=0}^{\infty} \frac{a_{j,\xi}}{j!} R_0^j \lesssim 1\,. \tag{5.41}$$

Next, define

$$f_\xi(t,y) = \sum_{j=0}^{\infty} \frac{\partial_y^j \omega^{\mathrm{E}}_\xi(t, \mathbb{R}\mathrm{e}\, y)}{j!} (y - \mathbb{R}\mathrm{e}\, y)^j\,, \qquad |y - \mathbb{R}\mathrm{e}\, y| < R_0, \quad \xi \in \mathbb{Z}\,.$$

By (5.41), we have

$$|f_\xi(t,y)| \leq \sum_{j=0}^{\infty} \frac{a_j}{j!} R_0^j \lesssim 1\,,$$

and thus the function $f_\xi(t,y)$ is holomorphic in the region

$$S_0 = \bigl\{ y \in \mathbb{C} : |\mathbb{I}\mathrm{m}\, y| \leq R_0, 0 \leq \mathbb{R}\mathrm{e}\, y < 1 \bigr\} \cup \bigl\{ y \in \mathbb{C} : |\mathbb{I}\mathrm{m}\, y| \leq 1 - R_0, 1 \leq \mathbb{R}\mathrm{e}\, y < 1 + R_0 \bigr\}\,.$$

Since $f_\xi(t,y) = \partial_y \omega^{\mathrm{E}}_\xi(t,y)$ on the segment $[0,1]$, by unique analytic continuation, we have

$$f_\xi(t,y) = \partial_y \omega^{\mathrm{E}}_\xi(t,y) \qquad \text{on} \quad S_0 \cap \Omega_\mu\,.$$

Now, choose $\mu_1$ sufficiently small so that the domain $\Omega_{\mu_1}$ lies inside the region $S_0$. For $y \in \Omega_{\mu_1}$, we then have

$$\sum_{\xi\in\mathbb{Z}} e^{\lambda(1+\mu_1)|\xi|}|\omega^{\mathrm{E}}_\xi(t,y)| \lesssim \sum_{\xi\in\mathbb{Z}} e^{\lambda(1+\mu_1)|\xi|}\sum_{j=0}^\infty \frac{a_{j,\xi}}{j!}R_0^j$$
$$\lesssim \sum_{i=0}^\infty \frac{(\lambda(1+\mu_1))^i|\xi|^i}{i!}\sum_{j=0}^\infty \frac{a_{j,\xi}}{j!}R_0^j \lesssim \sum_{i,j=0}^\infty \frac{|\xi|^i a_{j,\xi}}{(i+j)!}(\lambda(1+\mu_1))^i R_0^j \lesssim 1\,,$$

and the inequality (5.40) is proven provided $\lambda$ and $\mu$ are sufficiently small constants. □

## 6. Size of the Prandtl solution in analytic norms

The initial datum for the Prandtl equation $(\tilde{u}^{\mathrm{P}}_0, \bar{v}^{\mathrm{P}}_0)$ is given by the boundary layer part of the Navier-Stokes initial datum, cf. (2.13)–(2.14). In view of the definitions (2.6), (2.7), and (2.10), this initial Prandtl velocity may be computed from the tangential Euler trace $U^{\mathrm{E}}_0$ (which is known; cf. Section 5), and from the initial Prandtl vorticity $\Omega^{\mathrm{P}}_0$. We assume that the initial Prandtl vorticity is real-analytic and satisfies

$$\|\Omega^{\mathrm{P}}_0\|_{A_{\lambda_0/2}} \lesssim 1\,, \tag{6.1}$$

with $\lambda_0 > 0$ as in (5.1), and where we denote the analytic norm $A_\tau$ as

$$\|\Omega^{\mathrm{P}}\|^2_{A_\tau} = \sum_{|\alpha|\geq 0} \frac{\tau^{2|\alpha|}\kappa^{2\alpha_2}}{(|\alpha|-4)!^2}\|(1+Y)^\gamma Y^{\alpha_2} D^\alpha \Omega^{\mathrm{P}}\|^2_{L^2}. \tag{6.2}$$

At this stage, we also introduce a dissipative analytic norm $B_\tau$, given by

$$\|\Omega^{\mathrm{P}}\|^2_{B_\tau} = \sum_{|\alpha|\geq 5} \frac{|\alpha|\tau^{2|\alpha|}\kappa^{2\alpha_2}}{(|\alpha|-4)!^2}\|(1+Y)^\gamma Y^{\alpha_2} D^\alpha \Omega^{\mathrm{P}}\|^2_{L^2}\,. \tag{6.3}$$

The parameter $\kappa \in (0,1]$ is introduced in order to deal with the dissipative term $\partial_{YY}$ in the analytic estimate for the Prandtl system; one may for instance set $\kappa = 1/8$. The parameter $\tau > 0$ is related to the analyticity radius of $\Omega^{\mathrm{P}}$.

REMARK 6.1 (**Example of a compatible initial datum**). An example of a compatible Prandtl datum is given by (2.15), so that the initial vorticity equals

$$\Omega^{\mathrm{P}}_0(x,Y) = U^{\mathrm{E}}_0(x)\varphi''(Y)\,, \tag{6.4}$$

where the function $\varphi$ in (6.4) is assumed to satisfy

$$\sum_{n\geq 0} \frac{(\lambda_0\kappa)^{2n}}{(n-4)!^2}\|(1+Y)^\gamma Y^n \partial_Y^{n+2}\varphi\|^2_{L^2([0,\infty))} \lesssim 1\,, \tag{6.5}$$

and the parameter $\lambda_0$ is as in (5.1). With $\varphi$ satisfying (6.5) and with the assumption (5.1) for $\omega^{\mathrm{E}}_0$, which implies via the Biot-Savart law $U^{\mathrm{E}}_{0,\xi} = \int_0^\infty e^{-|\xi|z}\omega^{\mathrm{E}}_{0,\xi}(z)dz$ (see e.g. (7.35) with $g \equiv 0$ evaluated at $y = 0$) that $U^{\mathrm{E}}_0$ is real-analytic with respect to $x$ with radius $\lambda_0$, we obtain that $\Omega^{\mathrm{P}}_0$ in the definition (6.4) satisfies the condition (6.1).

Having assumed in (6.1) that the initial Prandtl vorticity is real-analytic, and since in Lemma 5.7 we have already shown that the Euler trace $U^{\mathrm{E}}$ is real-analytic on $[0,T_1]$, by using analytic energy estimates similar to those in [**35**] and [**37**] we may show that there exists $T_2 \in (0,T_1]$ and a real-analytic solution of the Prandtl system (2.8)–(2.9) on $[0,T_2]$. More precisely, in light of (6.1) and (5.39), we may set

$$\tau_0 = \frac{1}{2}\min\left\{\frac{\lambda_0}{2},\lambda_1\right\} = \frac{\lambda_1}{2}$$

and conclude that there exists $T_2 \in (0,T_1]$ and an analytic solution $\Omega^{\mathrm{P}}$ to the Prandtl equation (2.8)–(2.9) with analyticity properties quantified in the following way. There exists a decreasing function $\tau = \tau(t)$ (different than the one from Section 5) on $[0,T_2]$ such that $\tau(0) = \tau_0$ and

$$\tau(t) \geq \tau(T_2) \geq \frac{\tau_0}{2} = \frac{\lambda_1}{4}, \qquad t \in [0,T_2]\,, \tag{6.6}$$

with $\Omega^{\rm P}$ satisfying

$$\sup_{t\in[0,T_2]}\|\Omega^{\rm P}(t)\|_{A_{\tau(t)}}^2+\int_0^{T_2}\left(\|\partial_Y\Omega^{\rm P}(\tau)\|_{A_{\tau(t)}}^2+\|\Omega^{\rm P}(\tau)\|_{B_{\tau(t)}}^2\right)dt\lesssim 1\,. \tag{6.7}$$

The term involving $\|\partial_Y\Omega^{\rm P}(\tau)\|_{A_{\tau(t)}}^2$ results from the dissipation $\partial_{YY}\Omega^{\rm P}$ in (2.8), while the one with $\|\Omega^{\rm P}(\tau)\|_{B_{\tau(t)}}^2$ from the decay in analyticity radius. Note that since all constants are allowed to depend on $\lambda_1$, and since the lower bound (6.6) holds, we have

$$T_2\sim 1 \qquad\text{and}\qquad \tau(t)\sim\lambda_1\sim 1\,.$$

While the bound (6.7) provides analytic estimates for the Prandtl solution, these estimates are with respect to the $A_\tau$ and $B_\tau$ energy-type norms from (6.2)–(6.3). However, in order to bound the error vorticity $\omega_{\rm e}$, which is forced by the Prandtl solution via (3.17)–(3.19), we need to estimate the size of the Prandtl solution in the norm $P_{\lambda,\mu,\infty}$. This is achieved in the next statement.

LEMMA 6.1. *Let $\lambda_2=\frac{\lambda_1}{32}$, $\mu_2=\frac{\lambda_1\kappa}{32}\le 1$, and $\gamma\ge 4$, and assume that (6.7) holds. Then, for any $\lambda\in(0,\lambda_2]$, any $\mu\in(0,\mu_2]$, and for all $i,j\in\mathbb{N}_0$ the following bounds hold. For the classical Prandtl vorticity we have the pointwise in time estimates*

$$\|(1+Y)^{\gamma-1}Y^{j+1}\partial_x^i\partial_Y^j\Omega^{\rm P}\|_{P_{\lambda,\mu,\infty}}\lesssim 1\,, \tag{6.8}$$

*for the first component of the Prandtl velocity we have*

$$\|(1+Y)^{\gamma-\frac32}\partial_x^i\widetilde{u}^{\rm P}\|_{P_{\lambda,\mu,\infty}}+\sum_{\xi\in\mathbb{Z}}e^{\lambda(1+\mu)|\xi|}\int_0^\infty\left|(\partial_x^i\widetilde{u}^{\rm P})_\xi\right|dY\lesssim 1\,, \tag{6.9}$$

*while for the second component of the velocity*

$$\left\|\frac{\partial_x^i v^{\rm P}}{Y}\right\|_{P_{\lambda,\mu,\infty}}+\|(1+Y)^{\gamma-\frac52}\partial_x^i\bar{v}^{\rm P}\|_{P_{\lambda,\mu,\infty}}\lesssim 1\,, \tag{6.10}$$

*uniformly on $[0,T_2]$, where the implicit constants are allowed to depend on $i$ and $j$. In addition, we have the integrated in time estimate*

$$\int_0^{T_2}\|(1+Y)^{\gamma}Y^j\partial_x^i\partial_Y^j\Omega^{\rm P}\|_{P_{\lambda,\mu,\infty}}^4\,dt\lesssim 1\,, \tag{6.11}$$

*for $i,j\in\mathbb{N}_0$.*

Observe that the derivative $\partial_y^j$ is matched in (6.8) by the weight $Y^{j+1}$, while in (6.11) with $Y^j$.

PROOF OF LEMMA 6.1. Since the $P_{\lambda,\mu,\infty}$ norm is monotone in $\mu$ and $\lambda$, we assume throughout the proof that

$$\lambda=\lambda_2\qquad\text{and}\qquad \mu=\mu_2=\kappa\lambda_2\,,$$

and thus in view of (6.6) we have $8\lambda\le\tau(t)/2$, for any $t\in[0,T_2]$. It suffices to establish the bounds claimed in the lemma for the case $i=0$, as the cases $i\ge 1$ follow analogously (these bounds carry an additional factor of $\lambda_1^{-i}$, but since $\lambda_1\sim 1$, these factors are hidden in the implicit constant).

We start by establishing the $P_{\lambda,\mu,\infty}$ bounds for the first term in (6.8) by proving

$$\|(1+Y)^{\gamma-1}Y^{j+1}\partial_Y^j\Omega^{\rm P}\|_{P_{\lambda,\mu,\infty}}\lesssim 1\,. \tag{6.12}$$

For any weight function $\eta(Y)=\eta(\mathbb{R}{\rm e}\,Y)$, and a function $f$ which is analytic with respect to $Y$ in the domain $\widetilde{\Omega}_\mu$, form the Taylor series expansion for $f(Y)=f(\mathbb{R}{\rm e}\,Y+i\mathbb{I}{\rm m}\,Y)$ around $f(\mathbb{R}{\rm e}\,Y)$, and using that $|\mathbb{I}{\rm m}\,Y|\le\mu\mathbb{R}{\rm e}\,Y$ for $Y\in\widetilde{\Omega}_\mu$, we obtain

$$\sup_{\widetilde{\Omega}_\mu}|\eta(Y)f(Y)|\lesssim\sum_{m\ge0}\frac{1}{m!}\|\eta(Y)(\mu\mathbb{R}{\rm e}\,Y)^m\partial_Y^m f(Y)\|_{L_Y^\infty([0,\infty))}\,. \tag{6.13}$$

Applying this inequality with $f=(\partial_Y^j\Omega^{\rm P})_\xi$ and $\eta(Y)=(1+Y)^{\gamma-1}Y^{j+1}$ (for simplicity of notation we write $Y$ instead of $\mathbb{R}{\rm e}\,Y$ throughout this proof), we deduce

$$\sup_{\widetilde{\Omega}_\mu}\left|(1+Y)^{\gamma-1}Y^{j+1}(\partial_Y^j\Omega^{\rm P})_\xi\right|\lesssim\sum_{m\ge0}\frac{\mu^m}{m!}\|(1+Y)^{\gamma-1}Y^{m+j+1}(\partial_Y^{m+j}\Omega^{\rm P})_\xi\|_{L_Y^\infty([0,\infty))}\,.$$

Next, for a fixed $\xi \in \mathbb{Z}$, by expanding $e^{\lambda(1+\mu)|\xi|} \le e^{2\lambda|\xi|}$ into its power series, and using $(m+n)!/m!n! \le 2^{m+n}$ and $\mu = \lambda\kappa$, we get

$$\begin{aligned} e^{\lambda(1+\mu)|\xi|} \sup_{\widetilde{\Omega}_\mu} \left|(1+Y)^{\gamma-1} Y^{j+1} (\partial_Y^j \Omega^{\rm P})_\xi\right| &\lesssim \sum_{m,n\ge 0} \frac{(2\lambda)^n \mu^m}{m!n!} \|(1+Y)^{\gamma-1} Y^{m+j+1} (\partial_x^n \partial_Y^{m+j} \Omega^{\rm P})_\xi\|_{L^\infty_Y} \\ &\lesssim \sum_{m,n\ge 0} \frac{(4\lambda)^{m+n}\kappa^m}{(m+n)!} \|(1+Y)^{\gamma-1} Y^{m+j+1} (\partial_x^n \partial_Y^{m+j} \Omega^{\rm P})_\xi\|_{L^\infty_Y} . \end{aligned} \tag{6.14}$$

Taking the $\ell^1(\mathbb{Z})$ norm in $\xi$, and using a factor of $(1+\xi^2)2^{m+n}$ in order to obtain a bound in $\ell^2(\mathbb{Z}\times\mathbb{N}^2)$ in $(\xi, n, m)$ needed for Plancherel's identity, we estimate

$$\begin{aligned} &\|(1+Y)^{\gamma-1} Y^{j+1} \partial_Y^j \Omega^{\rm P}\|^2_{P_{\lambda,\mu,\infty}} \\ &\lesssim \sum_{\xi\in\mathbb{Z}} \sum_{n,m\ge 0} \frac{(1+\xi^2)(8\lambda)^{2(m+n)}\kappa^{2m}}{(m+n)!^2} \|(1+Y)^{\gamma-1} Y^{m+j+1} (\partial_x^n \partial_Y^{m+j} \Omega^{\rm P})_\xi\|^2_{L^\infty_Y} \\ &\lesssim \sum_{\xi\in\mathbb{Z}} \sum_{n,m\ge 0} \frac{(1+\xi^2)(8\lambda)^{2(m+n)}\kappa^{2m}}{(m+n)!^2} \|(1+Y)^{\gamma} Y^{m+j} (\partial_x^n \partial_Y^{m+j} \Omega^{\rm P})_\xi\|_{L^2_Y} \\ &\quad \times \Big( \|(1+Y)^{\gamma-1} Y^{m+j+1} (\partial_x^n \partial_Y^{m+j+1} \Omega^{\rm P})_\xi\|_{L^2_Y} + (m+j+1) \|(1+Y)^{\gamma-1} Y^{m+j} (\partial_x^n \partial_Y^{m+j} \Omega^{\rm P})_\xi\|_{L^2_Y} \Big) \\ &\lesssim \Big( \sum_{\xi\in\mathbb{Z}} \sum_{n,m\ge 0} \frac{(1+\xi^2)(8\lambda)^{2(m+n)}\kappa^{2m}}{(m+n)!^2} \|(1+Y)^{\gamma} Y^{m+j} (\partial_x^n \partial_Y^{m+j} \Omega^{\rm P})_\xi\|^2_{L^2_Y} \Big)^{1/2} \\ &\quad \times \Big( \sum_{\xi\in\mathbb{Z}} \sum_{n,m\ge 0} \frac{(1+\xi^2)(8\lambda)^{2(m+n)}\kappa^{2m}}{(m+n)!^2} \|(1+Y)^{\gamma-1} Y^{m+j+1} (\partial_x^n \partial_Y^{m+j+1} \Omega^{\rm P})_\xi\|^2_{L^2_Y} \\ &\qquad + \sum_{\xi\in\mathbb{Z}} \sum_{n,m\ge 0} \frac{(1+\xi^2)(8\lambda)^{2(m+n)}\kappa^{2m}}{(m+n)!^2} (m+j+1)^2 \|(1+Y)^{\gamma-1} Y^{m+j} (\partial_x^n \partial_Y^{m+j} \Omega^{\rm P})_\xi\|^2_{L^2_Y} \Big)^{1/2} . \end{aligned} \tag{6.15}$$

In the second inequality, we used Agmon's inequality in $Y$, along with the fact that $(1+Y)^{\gamma-1} Y^{m+j+1} \partial_y^{m+j} \Omega^{\rm P}$ vanishes at $Y = 0$ (recall that $m + j \ge 0$). Therefore, by Parseval's identity in the variable $x$,

$$\begin{aligned} &\|(1+Y)^{\gamma-1} Y^{j+1} \partial_Y^j \Omega^{\rm P}\|^2_{P_{\lambda,\mu,\infty}} \\ &\lesssim \Big( \sum_{n,m\ge 0} \frac{(\tau/2)^{2(m+n)}\kappa^{2m}}{(m+n)!^2} \|(1+Y)^{\gamma} Y^{m+j} (1-\partial_x^2) \partial_x^n \partial_Y^{m+j} \Omega^{\rm P}\|^2_{L^2} \Big)^{1/2} \\ &\quad \times \Big( \sum_{n,m\ge 0} \frac{(\tau/2)^{2(m+n)}\kappa^{2m}}{(m+n)!^2} \|(1+Y)^{\gamma-1} Y^{m+j+1} (1-\partial_x^2) \partial_x^n \partial_Y^{m+j+1} \Omega^{\rm P}\|^2_{L^2} \\ &\qquad + \sum_{n,m\ge 0} \frac{(\tau/2)^{2(m+n)}\kappa^{2m}}{(m+n)!^2} (m+j+1)^2 \|(1+Y)^{\gamma-1} Y^{m+j} (1-\partial_x^2) \partial_x^n \partial_Y^{m+j} \Omega^{\rm P}\|^2_{L^2} \Big)^{1/2} . \end{aligned} \tag{6.16}$$

Now, for $\alpha = (n, m)$ we use that $|\alpha|^r 2^{-|\alpha|} \lesssim_r 1$ holds whenever $r \ge 0$, and since $\tau \sim 1$ (meaning that our constants are allowed to depend on $\lambda_1$), that $\kappa = 1/8 \sim 1$, we obtain from (6.16) that

$$\|(1+Y)^{\gamma-1} Y^{j+1} \partial_Y^j \Omega^{\rm P}\|_{P_{\lambda,\mu,\infty}} \lesssim \left\|\Omega^{\rm P}\right\|_{A_\tau} , \tag{6.17}$$

where the implicit constant also depends on $i$ and $j$. The bound (6.12) now follows.

Next, we consider the bound (6.11), which is proven similarly to the arguments above, but with $Y^j$ and $(1+Y)^j$ replacing $Y^{j+1}$ and $(1+Y)^{\gamma-1}$, respectively. Agmon's inequality in $Y$ here reads

$$\begin{aligned}&\|(1+Y)^\gamma Y^{m+j}(\partial_x^n\partial_Y^{m+j}\Omega^{\rm P})_\xi\|_{L_Y^\infty}^2\\&\quad\lesssim\|(1+Y)^\gamma Y^{m+j}(\partial_x^n\partial_Y^{m+j}\Omega^{\rm P})_\xi\|_{L_Y^2}^2\\&\qquad+\|(1+Y)^\gamma Y^{m+j}(\partial_x^n\partial_Y^{m+j}\Omega^{\rm P})_\xi\|_{L_Y^2}\\&\qquad\quad\times\left(\|(1+Y)^\gamma Y^{m+j}(\partial_x^n\partial_Y^{m+j+1}\Omega^{\rm P})_\xi\|_{L_Y^2}+(m+j)\|(1+Y)^\gamma Y^{m+j-1}(\partial_x^n\partial_Y^{m+j}\Omega^{\rm P})_\xi\|_{L_Y^2}\right).\end{aligned}\tag{6.18}$$

When compared to (6.15), the main difference in (6.18) is that the terms on the last line contain factors of the type $Y^{\alpha_2}D^\alpha\partial_y\Omega^{\rm P}$, and thus we bound

$$\|(1+Y)^\gamma Y^j\partial_Y^j\Omega^{\rm P}\|_{P_{\lambda,\mu,\infty}}\lesssim\left\|\Omega^{\rm P}\right\|_{A_\tau}+\left\|\Omega^{\rm P}\right\|_{A_\tau}^{1/2}\left\|\partial_Y\Omega^{\rm P}\right\|_{A_\tau}^{1/2}.$$

The second term in the above inequality may only be estimated in $L^4$ in time, by appealing to the bound provided by the second term in (6.7); from this, the estimate (6.11) follows.

Next, we turn to the proof of the bound for the first term on the left side of (6.9). Using (6.13) with $f(y)=\widetilde{u}^{\rm P}$ and proceeding in the same way as in the first line of (6.15), we have

$$\begin{aligned}&\|(1+Y)^{\gamma-\frac32}\widetilde{u}^{\rm P}\|_{P_{\lambda,\mu,\infty}}^2\\&\quad\lesssim\sum_{\xi\in\mathbb{Z}}\sum_{n,m\ge0}\frac{(1+\xi^2)(8\lambda)^{2(m+n)}\kappa^{2m}}{(m+n)!^2}\|(1+Y)^{\gamma-\frac32}Y^m(\partial_x^n\partial_Y^m\widetilde{u}^{\rm P})_\xi\|_{L_Y^\infty}^2=I_1+I_2,\end{aligned}\tag{6.19}$$

where $I_1$ and $I_2$ correspond to the sums with $m=0$ and $m\ge1$, respectively. In order to estimate the first sum, we use the fundamental theorem of calculus on $[Y,\infty)$ and $\partial_Y\widetilde{u}^{\rm P}=\Omega^{\rm P}$ to obtain

$$\|(1+Y)^{\gamma-\frac12}\partial_x^n\widetilde{u}^{\rm P}\|_{L_Y^\infty}\lesssim\|(1+Y)^\gamma\partial_x^n\Omega^{\rm P}\|_{L_Y^2}$$

and thus, since $(1+Y)^{-1}\in L_Y^2$, we have

$$\|(1+Y)^{\gamma-\frac32}\partial_x^n\widetilde{u}^{\rm P}\|_{L_Y^2}\lesssim\|(1+Y)^\gamma\partial_x^n\Omega^{\rm P}\|_{L_Y^2}.$$

Therefore, using that $8\lambda\le\tau/2$ and $\tau\sim1$, we may use Plancherel's identity and (6.7) to obtain

$$I_1\lesssim\sum_{\xi\in\mathbb{Z}}\sum_{n\ge0}\frac{(1+\xi^2)(\tau/2)^{2n}}{n!^2}\|(1+Y)^\gamma(\partial_x^n\Omega^{\rm P})_\xi\|_{L_Y^2([0,\infty))}^2\lesssim\left\|\Omega^{\rm P}\right\|_{A_\tau}^2\lesssim1.$$

The bound for the $I_2$ term in (6.19) is more direct, and is obtained by replacing $\partial_Y\widetilde{u}^{\rm P}=\Omega^{\rm P}$ and repeating the proof of (6.17). This implies that $I_2\lesssim1$, and thus $\|(1+Y)^{\gamma-\frac32}\widetilde{u}^{\rm P}\|_{P_{\lambda,\mu,\infty}}\lesssim1$ holds.

The estimate on the first term in (6.10) follows from

$$\frac{1}{Y}v^{\rm P}(x,Y)=-\frac1Y\int_0^Y\partial_xu^{\rm P}(x,Y')dY'=\partial_xU^{\rm E}(x)-\frac1Y\int_0^Y\partial_x\widetilde{u}^{\rm P}(x,Y')dY',$$

the previously established bound (5.37) (which holds for a wider set of values for $\lambda,\mu$), the bound on the first term in (6.9), and the fact that $\gamma\ge3/2$.

To bound the second term in (6.10), we recall the identity

$$\widetilde{v}^{\rm P}(x,Y)=\int_Y^\infty\partial_x\widetilde{u}^{\rm P}(x,Y')dY',$$

which may be used in conjunction with the bound for the first term on the left side of (6.9), and integration in $Y$ (which is possible since $\gamma>5/2$), to yield the desired bound for the third term in (6.9).

In order to conclude the proof of the lemma, we need to estimate the second term on the right side of (6.9). For this, we have

$$\begin{aligned}\sum_{\xi\in\mathbb{Z}}e^{\lambda(1+\mu)|\xi|}\int_0^\infty\left|\widetilde{u}_\xi^{\rm P}\right|dY&\lesssim\int_0^\infty\sum_{\xi\in\mathbb{Z}}e^{\lambda(1+\mu)|\xi|}(1+Y)^{\gamma-\frac32}\left|\widetilde{u}_\xi^{\rm P}\right|(1+Y)^{\frac32-\gamma}dY\\&\lesssim\int_0^\infty\|(1+Y)^{\gamma-\frac32}\widetilde{u}^{\rm P}\|_{P_{\lambda,\mu,\infty}}(1+Y)^{\frac32-\gamma}dY\lesssim1,\end{aligned}$$

by appealing to the bound for the first term on the left side of (6.9), and the condition $\gamma>5/2$. ☐

We conclude this section by noting that the estimates obtained in Lemma 6.1 are all with respect to norms that are (weighted) $L^\infty_Y$. On the other hand, the a-priori bound (6.7) provides $L^2_Y$ information, and this may be used to improve the $Y_{\lambda,\mu}$ product estimate (4.11b), which in essence is an $L^1_y$ bound. In this direction we have the following.

LEMMA 6.2 (**Improved $Y_{\lambda,\mu}$ product estimate involving the Prandtl vorticity**). *Let $\lambda, \mu$ be as in Lemma 6.1, and assume that $g = g(x,y)$ is such that $\|g\|_{Y_{\lambda,\mu,\infty}} < \infty$. Then, we have the pointwise in time estimate*

$$\left\| g(x,y) Y^j \partial_x^i \partial_Y^j \Omega^{\mathrm{P}}(x,Y) \right\|_{Y_{\lambda,\mu}} \lesssim \epsilon \left\| g(x,y) \right\|_{Y_{\lambda,\mu,\infty}} , \tag{6.20}$$

*for any $i,j \in \mathbb{N}_0$.*

In comparison, (6.7) and (6.11) give a bound similar to (6.20), but which is valid only in $L^4$ with respect to time, as opposed to pointwise in time.

PROOF OF LEMMA 6.2. The statement follows from the first inequality in (4.14), if we are able to show that

$$\sum_\xi e^{\lambda(1+\mu)|\xi|} \left\| Y^j \partial_x^i \partial_Y^j \Omega^{\mathrm{P}}_\xi(Y) \right\|_{\mathcal{L}^1_\mu} \lesssim \epsilon \, . \tag{6.21}$$

Recall that the weight $Y^j$ in (6.21) is short hand notation for $(\mathbb{R}\mathrm{e}\, Y)^j$. At this stage we recall the definition of the $\mathcal{L}^1_\mu$ norm in (4.3), and note that this consists of $L^1$ norms over complex paths corresponding to the variable $y = \epsilon Y$. Moreover, we note that if $y \in \Omega_\mu$, then by the definitions (4.1) and (4.2), we have that $Y = y/\epsilon \in \widetilde{\Omega}_\mu$, for any $\epsilon \in (0,1]$. Lastly, we note that $dy = \epsilon dY$, and as such we have

$$\left\| (\mathbb{R}\mathrm{e}\, Y)^j \partial_x^i \partial_Y^j \Omega^{\mathrm{P}}_\xi(Y) \right\|_{\mathcal{L}^1_\mu} = \epsilon \sup_{0 \le \theta < \mu} \left\| (\mathbb{R}\mathrm{e}\, Y)^j \partial_x^i \partial_Y^j \Omega^{\mathrm{P}}_\xi(Y) \right\|_{L^1_Y(\Gamma_{\epsilon,\theta})} , \tag{6.22}$$

where $\Gamma_{\epsilon,\theta} = \{Y \in \mathbb{C} \colon \epsilon Y \in \partial\Omega_\theta\}$ consists of the union of the two complex paths $\Gamma^\pm_{\epsilon,\theta}$, where

$$\begin{aligned} \Gamma^\pm_{\epsilon,\theta} &= \{Y \in \widetilde{\Omega}_\mu \colon 0 \le \mathbb{R}\mathrm{e}\, Y \le 1/\epsilon, \mathbb{I}\mathrm{m}\, Y = \pm\theta \mathbb{R}\mathrm{e}\, Y\} \\ &\qquad \cup \{Y \in \widetilde{\Omega}_\mu \colon 1/\epsilon \le \mathbb{R}\mathrm{e}\, Y \le (1+\theta)/\epsilon, \mathbb{I}\mathrm{m}\, Y = \pm\theta/\epsilon \mp (\mathbb{R}\mathrm{e}\, Y - 1/\epsilon)\} \, . \end{aligned}$$

Note that for every $Y \in \Gamma^\pm_{\epsilon,\theta}$, we have that $|\mathbb{I}\mathrm{m}\, Y| \le \theta \mathbb{R}\mathrm{e}\, Y \le \mu \mathbb{R}\mathrm{e}\, Y$, independently of $\epsilon$, and for all $\theta \in [0,\mu)$. Due to this fact, using the Taylor expansion argument used to prove (6.13), we have that

$$\sup_{0 \le \theta < \mu} \left\| (\mathbb{R}\mathrm{e}\, Y)^j \partial_Y^j \Omega^{\mathrm{P}}_\xi(Y) \right\|_{L^1_Y(\Gamma_{\epsilon,\theta})} \lesssim \sum_{m \ge 0} \frac{\mu^m}{m!} \left\| Y^{j+m} \partial_Y^{m+j} \Omega^{\mathrm{P}}_\xi(Y) \right\|_{L^1_Y([0,\infty))} . \tag{6.23}$$

Using that $(1+Y)^{-1} \in L^2_Y$, we combine (6.22)–(6.23), and as in (6.14) we expand $e^{\lambda(1+\mu)|\xi|}$ into its power series, to arrive at

$$\sum_\xi e^{\lambda(1+\mu)|\xi|} \left\| Y^j \partial_x^i \partial_Y^j \Omega^{\mathrm{P}}_\xi(Y) \right\|_{\mathcal{L}^1_\mu} \lesssim \epsilon \sum_\xi \sum_{m,n \ge 0} \frac{\mu^m (\lambda(1+\mu))^n}{m! n!} \left\| (1+Y) Y^{j+m} \partial_Y^{m+j} (\partial_x^{n+i} \Omega^{\mathrm{P}})_\xi(Y) \right\|_{L^2_Y([0,\infty))} .$$

Since $\mu \le 1$, $(m+n)!/(m!n!) \le 2^{m+n}$, and as noted at the beginning of the proof of Lemma 6.1 by monotonicity in $\lambda$ and $\mu$ it suffices to consider $\lambda = \lambda_2$ and $\mu = \kappa\lambda_2$, where $\lambda_2 \le \tau(t)/16$ for all $t \in [0, T_2]$, it follows from the above bound that

$$\begin{aligned} \sum_\xi e^{\lambda(1+\mu)|\xi|} \left\| Y^j \partial_x^i \partial_Y^j \Omega^{\mathrm{P}}_\xi(Y) \right\|_{\mathcal{L}^1_\mu} &\lesssim \epsilon \sum_\xi \sum_{m,n \ge 0} \frac{(2\mu)^m (4\lambda)^n}{(m+n)!} \left\| (1+Y) Y^{j+m} \partial_Y^{m+j} (\partial_x^{n+i} \Omega^{\mathrm{P}})_\xi(Y) \right\|_{L^2_Y([0,\infty))} \\ &\lesssim \epsilon \sum_\xi \sum_{m,n \ge 0} \frac{\kappa^m (\tau/4)^{m+n}}{(m+n)!} \left\| (1+Y) Y^{j+m} \partial_Y^{m+j} (\partial_x^{n+i} \Omega^{\mathrm{P}})_\xi(Y) \right\|_{L^2_Y([0,\infty))} . \end{aligned}$$

The $\ell^1(\mathbb{Z} \times \mathbb{N}^2)$ norm taken above in $(\xi, n, m)$ may be converted into an $\ell^2(\mathbb{Z} \times \mathbb{N}^2)$ norm with respect to $(\xi, n, m)$, as in the transition from (6.14) to (6.15) earlier in the proof, at a cost of a factor of $(1+\xi^2)2^{m+n}$. After applying

Plancherel, recalling the definition of the $A_\tau$ norm in (6.2) and the fact that $\gamma \geq 1$, as in (6.16)–(6.17) we obtain

$$\begin{aligned}
&\sum_\xi e^{\lambda(1+\mu)|\xi|} \left\| Y^j \partial_x^i \partial_Y^j \Omega^{\mathrm{P}}_\xi(Y) \right\|_{\mathcal{L}^1_\mu} \\
&\lesssim \epsilon \left( \sum_\xi \sum_{m,n\geq 0} \frac{\kappa^{2m}(\tau/2)^{2(m+n)}}{(m+n)!^2} (1+\xi^2)^2 \left\| (1+Y) Y^{j+m} \partial_Y^{m+j} (\partial_x^{n+i}\Omega^{\mathrm{P}})_\xi(Y) \right\|^2_{L^2_Y([0,\infty))} \right)^{\frac12} \\
&\lesssim \epsilon \left( \sum_{m,n\geq 0} \frac{\kappa^{2m}(\tau/2)^{2(m+n)}}{(m+n)!^2} \left\| (1+Y) Y^{j+m} \partial_Y^{m+j} (1-\partial_x^2) \partial_x^{n+i} \Omega^{\mathrm{P}}(Y) \right\|^2_{L^2} \right)^{\frac12} \\
&\lesssim \epsilon \left\| \Omega^{\mathrm{P}} \right\|_{A_\tau} .
\end{aligned}$$

The desired estimate, (6.21), now follows from (6.7), concluding the proof of the Lemma. □

## 7. The $Y(t)$ norm estimate

We assume that the initial error vorticity obeys a bound consistent with the definitions of the $Y(t)$ and $Z$ norms in (4.9) and (4.8). More precisely, we assume that there exist $\epsilon$-independent constants $\lambda_3, \mu_3 \in (0,1]$ such that

$$\sum_{i+j\leq 2} \|\partial_x^i (y\partial_y)^j \omega_{\mathrm{e}0}\|_{Y_{\lambda_3,\mu_3}} + \sum_{i+j\leq 3} \left\| y \partial_x^i \partial_y^j \omega_{\mathrm{e}0} \right\|_{L^2(y\geq 1/4)} \lesssim 1 \,. \tag{7.1}$$

The goal of this section is to obtain an estimate for the $Y(t)$ norm of $\omega_{\mathrm{e}}$, by appealing to the assumption in the first sum in (7.1). The $Z$ norm estimate is performed in Section 8, cf. Proposition 8.4, and uses the finiteness of the second sum in (7.1).

REMARK 7.1 (**Example of compatible initial condition for the error vorticity**). The assumption (7.1) is for instance satisfied by $\omega_{\mathrm{e}0}$ as defined in (2.17), whenever there exists $\mu_3 \in (0,1]$ such that the function $\psi(y)$ satisfies

$$\sum_{0\leq j\leq 2} \|(y\partial_y)^j \psi\|_{\mathcal{L}^1_{\mu_3}} + \|(y\partial_y)^j \partial_y^2 \psi\|_{\mathcal{L}^1_{\mu_3}} + \sum_{0\leq j\leq 5} \left\| y \partial_y^j \psi \right\|_{L^2(y\geq 1/4)} \lesssim 1 \,, \tag{7.2}$$

where we recall that $\mathcal{L}^1_\mu$ is defined in (4.3) above. In order to see that (2.17) and (7.2) imply (7.1), we note that by the definition (4.4) and the previously established estimate (5.39), we have that for every $i,j \in \{0,1,2\}$,

$$\|\partial_x^i (y\partial_y)^j \omega_{\mathrm{e}0}\|_{Y_{\lambda_3,\mu_3}} \lesssim \sum_{\xi\in\mathbb{Z}} e^{\lambda_3(1+\mu_3)|\xi|} |\xi|^i (1+|\xi|^2) |U^{\mathrm{E}}_{0,\xi}| \lesssim 1$$

as soon as $\lambda_3(1+\mu_3) < \lambda_1$, where $\lambda_1$ is as in (5.7). The later condition is ensured by $\lambda_3 \leq \lambda_1/4$, since $\mu_3 \in (0,1]$. Similarly, the finiteness of the second sum in (7.2) and the estimate (5.39) gives that

$$\left\| \partial_x^i \partial_y^j \omega_{\mathrm{e}0} \right\|_{L^2(y\geq 1/4)} \lesssim \left\| \partial_x^i (1-\partial_x^2) U^{\mathrm{E}}_0 \right\|_{L^2} \lesssim 1 \,,$$

for every $0 \leq i+j \leq 3$. Thus, we have shown that (7.1) holds with $\mu_3$ as in (7.2), and with $\lambda_3 = \lambda_1/4$, $\lambda_1$ as in (5.7).

REMARK 7.2 (**The starred parameters**). Using the parameters $(T_1, \lambda_1, \mu_1)$ from Lemma 5.7, the parameters $(T_2, \lambda_2, \mu_2)$ from Lemma 6.1, and the parameters $(\mu_3, \lambda_3)$ from assumption (7.1), we define the parameters alluded to at the beginning of Section 4 by

$$\mu_* = \min\{\mu_1, \mu_2, \mu_3\}, \qquad \lambda_* = \min\{\lambda_1, \lambda_2, \lambda_3\}, \qquad T_* = \min\left\{ T_1, T_2, \frac{\mu_*}{2\gamma_*} \right\}, \tag{7.3}$$

where $\gamma_* \geq 2$ is the only free parameter left. We emphasize that the implicit constants in $\lesssim$ symbols are not allowed to depend on $\gamma_*$ or on $\epsilon$, but they are allowed to depend on $\mu_*, \lambda_* \in (0,1]$.

Having defined the parameters $\lambda_*, \mu_*$, and with $\gamma_*$ free, the norm $Y(t)$ in (4.9) is well-defined. The main result of this section is as follows.

PROPOSITION 7.1 (**The $Y(t)$ estimate**). *Assume that $\omega_{\rm e0}$ satisfies (7.1), that the Euler solution satisfies the conclusion of Lemma 5.7, and that the Prandtl solution satisfies the conclusions of Lemma 6.1. Let $\gamma_* \geq 2$ be arbitrary, and let $\mu_*, \lambda_*, T_*$ be as defined in (7.3). Then, for all $t \in [0, T_*]$ such that $\sup_{0\leq s\leq t} |||\omega_{\rm e}|||_s$ is finite, we have*

$$\|\omega_{\rm e}(t)\|_{Y(t)} \lesssim 1 + \frac{1}{\gamma_*}\left(\sup_{0\leq s\leq t} |||\omega_{\rm e}|||_s + \epsilon \sup_{0\leq s\leq t} |||\omega_{\rm e}|||_s^2\right), \tag{7.4}$$

*where the implicit constant is independent of $\gamma_*$ and $\epsilon$.*

The remainder of this section is dedicated to the proof of the above proposition, which is concluded in Section 7.4.

**7.1. Analytic estimates for the Stokes equation in the vorticity form.** The $Y_{\lambda,\mu}$ norm estimates for the error vorticity, necessary in order to prove Proposition 7.1, are obtained by using that $\omega_{\rm e}$ solves the Stokes equation (3.17)–(3.18). Applying the Fourier transform in the $x$ variable this Stokes system becomes

$$\partial_t \omega_{{\rm e},\xi} - \epsilon^2 \Delta_\xi \omega_{{\rm e},\xi} = F_\xi \qquad \text{in } \mathbb{H}\,, \tag{7.5}$$

$$\epsilon^2(\partial_y + |\xi|)\omega_{{\rm e},\xi} = B_\xi \qquad \text{on } \partial\mathbb{H}\,, \tag{7.6}$$

for $\xi \in \mathbb{Z}$, where $F_\xi$ denotes the tangential Fourier transform of the forcing term $F$ defined in (3.19) and $B_\xi$ denotes the tangential Fourier transform of the cumulative term appearing on the right side of (3.18), or alternatively, (3.23). The solution of (7.5)–(7.6) is given in terms of the Green's function $G_\xi(t, y, z)$ for this system as

$$\omega_{{\rm e},\xi}(t, y) = \int_0^t \int_0^\infty G_\xi(t-s, y, z) F_\xi(s, z)\, dz ds + \int_0^t G_\xi(t-s, y, 0) B_\xi(s)\, ds + \int_0^\infty G_\xi(t, y, z)\omega_{0{\rm e},\xi}(z)\, dz\,. \tag{7.7}$$

In turn, bounds on the Green's function $G_\xi$ are given in [**53**], and we recall these estimates here.

LEMMA 7.2. *The Green's function $G_\xi$ may be written as*

$$G_\xi = \widetilde{H}_\xi + R_\xi\,,$$

*where*

$$\widetilde{H}_\xi(t, y, z) = \frac{1}{\sqrt{\epsilon^2 t}}\left(e^{-\frac{(y-z)^2}{4\epsilon^2 t}} + e^{-\frac{(y+z)^2}{4\epsilon^2 t}}\right) e^{-\epsilon^2\xi^2 t}\,,$$

*and $R_\xi$ is a function of $y+z$, which obeys the bounds*

$$|\partial_z^k R_\xi(t, y, z)| \lesssim b^{k+1} e^{-\theta_0 b(y+z)} + \frac{1}{(\epsilon^2 t)^{(k+1)/2}} e^{-\theta_0\frac{(y+z)^2}{\epsilon^2 t}} e^{-\frac{\epsilon^2\xi^2 t}{8}}, \qquad k \in \mathbb{N}_0\,, \tag{7.8}$$

*where $\theta_0 > 0$ and*

$$b = b(\xi, \epsilon) = |\xi| + \frac{1}{\epsilon}\,.$$

*The implicit constant in (7.8) depends only on $k$ and $\theta_0$.*

Using the bounds stated in Lemma 7.2 and recalling the definition of $B_\xi$ in (3.18), we obtain the following $Y_{\lambda,\mu}$ analytic estimate for the error vorticity $\omega_{\rm e}$, as defined in (7.7).

LEMMA 7.3 (**The abstract $Y_{\lambda,\mu}$ analytic bound**). *Let $\gamma_* \geq 2$, and fix parameters $\lambda_*, \mu_*, T_* \in (0, 1]$ as in (7.3). Fix times $s, t$ such that $0 \leq s \leq t \leq T_*$, $\lambda \in (0, \lambda_*]$ arbitrary, a parameter $\mu \in (0, \mu_* - \gamma_* s)$, and let*

$$\bar\mu = \mu + \frac{1}{4}(\mu_* - \gamma_* s - \mu)\,, \tag{7.9}$$

*which obeys $\mu < \bar\mu < \mu_* - \gamma_* s$. Then, the forcing (first) term in (7.7) is bounded as*

$$\begin{aligned}
(\mu_* - \gamma_* s - \mu) &\sum_{i+j=2} \left\| \partial_x^i (y\partial_y)^j \int_0^\infty G(t-s, y, z) F(s, z)\, dz \right\|_{Y_{\lambda,\mu}} \\
&+ \sum_{i+j\leq 1} \left\| \partial_x^i (y\partial_y)^j \int_0^\infty G(t-s, y, z) F(s, z)\, dz \right\|_{Y_{\lambda,\bar\mu}} \\
&\lesssim \sum_{i+j\leq 1} \|\partial_x^i (y\partial_y)^j F(s)\|_{Y_{\lambda,\bar\mu}} + \sum_{i+j\leq 1} \|\partial_x^i \partial_y^j F(s)\|_{S_{\bar\mu}}\,.
\end{aligned} \tag{7.10}$$

*The trace kernel (second) term in (7.7) is estimated as*

$$
\begin{aligned}
&(\mu_* - \gamma_* s - \mu) \sum_{i+j=2} \|\partial_x^i (y\partial_y)^j G(t-s,y,0)B(s)\|_{Y_{\lambda,\mu}} + \sum_{i+j\le 1} \left\|\partial_x^i (y\partial_y)^j G(t-s,y,0)B(s)\right\|_{Y_{\lambda,\bar\mu}} \\
&\quad\lesssim \sum_{i\le 1} \left(\|\partial_x^i F(s)\|_{Y_{\lambda,\bar\mu}} + \|\partial_x^i F(s)\|_{S_\mu}\right) + \sum_{i\le 1}\sum_{\xi\in\mathbb{Z}} e^{\bar\mu|\xi|}|\xi|^i \left|\partial_t g_\xi(s)\right| . \qquad (7.11)
\end{aligned}
$$

*Lastly, for the initial datum (third) term in (7.7) we have*

$$
\begin{aligned}
&\sum_{i+j\le 2} \left\|\partial_x^i (y\partial_y)^j \int_0^\infty G(t,y,z)\omega_{0\mathrm{e}}(z)\,dz\right\|_{Y_{\lambda,\mu}} \\
&\quad\lesssim \sum_{i+j\le 2} \|\partial_x^i (y\partial_y)^j \omega_{0\mathrm{e}}\|_{Y_{\lambda,\mu}} + \sum_{i+j\le 2}\sum_{\xi} \left\|\xi^i \partial_y^j \omega_{0\mathrm{e},\xi}\right\|_{L^1(y\ge 1+\mu)} \lesssim 1 . \qquad (7.12)
\end{aligned}
$$

*We note that the second inequality in (7.12) is a direct consequence of the assumption (7.1) and the definition (7.3).*

In view of the integral representation (7.7) and the estimates in Lemma 7.3, it remains to bound the analytic and Sobolev norms of the forcing term $F$, which appears in both (7.10) and in (7.11), the analytic in $x$ norm of the trace term due to $\widetilde{u}^{\mathrm{P}}$ appearing on the right side of (7.11), and the analytic and Sobolev norms of the initial datum in (7.12). This is achieved in Lemma 7.4 below.

**7.2. Contribution of the forcing term.** In view of the representation formula for $\omega_{\mathrm{e}}$ given by (7.7), and of the abstract $Y_{\lambda,\mu}$ norm estimate provided by Lemma 7.3 for the three terms appearing on the right side of (7.7), in order to prove Theorem 3.1 we need to estimate the terms on the right side of (7.10)–(7.12) in terms of the $Y_{\lambda,\mu}$ norm of $\omega_{\mathrm{e}}$. This is the content of the following lemma.

LEMMA 7.4 (**Forcing and trace in $Y_{\lambda,\mu}$ analytic norms**). *Let $s \in [0,T_*]$, $\mu \in (0,\mu_* - \gamma_* s)$, and $\lambda \in (0,\lambda_*]$. For the forcing term in (3.19), we have the pointwise in time estimates*

$$
\begin{aligned}
\sum_{i+j\le 1} \|\partial_x^i (y\partial_y)^j F\|_{Y_{\lambda,\mu}} &\lesssim 1 + \epsilon \|(1+Y)^\gamma Y^j \partial_x^{i+2}\partial_Y^j \Omega^{\mathrm{P}}\|_{P_{\lambda,\mu,\infty}} + \sum_{i+j\le 2} \|\partial_x^i (y\partial_y)^j \omega_{\mathrm{e}}\|_{Y_{\lambda,\mu}\cap S_\mu} \\
&\quad + \epsilon \sum_{i+j\le 1} \|\partial_x^i (y\partial_y)^j \omega_{\mathrm{e}}\|_{Y_{\lambda,\mu}\cap S_\mu} \sum_{i+j\le 2} \left\|\partial_x^i (y\partial_y)^j \omega_{\mathrm{e}}\right\|_{Y_{\lambda,\mu}} \\
&\quad + \epsilon \left\|\partial_x^2 \omega_{\mathrm{e}}\right\|_{Y_{\lambda,\mu}\cap S_\mu} \|y\partial_y \omega_{\mathrm{e}}\|_{Y_{\lambda,\mu}} \qquad (7.13)
\end{aligned}
$$

*and*

$$
\begin{aligned}
\sum_{i+j\le 1} \|\partial_x^i \partial_y^j F\|_{S_\mu} &\lesssim 1 + \sum_{i+j\le 2} \left(\|\partial_x^i\partial_y^j u_{\mathrm{e}}\|_{L^\infty_{x,y}(y\ge 1+\mu)} + \|\partial_x^i\partial_y^j v_{\mathrm{e}}\|_{L^\infty_{x,y}(y\ge 1+\mu)}\right) + \sum_{i+j\le 2} \|\partial_x^i\partial_y^j \omega_{\mathrm{e}}\|_{S_\mu} \\
&\quad + \epsilon \left(\sum_{i+j\le 2} \left(\|\partial_x^i\partial_y^j u_{\mathrm{e}}\|_{L^\infty_{x,y}(y\ge 1+\mu)} + \|\partial_x^i\partial_y^j v_{\mathrm{e}}\|_{L^\infty_{x,y}(y\ge 1+\mu)}\right)\right) \sum_{i+j\le 2} \|\partial_x^i\partial_y^j \omega_{\mathrm{e}}\|_{S_\mu} , \qquad (7.14)
\end{aligned}
$$

*for all $s \in [0,T_*]$. Moreover, for $i \le 1$ we estimate the contribution of $\partial_t g$ appearing in (7.11) as*

$$
\left\|\sum_{\xi\in\mathbb{Z}} e^{\bar\mu|\xi|}|\xi|^i \left|\partial_t g_\xi\right|\right\|_{L^4(0,T_*)} \lesssim 1 , \qquad (7.15)
$$

*for all $t \in [0,T_0]$, with $T_0 \le 1$.*

Before proving the above lemma, we note that Lemma 7.4 immediately implies the following statement.

PROPOSITION 7.5. *Let $s \in [0,T_*]$, $\mu \in (\mu_* - \gamma_* s)$, and $\lambda \in [0,\lambda_*]$. The forcing term $F$ defined in (3.19) satisfies the pointwise estimates*

$$
\sum_{i+j\le 1} \|\partial_x^i (y\partial_y)^j F(s)\|_{Y_{\lambda,\mu}} \lesssim 1 + \|\!|\omega_{\mathrm{e}}|\!\|_s + \mathcal{E}(s) + (1+\epsilon\|\!|\omega_{\mathrm{e}}|\!\|_s)\frac{\|\!|\omega_{\mathrm{e}}|\!\|_s}{(\mu_* - \mu - \gamma_* s)^{1/3}} , \qquad (7.16)
$$

*where*

$$
\int_0^{T_*} (\mathcal{E}(s))^4 \, ds \lesssim \epsilon , \qquad (7.17)
$$

*and*

$$\sum_{i+j\le 1} \|\partial_x^i \partial_y^j F(s)\|_{S_\mu} \lesssim 1 + \|\!|\omega_{\rm e}|\!\|_s + \epsilon \|\!|\omega_{\rm e}|\!\|_s^2 \,. \tag{7.18}$$

*As stated in Remark 7.2, the implicit constants in the* $\lesssim$ *symbols do depend on* $\mu_*, \lambda_* \in (0,1]$*, but they are independent of* $\gamma_* \ge 2$*, and on* $\epsilon \in (0,1]$.

PROOF OF PROPOSITION 7.5. The bound (7.18) follows from (7.14) by appealing to the elliptic estimate (4.15), and noting that due to the inequality mentioned below (4.8) we have

$$\sum_{i+j\le 2} \|\partial_x^i \partial_y^j \omega_{\rm e}(s)\|_{S_\mu} \le \sum_{i+j\le 3} \|\partial_x^i \partial_y^j \omega_{\rm e}\|_{S} = \|\omega_{\rm e}\|_Z \le \|\!|\omega_{\rm e}|\!\|_s \,.$$

Similarly, the bound (7.16) follows from the estimate (7.13), the definition (4.9), which implies

$$\sum_{i+j=2} \|\partial_x^i (y\partial_y)^j \omega_{\rm e}(s)\|_{Y_{\lambda,\mu}} \le \sum_{i+j=2} \|\partial_x^i (y\partial_y)^j \omega_{\rm e}(s)\|_{Y_{\lambda_*,\mu}} \le \frac{\|\!|\omega_{\rm e}(s)|\!\|_s}{(\mu_* - \mu - \gamma_* s)^{1/3}} \,,$$

and the fact that by (6.11) the second term in (7.13), which defines $\mathcal{E}(s)$, may indeed be bounded as in (7.17). □

**7.3. The proof of Lemma 7.4.** The proof of this lemma is structured as follows. First, we establish the stand-alone estimate (7.15). Next, recalling the definition of the forcing term $F$ in (3.19), we estimate the contribution arising from the forcing terms $f_1$ and $f_2$ present in (3.9)–(3.10), as this term does not involve $(u_{\rm e}, v_{\rm e}, \omega_{\rm e})$. The next subsection provides analytic and Sobolev estimates for the error velocity $(u_{\rm e}, v_{\rm e})$ in terms of the error vorticity $\omega_{\rm e}$, via estimates for the inhomogeneous div-curl system (7.33). We conclude by estimating the remaining terms in (3.19).

7.3.1. *The proof of the estimate* (7.15)*.* In order to establish the bound (7.15), we prove the pointwise in time estimate

$$\sum_{\xi\in\mathbb{Z}} e^{\bar\mu|\xi|} |\xi|^i \, |\partial_t g_\xi(s)| \lesssim 1 + \|\partial_x^{i+1} \Omega^{\rm P}(s)\|_{Y_{\lambda_2,\mu_2,\infty}} \,, \tag{7.19}$$

where $\lambda_2$ and $\mu_2$ are as defined in Lemma 6.1. This estimate may then be combined with $L^4$ in time bound (6.11) with $j=0$, $\lambda=\lambda_2$, and $\mu=\mu_2$, to imply (7.15).

In order to prove (7.19), we first compute $\partial_t g$. Recall that $g = -\int_0^\infty \partial_x \widetilde{u}^{\rm P} dY$, and that in (3.22) we have computed a formula for $\int_0^\infty \partial_t \widetilde{u}^{\rm P} dY$. Combining these two identities, we arrive at

$$\begin{aligned} |\xi|^i \, |\partial_t g_\xi| &\lesssim \left| \left( \int_0^\infty \partial_x^{i+1} \partial_t \widetilde{u}^{\rm P} dY \right)_\xi \right| \\ &\lesssim |\xi|^{i+1} \left|\Omega_\xi^{\rm P}|_{Y=0}\right| + |\xi|^{i+1} |(U^{\rm E} g)_\xi| \\ &\quad + |\xi|^{i+1} \left| \left( \partial_x U^{\rm E} \int_0^\infty \widetilde{u}^{\rm P} dY \right)_\xi \right| + |\xi|^{i+2} \left| \left( \int_0^\infty (\widetilde{u}^{\rm P})^2 dY \right)_\xi \right| . \end{aligned} \tag{7.20}$$

Using $|\Omega_\xi^{\rm P}|_{Y=0}| \lesssim \|\Omega_\xi^{\rm P}\|_{L^\infty(\widetilde{\Omega}_{\bar\mu})}$ and the parameter inequality

$$\bar\mu < \mu_* \le \mu_2 = \kappa\lambda_2 \le \frac{\lambda_2(1+\mu_2)}{8} \,, \tag{7.21}$$

which holds by the definition (7.3), the parameter definitions in Lemma 4.1, and the choice $\kappa = 1/8$, we bound the contribution of the first term in (7.20) as

$$\sum_{\xi\in\mathbb{Z}} e^{\bar\mu|\xi|} |\xi|^{i+1} \left|\Omega_\xi^{\rm P}|_{Y=0}\right| \lesssim \sum_{\xi\in\mathbb{Z}} e^{\frac18 \lambda_2(1+\mu_2)|\xi|} |\xi|^{i+1} \|\Omega_\xi^{\rm P}\|_{L^\infty(\widetilde{\Omega}_{\bar\mu})} \lesssim \|\partial_x^{i+1}\Omega^{\rm P}\|_{P_{\lambda_2,\mu_2,\infty}} \,,$$

an expression which belongs to $L^4(0,T_*)$ according to (6.11), with the norm of constant size. For the second term in (7.20) we use that the Fourier transform of a product is a (discrete) convolution, which is well-estimated using $\ell^1_\xi$ norms. Therefore, by also appealing to the definition of $g$ in (3.14), to the bounds (5.39) and (6.9), and the parameter estimates (7.21) and $\bar\mu \le \lambda_2/8 \le \lambda_1/8$, we arrive at

$$\sum_{\xi\in\mathbb{Z}} e^{\bar\mu|\xi|} |\xi|^{i+1} \left|(U^{\rm E} g)_\xi\right| \lesssim \left( \sum_{\xi\in\mathbb{Z}} e^{\bar\mu|\xi|} (|\xi|+1)^{i+1} \left|U_\xi^{\rm E}\right| \right) \left( \sum_{\xi\in\mathbb{Z}} e^{\bar\mu|\xi|} (|\xi|+1)^{i+1} |g_\xi| \right) \lesssim 1 \,.$$

For the third term in (7.20), using the same parameter inequalities and appealing to (5.39) and (6.9) we similarly have

$$\sum_{\xi\in\mathbb{Z}} e^{\bar\mu|\xi|}|\xi|^{i+1}\left|\left(\partial_x U^{\rm E}\int_0^\infty \widetilde u^{\rm P} dY\right)_\xi\right| \lesssim \left(\sum_{\xi\in\mathbb{Z}} e^{\bar\mu|\xi|}(|\xi|+1)^{i+2}\left|U^{\rm E}_\xi\right|\right)\left(\sum_{\xi\in\mathbb{Z}} e^{\bar\mu|\xi|}(|\xi|+1)^{i+1}\int_0^\infty \left|\widetilde u^{\rm P}_\xi\right| dY\right)\lesssim 1\,.$$

The bound for the last term in (7.20) is similar, but also uses the estimate for the first term on the left side of (6.9):

$$\sum_{\xi\in\mathbb{Z}} e^{\bar\mu|\xi|}|\xi|^{i+2}\left|\left(\int_0^\infty (\widetilde u^{\rm P})^2 dY\right)_\xi\right| \lesssim \left(\sum_{\xi\in\mathbb{Z}} e^{\bar\mu|\xi|}\left|\widetilde u^{\rm P}_\xi\right|_{L^\infty_Y}\right)\left(\sum_{\xi\in\mathbb{Z}} e^{\bar\mu|\xi|}|\xi|^{i+2}\left|\int_0^\infty \widetilde u^{\rm P}_\xi dY\right|\right)\lesssim 1\,.$$

This concludes the proof of (7.19).

7.3.2. *Size of $\partial_x f_2 - \partial_y f_1$ in analytic and Sobolev norms.* According to (3.19), the last term in the definition of $F$ is the forcing term $\partial_x f_2 - \partial_y f_1$. In this section we provide a $Y_{\lambda,\mu}$ estimate for this term, which is needed in proving (7.13), and a $S_\mu$ estimate, which is required to prove (7.14).

LEMMA 7.6. *Let $0<\mu\le\mu_*$ and $\lambda\le\lambda_*$ be arbitrary. Then, for integers $i,j\ge 0$ such that $i+j\le 1$ we obtain*

$$\|\partial_x^i (y\partial_y)^j(\partial_x f_2 - \partial_y f_1)\|_{Y_{\lambda,\mu}} \lesssim 1 + \epsilon\|(1+Y)^{3/2}Y^j\partial_x^{i+2}\partial_Y^j\Omega^{\rm P}\|_{P_{\lambda,\mu,\infty}} \tag{7.22}$$

*and*

$$\|\partial_x^i (y\partial_y)^j(\partial_x f_2 - \partial_y f_1)\|_{S_\mu} \lesssim 1\,. \tag{7.23}$$

*By the estimate (6.11), we have that the second term on the right side of (7.22) is $\mathcal{O}(\epsilon)$ when measured in $L^4([0,T_*])$.*

PROOF OF LEMMA 7.6. We only consider the estimate (7.22) in the case $i=j=0$. The case $i+j=1$ follows mutatis mutandis. According to the definitions of $f_1$ and $f_2$ in (3.15)–(3.16), after taking into account incompressibility, the definitions of $\widetilde u^{\rm P}$ and $\bar v^{\rm P}$, and a number of cancellations, we have

$$\begin{aligned}\partial_x f_2 - \partial_y f_1 &= \frac{1}{\epsilon^2}\partial_x\Omega^{\rm P}(u^{\rm E}-U^{\rm E}) + \frac{1}{\epsilon^3}\partial_Y\Omega^{\rm P}(v^{\rm E}+y\partial_x U^{\rm E}) - \frac{1}{\epsilon}\widetilde u^{\rm P}\partial_x\omega^{\rm E} + \frac{1}{\epsilon}v^{\rm E}\partial_{xx}\widetilde u^{\rm P}\\ &\quad + \bar v^{\rm P}\partial_{yy}u^{\rm E} - 2\partial_{xx}\Omega^{\rm P} + \epsilon^2\partial_x^3\bar v^{\rm P} + \epsilon\Delta\omega^{\rm E} - \partial_t\partial_x\bar v^{\rm P} - (\widetilde u^{\rm P}+u^{\rm E})\partial_{xx}\bar v^{\rm P} + \bar v^{\rm P}\partial_{xx}(\widetilde u^{\rm P}+u^{\rm E})\,.\end{aligned} \tag{7.24}$$

Noting that $-\partial_t\partial_x\bar v^{\rm P} = -\partial_{xx}(\int_Y^\infty \partial_t\widetilde u^{\rm P} dY')$ by (3.21) and using the Prandtl evolution (3.20), we obtain

$$\begin{aligned}-\partial_t\partial_x\bar v^{\rm P} &= -\partial_{xx}\left(-\Omega^{\rm P} + v^P\widetilde u^{\rm P} - U^{\rm E}\bar v^{\rm P} - \partial_x\int_Y^\infty (\widetilde u^{\rm P})^2 dY' - 2\partial_x U^{\rm E}\int_Y^\infty \widetilde u^{\rm P} dY'\right)\\ &= \partial_{xx}\Omega^{\rm P} - \partial_{xx}(v^P\widetilde u^{\rm P}) + U^{\rm E}\partial_{xx}\bar v^{\rm P} + 5\partial_{xx}U^{\rm E}\bar v^{\rm P} + 4\partial_x U^{\rm E}\partial_x\bar v^{\rm P}\\ &\quad + \partial_x^3\int_Y^\infty (\widetilde u^{\rm P})^2 dY' + 2\partial_x^3 U^{\rm E}\int_Y^\infty \widetilde u^{\rm P} dY'\,.\end{aligned}$$

Combining the above two identities allows us to rewrite

$$\begin{aligned}-\partial_y f_1 + \partial_x f_2 &= \frac{1}{\epsilon^2}\partial_x\Omega^{\rm P}(u^{\rm E}-U^{\rm E}) + \frac{1}{\epsilon^3}\partial_Y\Omega^{\rm P}(v^{\rm E}+y\partial_x U^{\rm E}) - \frac{1}{\epsilon}\widetilde u^{\rm P}\partial_x\omega^{\rm E} + \frac{1}{\epsilon}v^{\rm E}\partial_{xx}\widetilde u^{\rm P}\\ &\quad + \epsilon^2\partial_x^3\bar v^{\rm P} + \epsilon\Delta\omega^{\rm E} - \partial_{xx}\Omega^{\rm P} - (\widetilde u^{\rm P}+u^{\rm E}-U^{\rm E})\partial_{xx}\bar v^{\rm P} + 4\partial_x\bar v^{\rm P}\partial_x U^{\rm E}\\ &\quad + \bar v^{\rm P}(\partial_{xx}\widetilde u^{\rm P} + \Delta u^{\rm E} + 5\partial_{xx}U^{\rm E}) - \partial_x^2(\widetilde u^{\rm P}v^{\rm P}) + 2\partial_x^3 U^{\rm E}\int_Y^\infty \widetilde u^{\rm P} dY' + \partial_x^3\int_Y^\infty (\widetilde u^{\rm P})^2 dY'\\ &= f_{\rm e,1} + \cdots + f_{\rm e,13}\,.\end{aligned} \tag{7.25}$$

For the $Y_{\lambda,\mu}$ estimate of $-\partial_y f_1 + \partial_x f_2$ we consider the thirteen terms in (7.25) individually. For the first term in (7.25), we have

$$f_{\rm e,1} = \frac{1}{\epsilon^2}\partial_x\Omega^{\rm P}(u^{\rm E}-U^{\rm E}) = \frac{1}{\epsilon}Y\partial_x\Omega^{\rm P}\frac{u^{\rm E}-U^{\rm E}}{y}\,.$$

Using (5.37), (6.8), and (4.11b) we thus obtain

$$\|f_{\rm e,1}\|_{Y_{\lambda,\mu}} \lesssim \frac{1}{\epsilon}\epsilon\|(1+Y)^{3/2}Y\partial_x\Omega^{\rm P}\|_{P_{\lambda,\mu,\infty}}\left\|\frac{u^{\rm E}-U^{\rm E}}{y}\right\|_{Y_{\lambda,\mu,\infty}} \lesssim 1\,. \tag{7.26}$$

Similarly, we have

$$f_{\rm e,2} = \frac{1}{\epsilon^3}\partial_Y \Omega^{\rm P}(v^{\rm E} + y\partial_x U^{\rm E}) = \frac{1}{\epsilon} Y^2 \partial_Y \Omega^{\rm P} \frac{v^{\rm E} + y\partial_x U^{\rm E}}{y^2}\,,$$

and so by appealing to (5.38), (6.8), and (4.11b) we may estimate

$$\|f_{\rm e,2}\|_{Y_{\lambda,\mu}} \lesssim \frac{1}{\epsilon}\epsilon \|(1+Y)^{3/2} Y^2 \partial_y \Omega^{\rm P}\|_{P_{\lambda,\mu,\infty}} \left\| \frac{v^{\rm E} + y\partial_x U^{\rm E}}{y^2}\right\|_{Y_{\lambda,\mu,\infty}} \lesssim 1\,. \tag{7.27}$$

In a similar fashion, from (5.36), (6.9), and (4.11b) we have

$$\|f_{\rm e,3}\|_{Y_{\lambda,\mu}} \lesssim \frac{1}{\epsilon}\epsilon \|(1+Y)^{3/2} \widetilde{u}^{\rm P}\|_{P_{\lambda,\mu,\infty}} \|\partial_x \omega^{\rm E}\|_{Y_{\lambda,\mu,\infty}} \lesssim 1\,, \tag{7.28}$$

while from (5.38), (6.9), and (4.11b) we have

$$\|f_{\rm e,4}\|_{Y_{\lambda,\mu}} \lesssim \frac{1}{\epsilon}\epsilon \|(1+Y)^{3/2} \partial_{xx}\widetilde{u}^{\rm P}\|_{P_{\lambda,\mu,\infty}} \|v^{\rm E}\|_{Y_{\lambda,\mu,\infty}} \lesssim 1\,. \tag{7.29}$$

This concludes the estimates for all the terms which have inverse powers of $\epsilon$ in (7.25). The next seven terms in (7.25) all have simple bounds in view of the bounds (4.11), (4.12), (4.13), and Lemmas 5.7, 6.1:

$$\begin{aligned}
\|f_{\rm e,5}\|_{Y_{\lambda,\mu}} &\lesssim \epsilon^2 \|\partial_x^3 \bar{v}^{\rm P}\|_{Y_{\lambda,\mu}} \lesssim \epsilon^3 \|(1+Y)^{\frac32} \partial_x^3 \bar{v}^{\rm P}\|_{P_{\lambda,\mu,\infty}} \lesssim \epsilon^3 \\
\|f_{\rm e,6}\|_{Y_{\lambda,\mu}} &\lesssim \epsilon \|\Delta \omega^{\rm E}\|_{Y_{\lambda,\mu}} \lesssim \epsilon \|\Delta \omega^{\rm E}\|_{Y_{\lambda,\mu,\infty}} \lesssim \epsilon \\
\|f_{\rm e,7}\|_{Y_{\lambda,\mu}} &\lesssim \|\partial_{xx}\Omega^{\rm P}\|_{Y_{\lambda,\mu}} \lesssim \epsilon \|(1+Y)^{3/2} \partial_x^2 \Omega^{\rm P}\|_{P_{\lambda,\mu,\infty}} \\
\|f_{\rm e,8}\|_{Y_{\lambda,\mu}} &\lesssim \epsilon \|(1+Y)^{3/2} \partial_{xx}\bar{v}^{\rm P}\|_{P_{\lambda,\mu,\infty}} \left(\|\widetilde{u}^{\rm P}\|_{P_{\lambda,\mu,\infty}} + \|u^{\rm E} - U^{\rm E}\|_{Y_{\lambda,\mu,\infty}}\right) \lesssim \epsilon \\
\|f_{\rm e,9}\|_{Y_{\lambda,\mu}} &\lesssim \epsilon \|(1+Y)^{3/2} \partial_x \bar{v}^{\rm P}\|_{P_{\lambda,\mu,\infty}} \|\partial_x U^{\rm E}\|_{Y_{\lambda,\mu,\infty}} \lesssim \epsilon \\
\|f_{\rm e,10}\|_{Y_{\lambda,\mu}} &\lesssim \epsilon \|(1+Y)^{3/2} \bar{v}^{\rm P}\|_{P_{\lambda,\mu,\infty}} \left(\|\partial_{xx}\widetilde{u}^{\rm P}\|_{P_{\lambda,\mu,\infty}} + \|\Delta u^{\rm E}\|_{Y_{\lambda,\mu,\infty}} + \|\partial_{xx} U^{\rm E}\|_{Y_{\lambda,\mu,\infty}}\right) \lesssim \epsilon \\
\|f_{\rm e,11}\|_{Y_{\lambda,\mu}} &\lesssim \epsilon \|(1+Y)^{3/2} \partial_{xx}(\widetilde{u}^{\rm P} v^{\rm P})\|_{P_{\lambda,\mu,\infty}} \lesssim \epsilon \sum_{i=0}^{2} \|(1+Y)^{5/2} \partial_x^i \widetilde{u}^{\rm P}\|_{P_{\lambda,\mu,\infty}} \left\|Y^{-1} \partial_x^{2-i} v^{\rm P}\right\|_{P_{\lambda,\mu,\infty}} \lesssim \epsilon\,.
\end{aligned} \tag{7.30}$$

We note that the above stated estimate for the term $f_{\rm e,7}$ is responsible for the second term on the right side of (7.22). It remains to consider the last two terms in (7.25). From Lemma 6.1 and using the bound

$$\begin{aligned}
&\sum_{\xi\in\mathbb{Z}} e^{\lambda(1+\mu)|\xi|} \sup_Y \left( (1+Y)^{\frac32} \int_Y^\infty \left|\widetilde{u}^{\rm P}_\xi\right| dY' \right) \\
&\lesssim \sum_{\xi\in\mathbb{Z}} e^{\lambda(1+\mu)|\xi|} \sup_Y \left( (1+Y)^{\frac32} \int_Y^\infty (1+Y')^{\gamma-\frac32} \left|\widetilde{u}^{\rm P}_\xi\right| \frac{dY'}{(1+Y')^{\gamma-\frac32}} \right) \\
&\lesssim \|(1+Y)^{\gamma-\frac32} \widetilde{u}^{\rm P}\|_{P_{\lambda,\mu,\infty}} \sup_Y (1+Y)^{\frac32} \int_Y^\infty \frac{dY'}{(1+Y')^{\gamma-\frac32}} \lesssim 1\,,
\end{aligned}$$

which holds since $\gamma \geq 4$, and combining with estimates (4.11c) and (4.12), we obtain

$$\|f_{\rm e,12}\|_{Y_{\lambda,\mu}} \lesssim \|\partial_x^3 U^{\rm E}\|_{Y_{\lambda,\mu,\infty}} \left\| \int_Y^\infty \widetilde{u}^{\rm P} dY' \right\|_{Y_{\lambda,\mu}} \lesssim \epsilon \left\| (1+Y)^{\frac32} \int_Y^\infty \widetilde{u}^{\rm P} dY' \right\|_{P_{\lambda,\mu,\infty}} \lesssim \epsilon\,. \tag{7.31}$$

From the product rule, estimate (4.12), and Lemma 6.1, we also obtain

$$\|f_{\rm e,13}\|_{Y_{\lambda,\mu}} \lesssim \epsilon \|(1+Y)^{\frac32} \partial_x^2 \bar{v}^{\rm P}\|_{P_{\lambda,\mu,\infty}} \|\widetilde{u}^{\rm P}\|_{P_{\lambda,\mu,\infty}} + \epsilon \|(1+Y)^{\frac32} \partial_x \bar{v}^{\rm P}\|_{P_{\lambda,\mu,\infty}} \|\partial_x \widetilde{u}^{\rm P}\|_{P_{\lambda,\mu,\infty}} \lesssim \epsilon\,. \tag{7.32}$$

Adding the upper bounds in (7.26)–(7.32), completes the proof of the $Y_{\lambda,\mu}$ estimate claimed in (7.22).

In order to complete the proof for the lemma, it remains to estimate the $S_\mu$ norm of $\partial_x^i (y\partial_y)^j(-\partial_y f_1 + \partial_x f_2)$; as noted earlier, we only give these details for the case $i = j = 0$. As before, we separately consider the thirteen terms in (7.25). We note that all terms that are a product of Prandtl part and Euler part are in fact small, in view of the product estimates (4.11d)–(4.11e), and the previously established estimates (5.2)–(5.3) and (5.39) for Euler, respectively (6.8)–(6.11) for Prandtl; however, since we only wish to obtain an $\mathcal{O}(1)$ upper bound, we do not attempt to estimate these terms in terms of optimal powers of $\epsilon$. Using (4.11e), (5.3), (5.39), and (6.8), we have

$$\begin{aligned}\|f_{\mathrm{e},1}\|_{S_\mu} &\le \epsilon^{-1}\left\|Y\partial_x\Omega^{\mathrm{P}}\frac{u^{\mathrm{E}}-U^{\mathrm{E}}}{y}\right\|_{S_\mu} \le \epsilon^{-1}\epsilon\|Y^2\partial_x\Omega^{\mathrm{P}}\|_{P_{\lambda,\mu,\infty}}\left\|\frac{u^{\mathrm{E}}-U^{\mathrm{E}}}{y}\right\|_{H^1_xL^2_y(y\ge 1/2)}\\ &\lesssim \|Y(1+Y)^{\gamma-1}\partial_x\Omega^{\mathrm{P}}\|_{P_{\lambda,\mu,\infty}}\left(\|u^{\mathrm{E}}\|_{H^1_xL^2_y(y\ge1/2)}+\|U^{\mathrm{E}}\|_{H^1_x}\right)\lesssim 1\end{aligned}$$

since $\gamma\ge 4$. For the next three terms, we similarly obtain

$$\begin{aligned}\|f_{\mathrm{e},2}\|_{S_\mu} &\lesssim \epsilon^{-1}\left\|Y^2\partial_Y\Omega^{\mathrm{P}}\frac{v^{\mathrm{E}}+y\partial_xU^{\mathrm{E}}}{y^2}\right\|_{S_\mu}\lesssim\|(1+Y)^2Y^2\partial_Y\Omega^{\mathrm{P}}\|_{P_{\lambda,\mu,\infty}}\left\|\frac{v^{\mathrm{E}}+y\partial_xU^{\mathrm{E}}}{y^2}\right\|_{H^1_xL^2_y(y\ge1/2)}\lesssim 1\\ \|f_{\mathrm{e},3}\|_{S_\mu} &\lesssim \epsilon^{-1}\epsilon\left\|Y\widetilde{u}^{\mathrm{P}}\right\|_{P_{\lambda,\mu,\infty}}\left\|\partial_x\omega^{\mathrm{E}}\right\|_{S_\mu}\lesssim\left\|(1+Y)^{\gamma-\frac32}\widetilde{u}^{\mathrm{P}}\right\|_{P_{\lambda,\mu,\infty}}\left\|y\partial_x\omega^{\mathrm{E}}\right\|_{H^1_xL^2_y(y\ge1/2)}\lesssim 1\\ \|f_{\mathrm{e},4}\|_{S_\mu} &\lesssim \left\|Y\partial_{xx}\widetilde{u}^{\mathrm{P}}\frac{v^{\mathrm{E}}}{y}\right\|_{S_\mu}\lesssim\|Y\partial_{xx}\widetilde{u}^{\mathrm{P}}\|_{P_{\lambda,\mu,\infty}}\left\|\frac{v^{\mathrm{E}}}{y}\right\|_{S_\mu}\lesssim\left\|(1+Y)^{\gamma-\frac32}\partial_{xx}\widetilde{u}^{\mathrm{P}}\right\|_{P_{\lambda,\mu,\infty}}\|v^{\mathrm{E}}\|_{H^1_xL^2_y(y\ge1/2)}\lesssim 1\end{aligned}$$

since $\gamma\ge 4$. For the fifth and seventh terms in the right side of (7.25), which are linear in Prandtl terms, we appeal to (4.11f) with $g\equiv 1\in L^2_xL^\infty_y$, to deduce

$$\|f_{\mathrm{e},5}\|_{S_\mu}+\|f_{\mathrm{e},7}\|_{S_\mu}\lesssim\epsilon^2\left\|\partial_x^3\bar{v}^{\mathrm{P}}\right\|_{P_{\lambda,\mu,\infty}}+\left\|\partial_x^2\Omega^{\mathrm{P}}\right\|_{P_{\lambda,\mu,\infty}}\lesssim 1\,.$$

For the only error term which is linear in the Euler solution, we note that

$$\|f_{\mathrm{e},6}\|_{S_\mu}=\epsilon\left\|\Delta\omega^{\mathrm{E}}\right\|_{S_\mu}\le\epsilon\left\|y\Delta\omega^{\mathrm{E}}\right\|_{L^2(y\ge1+\mu)}+\epsilon\left\|y\partial_x\Delta\omega^{\mathrm{E}}\right\|_{L^2(y\ge1+\mu)}\lesssim\epsilon\lesssim 1$$

in view of (5.2). The remaining terms consist of Euler–Prandtl products, which are estimated using (4.11d)–(4.11f), and Prandtl–Prandtl products, which are bounded using (4.11f) and the fact that

$$\|g(x,Y)\|_{H^1_xL^\infty_y(y\ge1/2)}\le\|g(x,Y)\|_{H^1_xL^\infty_Y(Y\ge1/(2\epsilon))}\lesssim\epsilon^\theta\|Y^\theta g\|_{P_{\lambda,\mu,\infty}}\,,$$

for any $\lambda,\mu>0$ and any $\theta\ge0$. We may thus show that

$$\begin{aligned}\|f_{\mathrm{e},8}\|_{S_\mu}&\lesssim\|\partial_x^2\bar{v}^{\mathrm{P}}\|_{P_{\lambda,\mu,\infty}}\left(\left\|\widetilde{u}^{\mathrm{P}}\right\|_{H^1_xL^\infty_y(y\ge1/2)}+\left\|u^{\mathrm{E}}\right\|_{H^1_xL^2_y(y\ge1/2)}+\left\|U^{\mathrm{E}}\right\|_{H^1_x}\right)\lesssim1\\ \|f_{\mathrm{e},9}\|_{S_\mu}&\lesssim\left\|\partial_x\bar{v}^{\mathrm{P}}\right\|_{P_{\lambda,\mu,\infty}}\left\|\partial_xU^{\mathrm{E}}\right\|_{H^1_x}\lesssim1\\ \|f_{\mathrm{e},10}\|_{S_\mu}&\lesssim\left\|\bar{v}^{\mathrm{P}}\right\|_{P_{\lambda,\mu,\infty}}\left(\left\|\partial_x^2\widetilde{u}^{\mathrm{P}}\right\|_{H^1_xL^\infty_y(y\ge1/2)}+\left\|\Delta u^{\mathrm{E}}\right\|_{H^1_xL^2_y(y\ge1/2)}+\left\|\partial_x^2U^{\mathrm{E}}\right\|_{H^1_x}\right)\lesssim1\\ \|f_{\mathrm{e},11}\|_{S_\mu}&\lesssim\left\|\partial_x^2(\widetilde{u}^{\mathrm{P}}v^{\mathrm{P}})\right\|_{P_{\lambda,\mu,\infty}}\lesssim1\\ \|f_{\mathrm{e},12}\|_{S_\mu}&\lesssim\left\|\int_Y^\infty\widetilde{u}^{\mathrm{P}}\right\|_{P_{\lambda,\mu,\infty}}\left\|\partial_x^3U^{\mathrm{E}}\right\|_{H^1_x}\lesssim\left\|(1+Y)^{\gamma-\frac32}\widetilde{u}^{\mathrm{P}}\right\|_{P_{\lambda,\mu,\infty}}\left\|\partial_x^3U^{\mathrm{E}}\right\|_{H^1_x}\lesssim1\\ \|f_{\mathrm{e},13}\|_{S_\mu}&\lesssim\left\|\partial_x^3\int_Y^\infty(\widetilde{u}^{\mathrm{P}})^2\right\|_{P_{\lambda,\mu,\infty}}\lesssim1\,,\end{aligned}$$

where in the last two inequalities we have used that $(1+Y)^{\frac32-\gamma}\in L^1_Y$, since $\gamma\ge4$. This completes the proof of Lemma 7.6. ☐

7.3.3. *Modified Biot-Savart law.* The first, third, and fourth terms in the definition of $F$ in (3.19) all involve the vector $(u_{\mathrm{e}},v_{\mathrm{e}})$, which is obtained from the error vorticity $\omega_{\mathrm{e}}$ and the Prandtl boundary vertical velocity $g=-\bar{v}^{\mathrm{P}}|_{Y=0}$ (see (3.14)), via the div-curl system

$$\begin{aligned}-\partial_yu_{\mathrm{e}}+\partial_xv_{\mathrm{e}}&=\omega_{\mathrm{e}} &&\text{in}\quad\mathbb{H}\\ \partial_xu_{\mathrm{e}}+\partial_yv_{\mathrm{e}}&=0 &&\text{in}\quad\mathbb{H}\\ v_{\mathrm{e}}=g&=\partial_xh &&\text{on}\quad\partial\mathbb{H}\,.\end{aligned}\tag{7.33}$$

The representation formula for the system (7.33) is as follows. With $\nabla^\perp = (-\partial_y, \partial_x)$, we define the corrector $\nabla^\perp(e^{-y|\partial_x|}h(x))$, which is curl-free, divergence-free, and its second component equals $\partial_x h = g$ on $\partial\mathbb{H}$. Therefore,

$$
\begin{aligned}
-\partial_y\left(u_{\mathrm{e}} + \frac{\partial_x}{|\partial_x|}e^{-|\partial_x|y}g\right) + \partial_x\left(v_{\mathrm{e}} - e^{-|\partial_x|y}g\right) &= \omega_{\mathrm{e}} \qquad \text{in} \qquad \mathbb{H}\\
\partial_x\left(u_{\mathrm{e}} + \frac{\partial_x}{|\partial_x|}e^{-|\partial_x|y}g\right) + \partial_y\left(v_{\mathrm{e}} - e^{-|\partial_x|y}g\right) &= 0 \qquad \text{in} \qquad \mathbb{H}\\
v_{\mathrm{e}} - e^{-|\partial_x|y}g &= 0 \qquad \text{on} \qquad \partial\mathbb{H}\,.
\end{aligned}
\tag{7.34}
$$

Using the classical Biot-Savart law (cf. [47], or (6.2)–(6.3) in [38]), upon taking the Fourier transforms in $x$ we deduce

$$
\begin{aligned}
u_{\mathrm{e},\xi}(y) &= -\frac{\hat{\imath}\xi}{|\xi|}e^{-|\xi|y}g_\xi\\
&\quad + \frac{1}{2}\left(-\int_0^y e^{-|\xi|(y-z)}(1-e^{-2|\xi|z})\omega_{\mathrm{e},\xi}(z)\,dz + \int_y^\infty e^{-|\xi|(z-y)}(1+e^{-2|\xi|y})\omega_{\mathrm{e},\xi}(z)\,dz\right)
\end{aligned}
\tag{7.35}
$$

and

$$
\begin{aligned}
v_{\mathrm{e},\xi}(y) &= e^{-|\xi|y}g_\xi\\
&\quad - \frac{\hat{\imath}\xi}{2|\xi|}\left(\int_0^y e^{-|\xi|(y-z)}(1-e^{-2|\xi|z})\omega_{\mathrm{e},\xi}(z)\,dz + \int_y^\infty e^{-|\xi|(z-y)}(1-e^{-2|\xi|y})\omega_{\mathrm{e},\xi}(z)\,dz\right)\,.
\end{aligned}
\tag{7.36}
$$

As a direct consequence of the above formulae, we obtain an inequality for the velocity in a $L^\infty$-based analytic norm in terms of the vorticity in a $L^1_y$-based analytic norm.

LEMMA 7.7 (**$Y_{\lambda,\mu,\infty}$ norm estimates for the modified Biot-Savart law**). *Let $\mu \in (0, \mu_* - \gamma_* t)$ and $\lambda \in (0, \lambda_*]$. Then, the functions $u_{\mathrm{e}}$ and $v_{\mathrm{e}}$ defined via the modified Biot-Savart law (7.35)–(7.36), satisfy the estimates*

$$
\|\partial_x^i (y\partial_y)^j u_{\mathrm{e}}\|_{Y_{\lambda,\mu,\infty}} \lesssim \|\partial_x^{i+j}\omega_{\mathrm{e}}\|_{Y_{\lambda,\mu}\cap S_\mu} + j\left(\|\omega_{\mathrm{e}}\|_{Y_{\lambda,\mu}} + \|y\partial_y\omega_{\mathrm{e}}\|_{Y_{\lambda,\mu}}\right) + \|\partial_x^i g\|_{P_{\lambda,\mu,\infty}}
\tag{7.37}
$$

*and*

$$
\left\|(y\partial_y)^j\left(\frac{\partial_x^i(v_{\mathrm{e}}-g)}{y}\right)\right\|_{Y_{\lambda,\mu,\infty}} \lesssim \|\partial_x^{i+1}\omega_{\mathrm{e}}\|_{Y_{\lambda,\mu}\cap S_\mu} + \|\partial_x^{i+1} g\|_{P_{\lambda,\mu,\infty}}
$$

*for all integers $i,j \geq 0$ such that $i+j \leq 1$. Lastly, for $0 \leq i \leq 1$ we have*

$$
\|\partial_x^i v_{\mathrm{e}}\|_{Y_{\lambda,\mu,\infty}} \lesssim \|\partial_x^i\omega_{\mathrm{e}}\|_{Y_{\lambda,\mu}\cap S_\mu} + \|\partial_x^i g\|_{P_{\lambda,\mu,\infty}}\,.
\tag{7.38}
$$

PROOF OF LEMMA 7.7. The proof follows closely estimates in Section 6 of [38] and Section 4 of [39]. For simplicity, we only provide estimates for the real values in definition (4.5); the bounds along complex contour integrals follow along the same lines. From (7.35) and (7.36), the velocity field $(u_{\mathrm{e}}, v_{\mathrm{e}})$ can be decomposed as

$$
(u_{\mathrm{e}}, v_{\mathrm{e}}) = \left(-\frac{\hat{\imath}\xi}{|\xi|}e^{-|\xi|y}g_\xi, e^{-|\xi|y}g_\xi\right) + (\widetilde{u}_{\mathrm{e}}, \widetilde{v}_{\mathrm{e}})\,,
\tag{7.39}
$$

where $(\widetilde{u}_{\mathrm{e}}, \widetilde{v}_{\mathrm{e}})$ is obtained from the vorticity $\omega_{\mathrm{e}}$ by the usual Biot-Savart law on $\mathbb{T}\times\mathbb{R}_+$ (cf. (7.34)).

The first term on the right of (7.39) contributes the $g$ terms on the right sides of (7.37)–(7.38) thanks to the inequalities

$$
\left|(y|\xi|)^j e^{\lambda(1+\mu-y)|\xi|}e^{-|\xi|y}\right| \lesssim e^{\lambda(1+\mu)|\xi|}, \qquad \left|(y\partial_y)^j\left(\frac{1-e^{-|\xi|y}}{y}\right)\right| \lesssim |\xi|,
$$

which hold for $0 \leq \mathrm{Re}\, y \leq 1+\mu$.

For the second term on the right of (7.39), the estimates corresponding to (7.37)–(7.39) are given by the elliptic estimates in Lemma 4.3, since the map $\omega_{\mathrm{e}} \mapsto (\widetilde{u}_{\mathrm{e}}, \widetilde{v}_{\mathrm{e}})$ is the usual Biot-Savart law on $\mathbb{T}\times\mathbb{R}_+$. The estimate claimed in (7.38) for $\widetilde{v}_{\mathrm{e}}$ is immediate upon inspecting the second line in (7.36), and recalling the definitions (4.4), (4.5), (4.7). □

The estimate provided by Lemma 7.7 contains tangentially analytic norms of the trace term $g$, which we recall is given in terms of the Prandtl solution as $g(t,x) = -\int_0^\infty \partial_x \widetilde{u}^{\mathrm{P}}(x,Y,t)dY$, where $\widetilde{u}^{\mathrm{P}} = u^{\mathrm{P}} - U^{\mathrm{E}}$. However, this is precisely the term which was bounded in estimate (6.9) of Lemma 6.1. By combining these estimates we obtain:

COROLLARY 7.8. *For $s \in [0, T_*]$ and $\mu \in (0, \mu_* - \gamma_* s)$, we have*

$$\|\partial_x^i (y\partial_y)^j u_{\rm e}\|_{Y_{\lambda,\mu,\infty}} \lesssim \|\partial_x^{i+j}\omega_{\rm e}\|_{Y_{\lambda,\mu}\cap S_\mu} + j\left(\|\omega_{\rm e}\|_{Y_{\lambda,\mu}} + \|y\partial_y\omega_{\rm e}\|_{Y_{\lambda,\mu}}\right) + 1\,, \tag{7.40}$$

$$\left\|(y\partial_y)^j\left(\frac{\partial_x^i(v_{\rm e}-g)}{y}\right)\right\|_{Y_{\lambda,\mu,\infty}} \lesssim \|\partial_x^{i+1}\omega_{\rm e}\|_{Y_{\lambda,\mu}\cap S_\mu} + 1\,, \tag{7.41}$$

$$\|\partial_x^i v_{\rm e}\|_{Y_{\lambda,\mu,\infty}} \lesssim \|\partial_x^i\omega_{\rm e}\|_{Y_{\lambda,\mu}\cap S_\mu} + 1\,, \tag{7.42}$$

*for integers $i, j \geq 0$ such that $i + j \leq 1$.*

7.3.4. *Proof of Lemma 7.4, the forcing term.* In this section, we establish the $Y_{\lambda,\mu}$ and $S_\mu$ estimates for $F$ and its first order tangential and conormal derivatives, as claimed in (7.13) and (7.14). We recall that $F$ is given by (3.19), which we re-arrange by appealing to (3.7) as

$$\begin{aligned} F &= -u_{\rm e}\partial_x\omega_{\rm a} - \left(v_{\rm e}\partial_y\omega_{\rm a} + \frac{1}{\epsilon^2}g\partial_Y\Omega^{\rm P}\right) - (u_{\rm a}\partial_x + v^{\rm E}\partial_y)\omega_{\rm e} - \epsilon\left(u_{\rm e}\partial_x + (\bar v^{\rm P} + v_{\rm e})\partial_y\right)\omega_{\rm e} + (\partial_x f_2 - \partial_y f_1) \\ &= F^{(1)} + \cdots + F^{(5)}\,. \end{aligned} \tag{7.43}$$

The estimate for the last term in (7.43), namely $F^{(5)}$, was given earlier in Lemma 7.6, and these bounds are already consistent with (7.13) and (7.14). We divide this section into four steps, in which we bound $\{F^{(i)}\}_{i=1}^4$.

**Step 1. Bounding $F^{(1)}$ in (7.43).** We recall the definitions (3.7)–(3.8), which give that

$$F^{(1)} = -u_{\rm e}\left(\partial_x\omega^{\rm E} - \frac{1}{\epsilon}\partial_x\Omega^{\rm P} + \epsilon\partial_x^2\bar v^{\rm P}\right)\,.$$

We apply Lemma 4.2, the improved product estimate in Lemma 6.2 for the term containing $\partial_x\Omega^{\rm P}$, the estimates (4.13), (5.36), (6.8), (6.10), and Corollary 7.8 to obtain

$$\begin{aligned} \|F^{(1)}\|_{Y_{\lambda,\mu}} = \|u_{\rm e}\partial_x\omega_{\rm a}\|_{Y_{\lambda,\mu}} &\lesssim \|u_{\rm e}\|_{Y_{\lambda,\mu,\infty}}\left(\|\partial_x\omega^{\rm E}\|_{Y_{\lambda,\mu}} + 1 + \epsilon\|\partial_{xx}\bar v^{\rm P}\|_{P_{\lambda,\mu,\infty}}\right) \\ &\lesssim 1 + \|\omega_{\rm e}\|_{Y_{\lambda,\mu}\cap S_\mu}\,, \end{aligned} \tag{7.44}$$

where we used $\epsilon \leq 1$. The estimate for $\partial_x F^{(1)}$ is essentially the same and gives

$$\|\partial_x F^{(1)}\|_{Y_{\lambda,\mu}} \lesssim 1 + \sum_{i\leq 1}\|\partial_x^i\omega_{\rm e}\|_{Y_{\lambda,\mu}\cap S_\mu}\,. \tag{7.45}$$

Similarly, the application of $y\partial_y$ results in two terms: When this operator acts on $u_{\rm e}$ we use (7.40); on the other hand, when this operator acts on $\omega_{\rm a}$, we use that $y\partial_y = Y\partial_Y$, the identity $\partial_Y\bar v^{\rm P} = -\partial_x\widetilde u^{\rm P}$, the bounds (5.36), (6.8), (6.9), (6.10), and (6.20); in summary

$$\|y\partial_y F^{(1)}\|_{Y_{\lambda,\mu}} \lesssim 1 + \sum_{i+j\leq 1}\|\partial_x^i(y\partial_y)^j\omega_{\rm e}\|_{Y_{\lambda,\mu}\cap S_\mu}\,. \tag{7.46}$$

The above three estimates are all consistent with (7.13).

Next, we bound the $S_\mu$ norm of the first term in (7.43). For $(i,j) = (0,0)$, by appealing to Lemma 4.2 and the bounds (5.2), (5.3), (6.8), (6.10), we obtain

$$\begin{aligned} \|F^{(1)}\|_{S_\mu} &\leq \|u_{\rm e}\|_{H_x^1 L_y^\infty(y\geq 1+\mu)}\left(\left\|\partial_x\omega^{\rm E}\right\|_{S_\mu} + \left\|Y\partial_x\Omega^{\rm P}\right\|_{P_{\lambda,\mu,\infty}} + \epsilon\left\|\partial_x^2\bar v^{\rm P}\right\|_{P_{\lambda,\mu,\infty}}\right) \\ &\lesssim \sum_{i\leq 1}\left\|\partial_x^i u_{\rm e}\right\|_{L_{x,y}^\infty(y\geq 1+\mu)}\,. \end{aligned} \tag{7.47}$$

Here we have implicitly used $L_x^\infty(\mathbb{T}) \subset L_x^2(\mathbb{T})$. For $(i,j) = (1,0)$, by a similar argument, we obtain

$$\|\partial_x F^{(1)}\|_{S_\mu} \lesssim \sum_{i\leq 2}\left\|\partial_x^i u_{\rm e}\right\|_{L_{x,y}^\infty(y\geq 1+\mu)}\,. \tag{7.48}$$

Lastly, for $(i,j) = (0,1)$, we have

$$\partial_y F^{(1)} = y\partial_y(u_{\rm e}\partial_x\omega_a) = (\partial_y u_{\rm e})\partial_x\omega_a + u_{\rm e}\partial_x\partial_y\omega_a\,,$$

and thus by using the identity $\partial_Y \bar{v}^{\rm P} = -\partial_x \widetilde{u}^{\rm P}$ and a similar argument to the bound (7.47), we have

$$\begin{aligned}\|\partial_y F^{(1)}\|_{S_\mu} &\lesssim \|\partial_y u_{\rm e}\|_{H^1_x L^\infty_y(y\geq 1+\mu)} \left(\left\|\partial_x \omega^{\rm E}\right\|_{S_\mu} + \left\|Y\partial_x \Omega^{\rm P}\right\|_{P_{\lambda,\mu,\infty}} + \epsilon \left\|\partial_x^2 \bar{v}^{\rm P}\right\|_{P_{\lambda,\mu,\infty}}\right)\\ &\quad + \|u_{\rm e}\|_{H^1_x L^\infty_y(y\geq 1+\mu)} \left(\left\|\partial_x \partial_y \omega^{\rm E}\right\|_{S_\mu} + \left\|Y^2 \partial_x \partial_Y \Omega^{\rm P}\right\|_{P_{\lambda,\mu,\infty}} + \left\|\partial_x^3 \widetilde{u}^{\rm P}\right\|_{P_{\lambda,\mu,\infty}}\right)\\ &\lesssim \sum_{i+j\leq 2} \left\|\partial_x^i \partial_y^j u_{\rm e}\right\|_{L^\infty_{x,y}(y\geq 1+\mu)} . \end{aligned} \tag{7.49}$$

**Step 2. Bounding $F^{(2)}$ in (7.43).** Appealing to (3.8) and $\partial_Y \bar{v}^{\rm P} = -\partial_x \widetilde{u}^{\rm P}$, we write the second term in (7.43) as

$$\begin{aligned}-F^{(2)} = v_{\rm e} \partial_y \omega_{\rm a} + \frac{1}{\epsilon^2} g \partial_Y \Omega^{\rm P} &= v_{\rm e} \partial_y \omega^{\rm E} - \frac{1}{\epsilon^2}(v_{\rm e} - g)\partial_Y \Omega^{\rm P} - v_{\rm e} \partial_x^2 \widetilde{u}^{\rm P}\\ &= v_{\rm e} \partial_y \omega^{\rm E} - \frac{1}{\epsilon} \frac{v_{\rm e} - g}{y} Y \partial_Y \Omega^{\rm P} - v_{\rm e} \partial_x^2 \widetilde{u}^{\rm P} . \end{aligned} \tag{7.50}$$

When $(i,j) = (0,0)$, using the above decomposition, and appealing to Lemmas 4.2, 5.7, 6.1, Corollary 7.8, and Lemma 6.2 for the term containing $Y\partial_Y \Omega^{\rm P}$, we obtain

$$\begin{aligned}\|F^{(2)}\|_{Y_{\lambda,\mu}} &\lesssim \|v_{\rm e}\|_{Y_{\lambda,\mu,\infty}} \left(\|\partial_y \omega^{\rm E}\|_{Y_{\lambda,\mu}} + \epsilon \|(1+Y)^{3/2} \partial_x^2 \widetilde{u}^{\rm P}\|_{P_{\lambda,\mu,\infty}}\right) + \left\|\frac{v_{\rm e} - g}{y}\right\|_{Y_{\lambda,\mu,\infty}}\\ &\lesssim \left(\|\omega_{\rm e}\|_{Y_{\lambda,\mu}\cap S_\mu} + 1\right)(1+\epsilon) + \left(\|\partial_x \omega_{\rm e}\|_{Y_{\lambda,\mu}\cap S_\mu} + 1\right)\\ &\lesssim 1 + \sum_{i\leq 1} \|\partial_x^i \omega_{\rm e}\|_{Y_{\lambda,\mu}\cap S_\mu} . \end{aligned} \tag{7.51}$$

Applying $\partial_x^i (y\partial_y)^j = \partial_x^i (Y\partial_Y)^j$, with $i + j = 1$, to the definition of $F^{(2)}$ in (7.50), and using that $\mathbb{R}{\rm e}\, y \lesssim 1$ for $y \in \Omega_\mu$, yields a similar bound

$$\begin{aligned}\|\partial_x^i (y\partial_y)^j F^{(2)}\|_{Y_{\lambda,\mu}} &\lesssim \|\partial_x^i (y\partial_y)^j v_{\rm e}\|_{Y_{\lambda,\mu,\infty}} \left(\|\partial_y \omega^{\rm E}\|_{Y_{\lambda,\mu}} + \epsilon \|(1+Y)^{3/2} \partial_x^2 \widetilde{u}^{\rm P}\|_{P_{\lambda,\mu,\infty}}\right)\\ &\quad + \|v_{\rm e}\|_{Y_{\lambda,\mu,\infty}} \left(\|\partial_x^i \partial_y^{j+1} \omega^{\rm E}\|_{Y_{\lambda,\mu}} + \epsilon \|(1+Y)^{3/2} \partial_x^{i+2} (Y\partial_Y)^j \widetilde{u}^{\rm P}\|_{P_{\lambda,\mu,\infty}}\right)\\ &\quad + \left\|(y\partial_y)^j \left(\frac{\partial_x^i (v_{\rm e} - g)}{y}\right)\right\|_{Y_{\lambda,\mu,\infty}} + \left\|\frac{v_{\rm e} - g}{y}\right\|_{Y_{\lambda,\mu,\infty}}\\ &\lesssim \left(\|\partial_x \omega_{\rm e}\|_{Y_{\lambda,\mu}\cap S_\mu} + 1\right) + \left(\|\partial_x^{i+1} \omega_{\rm e}\|_{Y_{\lambda,\mu}\cap S_\mu} + 1\right) + \left(\|\omega_{\rm e}\|_{Y_{\lambda,\mu}\cap S_\mu} + 1\right)\\ &\lesssim 1 + \sum_{i\leq 2} \|\partial_x^i \omega_{\rm e}\|_{Y_{\lambda,\mu}\cap S_\mu} . \end{aligned} \tag{7.52}$$

Here we have used $\partial_y v_{\rm e} = -\partial_x u_{\rm e}$ and $\partial_Y \widetilde{u}^{\rm P} = \Omega^{\rm P}$. All these terms are bounded by the right side of (7.13).

Next, we bound the $S_\mu$ norm of $F^{(2)}$, as defined in (7.50). Using Lemma 4.2, Theorem 5.1, Lemma 6.1, and Corollary 7.8, we obtain

$$\begin{aligned}\|F^{(2)}\|_{S_\mu} &\lesssim \|v_{\rm e}\|_{H^1_x L^\infty_y(y\geq 1+\mu)} \left(\|\partial_y \omega_{\rm e}\|_{S_\mu} + \left\|\partial_x^2 \widetilde{u}^{\rm P}\right\|_{P_{\lambda,\mu,\infty}}\right) + \left(\|v_{\rm e}\|_{H^1_x L^\infty_y(y\geq 1+\mu)} + \|g\|_{H^1_x}\right) \left\|Y^2 \partial_Y \Omega^{\rm P}\right\|_{P_{\lambda,\mu,\infty}}\\ &\lesssim \|v_{\rm e}\|_{H^1_x L^\infty_y(y\geq 1+\mu)} \left(\|\partial_y \omega_{\rm e}\|_S + \|\partial_x \partial_y \omega_{\rm e}\|_S + \left\|\partial_x^2 \widetilde{u}^{\rm P}\right\|_{P_{\lambda,\mu,\infty}}\right)\\ &\quad + \left(\|v_{\rm e}\|_{H^1_x L^\infty_y(y\geq 1+\mu)} + \sum_\xi (1+|\xi|^2) \int_0^\infty |(\partial_x \widetilde{u}^{\rm P})_\xi| dY\right) \left\|Y^2 \partial_Y \Omega^{\rm P}\right\|_{P_{\lambda,\mu,\infty}}\\ &\lesssim 1 + \sum_{i\leq 1} \left\|\partial_x^i v_{\rm e}\right\|_{L^\infty_{x,y}(y\geq 1+\mu)} . \end{aligned} \tag{7.53}$$

The estimates for the $S_\mu$ norm of $\partial_x^i \partial_y^j F^{(2)}$ follow similarly to (7.53) by applying the Leibniz rule, resulting in

$$\sum_{i+j=1} \|\partial_x^i \partial_y^j F^{(2)}\|_{S_\mu} \lesssim 1 + \sum_{i+j\leq 2} \left\|\partial_x^i v_{\rm e}\right\|_{L^\infty_{x,y}(y\geq 1+\mu)} , \tag{7.54}$$

and thus we omit the details.

**Step 3. Bounding $F^{(3)}$ in (7.43).** Recalling (3.7), we return to the third term in (7.43), which we re-write as

$$-F^{(3)} = (u_{\rm a}\partial_x + v^{\rm E}\partial_y)\omega_{\rm e} = u^{\rm E}\partial_x\omega_{\rm e} + \widetilde{u}^{\rm P}\partial_x\omega_{\rm e} + v^{\rm E}\partial_y\omega_{\rm e}\,.$$

First, we bound the $Y_{\lambda,\mu}$ norm of $F^{(3)}$, i.e., for $(i,j)=(0,0)$. By Lemma 4.2, Lemma 5.7, and Lemma 6.1, we have

$$\begin{aligned}\|F^{(3)}\|_{Y_{\lambda,\mu}} &\lesssim \left(\|u^{\rm E}\|_{Y_{\lambda,\mu,\infty}} + \|\widetilde{u}^{\rm P}\|_{P_{\lambda,\mu,\infty}}\right)\|\partial_x\omega_{\rm e}\|_{Y_{\lambda,\mu}} + \left\|\frac{1}{y}v^{\rm E}\right\|_{Y_{\lambda,\mu,\infty}}\|y\partial_y\omega_{\rm e}\|_{Y_{\lambda,\mu}}\\ &\lesssim \|\partial_x\omega_{\rm e}\|_{Y_{\lambda,\mu}} + \|y\partial_y\omega_{\rm e}\|_{Y_{\lambda,\mu}}\,. \end{aligned}\tag{7.55}$$

Similarly, for $i+j=1$, since $\mathbb{R}\mathrm{e}\, y \lesssim 1$ for $y\in\Omega_\mu$, by Lemma 4.2, Lemma 5.7, and Lemma 6.1, we have

$$\begin{aligned}\|\partial_x^i(y\partial_y)^jF^{(3)}\|_{Y_{\lambda,\mu}} &\lesssim \left(\|\partial_x^i\partial_y^ju^{\rm E}\|_{Y_{\lambda,\mu,\infty}} + \|\partial_x^i(Y\partial_Y)^j\widetilde{u}^{\rm P}\|_{P_{\lambda,\mu,\infty}}\right)\|\partial_x\omega_{\rm e}\|_{Y_{\lambda,\mu}}\\ &\quad+\left(i\left\|\frac{1}{y}\partial_x^iv^{\rm E}\right\|_{Y_{\lambda,\mu,\infty}} + j\|\partial_xu_{\rm e}\|_{Y_{\lambda,\mu,\infty}}\right)\|y\partial_y\omega_{\rm e}\|_{Y_{\lambda,\mu}}\\ &\quad+\left(\|u^{\rm E}\|_{Y_{\lambda,\mu,\infty}} + \|\widetilde{u}^{\rm P}\|_{P_{\lambda,\mu,\infty}}\right)\|\partial_x^{i+1}(y\partial_y)^j\omega_{\rm e}\|_{Y_{\lambda,\mu}} + \left\|\frac{1}{y}v^{\rm E}\right\|_{Y_{\lambda,\mu,\infty}}\left\|\partial_x^i(y\partial_y)^{j+1}\omega_{\rm e}\right\|_{Y_{\lambda,\mu}}\\ &\lesssim \sum_{i+j\le 2}\left\|\partial_x^i(y\partial_y)^j\omega_{\rm e}\right\|_{Y_{\lambda,\mu}}\,. \end{aligned}\tag{7.56}$$

Next, we bound the $S_\mu$ norm of $F^{(3)}$ and its first tangential and conormal derivatives. When $(i,j)=(0,0)$, by Lemma 4.2, Theorem 5.1, and Lemma 6.1 we obtain

$$\begin{aligned}\|F^{(3)}\|_{S_\mu} &\lesssim \left(\left\|u^{\rm E}\right\|_{H^1_xL^\infty_y(y\ge 1+\mu)} + \left\|\widetilde{u}^{\rm P}\right\|_{P_{\lambda,\mu,\infty}}\right)\|\partial_x\omega_{\rm e}\|_{S_\mu} + \left\|v^{\rm E}\right\|_{H^1_xL^\infty_y(y\ge 1+\mu)}\|\partial_y\omega_{\rm e}\|_{S_\mu}\\ &\lesssim \sum_{i+j=1}\left\|\partial_x^i\partial_y^j\omega_{\rm e}\right\|_{S_\mu}\,. \end{aligned}\tag{7.57}$$

By a very similar argument, for $i+j=1$ we get

$$\|\partial_x^i\partial_y^jF^{(3)}\|_{S_\mu} \lesssim \sum_{i+j\le 2}\left\|\partial_x^i\partial_y^j\omega_{\rm e}\right\|_{S_\mu}\,. \tag{7.58}$$

**Step 4: Bounding $F^{(4)}$ in (7.43).** It remains to consider the fourth term in (7.43), which we recall is given by

$$F^{(4)} = -\epsilon\left(u_{\rm e}\partial_x\omega_{\rm e} + (\overline{v}^P + v_{\rm e})\partial_y\omega_{\rm e}\right)\,. \tag{7.59}$$

This term is the only one which is nonlinear in $\omega_{\rm e}$, but it has the added benefit that it has a power of $\epsilon$ as a multiplying factor. Using that (3.14) gives $\overline{v}^{\rm P}|_{Y=0} = -g$, and recalling the definition of $\overline{v}^{\rm P}$ in (2.10), we rewrite

$$\begin{aligned}\epsilon(\overline{v}^{\rm P} + v_{\rm e})\partial_y\omega_{\rm e} = (v_{\rm e} - g) + (\overline{v}^{\rm P} + g) &= \epsilon(v_{\rm e} - g)\partial_y\omega_{\rm e} - \epsilon\partial_y\omega_{\rm e}\int_0^Y \partial_x\widetilde{u}^{\rm P}dY'\\ &= \left(\epsilon\frac{v_{\rm e}-g}{y} - \frac{1}{Y}\int_0^Y\partial_x\widetilde{u}^{\rm P}dY'\right)y\partial_y\omega_{\rm e}\,. \end{aligned}\tag{7.60}$$

Using (7.59) and (7.60), we appeal to Lemma 4.2, Lemma 6.1, and Corollary 7.8, to arrive at

$$\begin{aligned}\|F^{(4)}\|_{Y_{\lambda,\mu}} &\lesssim \epsilon\|u_{\rm e}\|_{Y_{\lambda,\mu,\infty}}\|\partial_x\omega_{\rm e}\|_{Y_{\lambda,\mu}} + \left(\epsilon\left\|\frac{v_{\rm e}-g}{y}\right\|_{Y_{\lambda,\mu,\infty}} + \left\|\partial_x\widetilde{u}^{\rm P}\right\|_{P_{\lambda,\mu,\infty}}\right)\|y\partial_y\omega_{\rm e}\|_{Y_{\lambda,\mu}}\\ &\lesssim \epsilon\left(1+\|\omega_{\rm e}\|_{Y_{\lambda,\mu}\cap S_\mu}\right)\|\partial_x\omega_{\rm e}\|_{Y_{\lambda,\mu}} + \left(1+\epsilon\|\partial_x\omega_{\rm e}\|_{Y_{\lambda,\mu}\cap S_\mu}\right)\|y\partial_y\omega_{\rm e}\|_{Y_{\lambda,\mu}}\,, \end{aligned}\tag{7.61}$$

a bound which is consistent with (7.13). Similarly, for $(i,j)=(1,0)$ we get

$$\begin{aligned}\|\partial_xF^{(4)}\|_{Y_{\lambda,\mu}} &\lesssim \epsilon\left(1+\|\partial_x\omega_{\rm e}\|_{Y_{\lambda,\mu}\cap S_\mu}\right)\|\partial_x\omega_{\rm e}\|_{Y_{\lambda,\mu}} + \left(1+\epsilon\|\partial_x^2\omega_{\rm e}\|_{Y_{\lambda,\mu}\cap S_\mu}\right)\|y\partial_y\omega_{\rm e}\|_{Y_{\lambda,\mu}}\\ &\quad+\epsilon\left(1+\|\omega_{\rm e}\|_{Y_{\lambda,\mu}\cap S_\mu}\right)\|\partial_x^2\omega_{\rm e}\|_{Y_{\lambda,\mu}} + \left(1+\epsilon\|\partial_x\omega_{\rm e}\|_{Y_{\lambda,\mu}\cap S_\mu}\right)\|y\partial_y\partial_x\omega_{\rm e}\|_{Y_{\lambda,\mu}}\,. \end{aligned}\tag{7.62}$$

On the other hand, for $(i,j) = (0,1)$ we obtain

$$\begin{aligned}\|y\partial_y F^{(4)}\|_{Y_{\lambda,\mu}} &\lesssim \epsilon\left(1 + \|\omega_{\rm e}\|_{Y_{\lambda,\mu}\cap S_\mu}\right)\|y\partial_y\partial_x\omega_{\rm e}\|_{Y_{\lambda,\mu}} + \left(1+\epsilon\|\partial_x\omega_{\rm e}\|_{Y_{\lambda,\mu}\cap S_\mu}\right)\|(y\partial_y)^2\omega_{\rm e}\|_{Y_{\lambda,\mu}} \\ &\quad + \epsilon\left(1 + \|\partial_x\omega_{\rm e}\|_{Y_{\lambda,\mu}\cap S_\mu} + \|\omega_{\rm e}\|_{Y_{\lambda,\mu}} + \|y\partial_y\omega_{\rm e}\|_{Y_{\lambda,\mu}}\right)\|\partial_x\omega_{\rm e}\|_{Y_{\lambda,\mu}} \\ &\quad + \left(1+\epsilon\|\partial_x\omega_{\rm e}\|_{Y_{\lambda,\mu}\cap S_\mu}\right)\|y\partial_y\omega_{\rm e}\|_{Y_{\lambda,\mu}}\,. \end{aligned} \tag{7.63}$$

To conclude, it remains to estimate $\partial_x^i(y\partial_y)^j F^{(4)}$ with respect to the $S_\mu$ norm. For $(i,j)=(0,0)$, using (7.59), Lemma 4.2, and Lemma 6.1, we have

$$\begin{aligned}\|F^{(4)}\|_{S_\mu} &\lesssim \epsilon\,\|u_{\rm e}\|_{H^1_x L^\infty_y(y\geq 1+\mu)}\|\partial_x\omega_{\rm e}\|_{S_\mu} + \epsilon\left(\|v_{\rm e}\|_{H^1_x L^\infty_y(y\geq 1+\mu)} + \left\|\bar{v}^{\rm P}\right\|_{P_{\lambda,\mu,\infty}}\right)\|\partial_y\omega_{\rm e}\|_{S_\mu} \\ &\lesssim \epsilon\left(1 + \sum_{i\leq 1}\left\|\partial_x^i u_{\rm e}\right\|_{L^\infty(y\geq 1+\mu)} + \left\|\partial_x^i v_{\rm e}\right\|_{L^\infty(y\geq 1+\mu)}\right)\sum_{i+j=1}\left\|\partial_x^i\partial_y^j\omega_{\rm e}\right\|_{S_\mu}\,. \end{aligned} \tag{7.64}$$

The estimate for $\partial_x F^{(4)}$ is nearly identical, upon applying the Leibniz rule in $x$. For the $\partial_y F^{(4)}$ estimate, the only special term is $\partial_y\bar{v}^{\rm P}\partial_y\omega_{\rm e} = -\epsilon^{-1}\partial_x\widetilde{u}^{\rm P}\partial_y\omega_{\rm e}$, which nonetheless may be bounded using (4.11d) with $\theta = 1$. In analogy to (7.64), for $i+j=1$ the resulting estimate is

$$\|\partial_x^i\partial_y^j F^{(4)}\|_{S_\mu} \lesssim \epsilon\left(1 + \sum_{i+j\leq 2}\left\|\partial_x^i\partial_y^j u_{\rm e}\right\|_{L^\infty(y\geq 1+\mu)} + \left\|\partial_x^i\partial_y^j v_{\rm e}\right\|_{L^\infty(y\geq 1+\mu)}\right)\sum_{i+j\leq 2}\left\|\partial_x^i\partial_y^j\omega_{\rm e}\right\|_{S_\mu}\,. \tag{7.65}$$

**Step 5: Conclusion of the proof of Lemma 7.4.** By adding the upper bounds obtained in (7.44), (7.45), and (7.46) for $F^{(1)}$, the estimates (7.51) and (7.52) for $F^{(2)}$, the upper bounds (7.55) and (7.56) for $F^{(3)}$, the estimates (7.61), (7.62) and (7.63) for $F^{(4)}$, and the bound (7.22) for $F^{(5)}$, we obtain the proof of (7.13).

By adding the upper bounds obtained in (7.47), (7.48), and (7.49) for $F^{(1)}$, the estimates (7.53) and (7.54) for $F^{(2)}$, the upper bounds (7.57) and (7.58) for $F^{(3)}$, the estimates (7.64) and (7.65) for $F^{(4)}$, and the bound (7.23) for $F^{(5)}$, we obtain the proof of (7.13).

Lastly, we recall that bound (7.15) was previously established in Section 7.3.1, thereby establishing Lemma 7.4.

**7.4. Proof of Proposition 7.1.** According to definition (4.9), we fix $0 \leq t \leq T_*$ and let $\mu \in (0, \mu_* - \gamma_* t)$. Using the mild formulation (7.7), and applying Lemma 7.3, we obtain

$$\begin{aligned}\sum_{i+j=2}\|\partial_x^i(y\partial_y)^j\omega_{\rm e}(t)\|_{Y_{\lambda,\mu}} &\lesssim 1 + \int_0^t\left(\sum_{i+j\leq 1}\|\partial_x^i(y\partial_y)^j F(s)\|_{Y_{\lambda,\bar\mu}} + \sum_{i+j\leq 1}\|\partial_x^i\partial_y^j F(s)\|_{S_{\bar\mu}}\right)\frac{ds}{\mu_*-\mu-\gamma_* s} \\ &\quad + \int_0^t\left(\sum_{i\leq 1}(\|\partial_x^i F(s)\|_{Y_{\lambda,\bar\mu}} + \|\partial_x^i F(s)\|_{S_\mu}) + \sum_{i\leq 1}\sum_\xi e^{\bar\mu|\xi|}|\xi|^i|\partial_t g_\xi(s)|\right)\frac{ds}{\mu_*-\mu-\gamma_* s}\,,\end{aligned}$$

where $\bar\mu$ is as defined in (7.9). In particular, $\mu_* - \bar\mu - \gamma_* s = (3/4)(\mu_* - \mu - \gamma_* s)$. Applying Lemma 7.4 and Proposition 7.5, we deduce

$$\begin{aligned}\sum_{i+j=2}\|\partial_x^i(y\partial_y)^j\omega_{\rm e}(t)\|_{Y_{\lambda,\mu}} &\lesssim 1 + \int_0^t\left(\frac{1+|\!|\!|\omega_{\rm e}|\!|\!|_s}{\mu_*-\mu-\gamma_* s} + \frac{(1+\epsilon|\!|\!|\omega_{\rm e}|\!|\!|_s)|\!|\!|\omega_{\rm e}|\!|\!|_s}{(\mu_*-\mu-\gamma_* s)^{4/3}}\right)ds \\ &\quad + \int_0^t\frac{1+|\!|\!|\omega_{\rm e}|\!|\!|_s + \epsilon|\!|\!|\omega_{\rm e}|\!|\!|_s^2}{\mu_*-\mu-\gamma_* s}ds + \frac{1}{\gamma_*^{3/4}(\mu_*-\mu-\gamma_* t)^{1/4}}\,.\end{aligned} \tag{7.66}$$

In the above estimate, we have used the inequalities (7.15) and (7.17), applied the Hölder inequality in time, and have used the estimate

$$\int_0^t \frac{ds}{(\mu_*-\mu-\gamma_* s)^{1+\alpha}} \lesssim \frac{1}{\gamma_*(\mu_*-\mu-\gamma_* t)^\alpha}\,,$$

which holds for $\alpha \geq 0$ and $\mu < \mu_* - \gamma_* t$. Now, using the definition of $Y(t)$ norm, and fact that $(\mu_* - \mu - \gamma_* t)^{1/3} \leq (\mu_* - \mu - \gamma_* s)^{1/3}$, and the fact that $\gamma_* \geq 2$, we get

$$\begin{aligned}
&(\mu_* - \mu - \gamma_* t)^{1/3} \sum_{i+j=2} \|\partial_x^i (y\partial_y)^j \omega_{\rm e}(t)\|_{Y_{\lambda,\mu}} \\
&\lesssim 1 + \left(1 + \sup_{0\leq s\leq t} \|\!|\omega_{\rm e}|\!\|_s + \epsilon \sup_{0\leq s\leq t} \|\!|\omega_{\rm e}|\!\|_s^2\right) \left(\int_0^t \frac{ds}{(\mu_* - \mu - \gamma_* s)^{2/3}} + \int_0^t \frac{(\mu_* - \mu - \gamma_* t)^{1/3} ds}{(\mu_* - \mu - \gamma_* s)^{4/3}}\right) \\
&\lesssim 1 + \frac{1}{\gamma_*} \left(\sup_{0\leq s\leq t} \|\!|\omega_{\rm e}|\!\|_s + \epsilon \sup_{0\leq s\leq t} \|\!|\omega_{\rm e}|\!\|_s^2\right) . \qquad (7.67)
\end{aligned}$$

Similarly to the argument leading to (7.66), using that $\mu_* \lesssim 1$ we also may show that

$$\begin{aligned}
\sum_{i+j\leq 1} \|\partial_x^i (y\partial_y)^j \omega_{\rm e}(t)\|_{Y_{\lambda,\mu}} &\lesssim 1 + \int_0^t \left(\sum_{i+j\leq 1} \|\partial_x^i (y\partial_y)^j F(s)\|_{Y_{\lambda,\bar\mu}} + \sum_{i+j\leq 1} \|\partial_x^i \partial_y^j F(s)\|_{S_{\bar\mu}}\right) ds \\
&\quad + \int_0^t \left(\sum_{i\leq 1} (\|\partial_x^i F(s)\|_{Y_{\lambda,\bar\mu}} + \|\partial_x^i F(s)\|_{S_\mu}) + \sum_{i\leq 1}\sum_{\xi} e^{\bar\mu|\xi|} |\xi|^i |\partial_t g_\xi(s)|\right) ds \\
&\lesssim 1 + \int_0^t \frac{1 + \|\!|\omega_{\rm e}|\!\|_s + \epsilon \|\!|\omega_{\rm e}|\!\|_s^2}{(\mu_* - \mu - \gamma_* s)^{1/3}} ds + t^{3/4} \\
&\lesssim 1 + \frac{1}{\gamma_*} \left(\sup_{0\leq s\leq t} \|\!|\omega_{\rm e}|\!\|_s + \epsilon \sup_{0\leq s\leq t} \|\!|\omega_{\rm e}|\!\|_s^2\right) . \qquad (7.68)
\end{aligned}$$

Combining (7.67) and (7.68), taking a supremum over all $\mu \in (0, \mu_* - \gamma_* t)$, and appealing to the definition of the $Y(t)$ norm in (4.9), concludes the proof of (7.4).

## 8. The $Z$ norm estimate

In this section, we obtain a bound on the $Z$ norm, defined in (4.8), for $\omega_{\rm e}$. From (3.17), we recall that $\omega_{\rm e}$ satisfies

$$\partial_t \omega_{\rm e} - \epsilon^2 \Delta \omega_{\rm e} + \epsilon(u_{\rm e}\partial_x + v_{\rm e}\partial_y)\omega_{\rm e} + (u_a \partial_x \omega_{\rm e} + v_a \partial_y \omega_{\rm e}) + (u_{\rm e}\partial_x \omega_{\rm a} + v_{\rm e}\partial_y \omega_{\rm a}) = \widetilde{F} ,$$

where

$$\widetilde{F} = -\frac{1}{\epsilon^2} g \partial_Y \Omega^{\rm P} - \partial_y f_1 + \partial_x f_2.$$

Denote $\phi(y) = y\psi(y)$ where $\psi \in C^\infty$ is a non-decreasing function such that $\psi = 0$ for $0 \leq y \leq \frac14$ and $\psi = 1$ for $y \geq \frac12$. Observe that $\|y\omega_{\rm e}\|_{L^2(y\geq\frac12)} \leq \|\phi\omega_{\rm e}\|_{L^2}$. The function

$$Q(t) = \sum_{i+j\leq 3} \|\phi \partial_x^i \partial_y^j \omega_{\rm e}\|_{L^2}^2$$

satisfies

$$\begin{aligned}
\frac{dQ}{dt} &\lesssim \left(\epsilon^2 + \epsilon \|v_{\rm e}\|_{L^\infty(y\geq 1/4)} + \|v_a\|_{L^\infty(y\geq 1/4)}\right) Q \\
&\quad + \left(\epsilon \sum_{1\leq i+j\leq 2} \|\partial_x^i \partial_y^j \boldsymbol{u}_{\rm e}\|_{L^\infty(y\geq 1/4)} + \sum_{1\leq i+j\leq 2} \|\partial_x^i \partial_y^j \boldsymbol{u}_{\rm a}\|_{L^\infty(y\geq 1/4)}\right) Q \\
&\quad + \left(\epsilon \sum_{i+j=3} \left\|\partial_x^i \partial_y^j \boldsymbol{u}_{\rm e}(t)\right\|_{L^2(y\geq 1/4)} + \sum_{i+j=3} \left\|\partial_x^i \partial_y^j \boldsymbol{u}_{\rm a}(t)\right\|_{L^2(y\geq 1/4)}\right) \|\phi \nabla \omega_{\rm e}\|_{L^\infty(\mathbb{H})} Q^{1/2} \\
&\quad + (\epsilon^2 + \epsilon \|v_{\rm e}\|_{L^\infty(1/4\leq y\leq 1/2)} + \|v_a\|_{L^\infty(1/4\leq y\leq 1/2)}) \sum_{i+j\leq 3} \|\partial_x^i \partial_y^j \omega_{\rm e}\|_{L^2_{x,y}(1/4\leq y\leq 1/2)}^2 \\
&\quad + \sum_{0\leq i+j\leq 2} \|\partial_x^i \partial_y^j \boldsymbol{u}_{\rm e}\|_{L^\infty(y\geq 1/4)} \sum_{0\leq i+j\leq 4} \|\partial_x^i \partial_y^j \omega_{\rm a}\|_{L^2(y\geq 1/4)} \\
&\quad + \sum_{i+j=3} \left\|\partial_x^i \partial_y^j \boldsymbol{u}_{\rm e}(t)\right\|_{L^2(y\geq 1/4)} \sum_{i+j=1} \left\|\partial_x^i \partial_y^j \omega_{\rm a}(t)\right\|_{L^\infty(y\geq 1/4)} + \sum_{i+j\leq 3} \|\partial_x^i \partial_y^j \widetilde{F}\|_{L^2_{x,y}(y\geq 1/4)} Q^{1/2} , \qquad (8.1)
\end{aligned}$$

where $\boldsymbol{u}_{\mathrm{e}} = (u_{\mathrm{e}}, v_{\mathrm{e}})$ and $\boldsymbol{u}_{\mathrm{a}} = (u_{\mathrm{a}}, v_{\mathrm{a}})$. Also, by (7.1), we have

$$Q(0) \lesssim 1. \tag{8.2}$$

Our next goal is to estimate the right hand side of the inequality (8.1). First, we estimate the error velocity $\boldsymbol{u}_{\mathrm{e}}$ in terms of the error vorticity, which is needed in several terms in (8.1).

LEMMA 8.1. *For all $\delta \in (0, 1/2)$, we have*

$$\sum_{0\le i+j\le 2} \|\partial_x^i \partial_y^j \boldsymbol{u}_{\mathrm{e}}(t)\|_{L^\infty_{x,y}(y\ge\delta)} + \sum_{i+j=3} \left\|\partial_x^i \partial_y^j \boldsymbol{u}_{\mathrm{e}}(t)\right\|_{L^2_{x,y}(y\ge\delta)} \lesssim 1 + |\!|\!|\omega|\!|\!|_t \,,$$

*where the implicit constants depend on $\delta$. Also, we have the bound*

$$\sum_{0\le i+j\le 2} \|\partial_x^i \partial_y^j \omega_{\mathrm{e}}(t)\|_{L^\infty_{x,y}(\delta\le y\le 3/4)} + \sum_{i+j=3} \left\|\partial_x^i \partial_y^j \omega_{\mathrm{e}}(t)\right\|_{L^2_{x,y}(y\ge\delta)} \lesssim 1 + |\!|\!|\omega_{\mathrm{e}}|\!|\!|_t \,.$$

PROOF OF LEMMA 8.1. Recall from (3.14) that $g = -\bar{v}^{\mathrm{P}}|_{Y=0}$. By the estimate (6.10), we obtain

$$\|\partial_i g\|_{L^\infty(\mathbb{T})} \lesssim 1, \qquad i \in \mathbb{N}_0 \,, \tag{8.3}$$

where the implicit constant depends on $i$, as long as $\gamma \ge 5/2$. The rest of the proof proceeds exactly as in the proof of [**39**, Lemma 5.1]. Note that the proof depends on the Biot-Savart law (7.35)–(7.36), and the only difference between the Biot-Savart law here and in [**39**] is the presence of $g$, which is simply bounded by (8.3). □

Next, we bound the Sobolev norms of the approximate velocity $\boldsymbol{u}_{\mathrm{a}}$ and vorticity $\omega_{\mathrm{a}}$ from (3.7)–(3.8).

LEMMA 8.2. *Assume $\gamma > 5/2$. For all $\delta \in (0, 1/2)$*

$$\sum_{0\le i+j\le 3} \|\partial_x^i \partial_y^j \boldsymbol{u}_{\mathrm{a}}(t)\|_{L^\infty_{x,y}(y\ge\delta)} + \sum_{i+j=4} \left\|\partial_x^i \partial_y^j \boldsymbol{u}_{\mathrm{a}}(t)\right\|_{L^2_{x,y}(y\ge\delta)} \lesssim 1 \,, \tag{8.4}$$

*where the implicit constants depend on $\delta$. Also, we have the bound*

$$\sum_{0\le i+j\le 2} \|\partial_x^i \partial_y^j \omega_{\mathrm{a}}(t)\|_{L^\infty_{x,y}(\delta\le y\le 3/4)} + \sum_{i+j=3} \left\|\partial_x^i \partial_y^j \omega_{\mathrm{a}}(t)\right\|_{L^2_{x,y}(y\ge\delta)} \lesssim 1 \,. \tag{8.5}$$

PROOF OF LEMMA 8.2. Recall that $u_{\mathrm{a}} = u^{\mathrm{E}} + \widetilde{u}^{\mathrm{P}}$ and $v_{\mathrm{a}} = v^{\mathrm{E}} + \epsilon\bar{v}^{\mathrm{P}}$. Since (5.16) holds, in order to prove the claimed upper bound for the first term in (8.4), we only need to prove

$$\sum_{0\le i+j\le 2} \|\partial_x^i \partial_y^j \widetilde{u}^{\mathrm{P}}(t)\|_{L^\infty_{x,y}(y\ge\delta)} + \|\partial_x^i \partial_y^j \bar{v}^{\mathrm{P}}(t)\|_{L^\infty_{x,y}(y\ge\delta)} \lesssim 1. \tag{8.6}$$

Note that the bound on the first term in (6.9) implies

$$|\partial_x^i \partial_y^j \widetilde{u}^{\mathrm{P}}(Y)| \lesssim \frac{1}{Y^{j+\gamma-3/2}}, \qquad i, j \in \mathbb{N}_0 \,,$$

where the implicit constant depends on $i$ and $j$. The bound for the first term in (8.6) then holds if we assume $\gamma > 3/2$. The bound for the second term in (8.6) is the same, except that we use (6.10) instead of (6.9) and we assume $\gamma > 5/2$.

For the bound on the second term in (8.4), recall that $Y = y/\epsilon$ and thus

$$\|f(Y)\|_{L^2_y(y\ge\delta)} = \epsilon^{1/2} \|f(Y)\|_{L^2_Y(Y\ge\delta/\epsilon)} \,.$$

The bound on the first component of the velocity then holds if $\gamma > 2$ and for the second component if $\gamma > 3$.

In order to prove the estimate (8.5), we use (5.16) for the Euler part, while for the Prandtl part we have the bound

$$|\partial_x^i \partial_y^j \Omega^{\mathrm{P}}(Y)| \lesssim \frac{1}{Y^{\gamma+j}}, \qquad i, j \in \mathbb{N}_0 \,, \tag{8.7}$$

which follows from (6.8). The bound for the first term in (8.5) then holds if $\gamma \ge 1$ while the bound for the second term in (8.5) follows if $\gamma > 3/2$. □

Finally, we state the bound for the forcing term $\widetilde{F}$.

LEMMA 8.3. *Assume that $\gamma > 2$. For every $\delta > 0$, we have*

$$\|\partial_x^i \partial_y^j \widetilde{F}\|_{L^2_{x,y}(y\ge\delta)} \lesssim 1, \qquad i, j \in \mathbb{N}_0 \,,$$

*where the implicit constant depends on $i$, $j$, and $\delta > 0$.*

PROOF OF LEMMA 8.3. Observing the expansion (7.24) for $-\partial_y f_1 + \partial_x f_2$, we note that all terms contain products of Prandtl and Euler velocities and vorticities. To avoid repetition, we only estimate the higher order term, which is the first term in (7.24) and requires bounding $\epsilon^{-2}\partial_x^i \partial_y^j \int_i \Omega^{\mathrm{P}}$ in $L^2(y \geq \delta)$. Using (8.7), we get

$$\|\partial_x^i \partial_y^j \Omega^{\mathrm{P}}\|_{L^2(y\geq\delta)} \lesssim 1$$

provided $\gamma \geq 2$ since

$$\|Y^{-\gamma}\|_{L^2(y\geq\delta)} = \epsilon^{1/2}\|Y^{-\gamma}\|_{L^2(Y\geq\delta/\epsilon)} \lesssim \epsilon^{\gamma},$$

where the constant depends on $\epsilon$. □

Next, we give the bound for the $Z$ norm of the error vorticity $\omega_{\mathrm{e}}$, which we recall, cf. (4.8), is given by

$$\|\omega_{\mathrm{e}}\|_Z = \sum_{i+j\leq 3} \|\omega_{\mathrm{e}}\|_S = \sum_{i+j\leq 3} \|y\partial_x^i \partial_y^j \omega_{\mathrm{e}}\|_{L^2(y\geq\frac{1}{2})}\,.$$

PROPOSITION 8.4 (**The $Z$ norm estimate**). *Assume that $\sup_{t\in[0,T]} |\!|\!|\omega_{\mathrm{e}}|\!|\!|_t$ is finite. Then we have the bound*

$$\|\omega_{\mathrm{e}}(t)\|_Z \lesssim \left(1 + \int_0^t (1 + |\!|\!|\omega_{\mathrm{e}}(s)|\!|\!|_s)^3\, ds\right) \exp\left(C\int_0^t (1 + |\!|\!|\omega_{\mathrm{e}}(s)|\!|\!|_s) ds\right), \tag{8.8}$$

*provided $\gamma > 5/2$.*

PROOF OF PROPOSITION 8.4. Applying the bounds in Lemmas 8.1, 8.2, and 8.3 in (8.1), we get

$$\frac{dQ}{dt} \lesssim (1 + |\!|\!|\omega|\!|\!|_t)Q + (1 + |\!|\!|\omega|\!|\!|_t)^3\,.$$

Using also (8.2) and applying the Grönwall lemma, we obtain (8.8). □

## 9. Proof of Theorem 3.1

The main result of the paper follows from the definition (2.12) and the following result:

THEOREM 9.1. *Assume that the Navier-Stokes initial vorticity $\omega_0^{\mathrm{NS}}$ is given by (2.12), where the Euler initial vorticity satisfies (5.1) and the Prandtl initial vorticity satisfies (6.1), for some $\lambda_0 > 0$, independent of $\epsilon$. Moreover, assume that $\omega_{\mathrm{e}0}$ that satisfies (7.1) for some $\lambda_3, \mu_3 > 0$, independent of $\epsilon$. Then, there exists a $\gamma_* \geq 2$ sufficiently large, independent of $\epsilon$, such that with the parameters $\mu_*, T_* \in (0,1]$ defined in (7.3) we have that*

$$\sup_{t\in[0,T_*]} |\!|\!|\omega_{\mathrm{e}}(\cdot,t)|\!|\!|_t \leq C_*\,,$$

*for a constant $C_* > 0$ independent of $\epsilon$.*

PROOF OF THEOREM 9.1. Under the assumption (5.1), the Euler solution satisfies the estimates in Lemma 5.7, for suitable $(T_1, \lambda_1, \mu_1)$. Assuming (6.1), and using that the Euler trace $U^{\mathrm{E}}$ is known to be real-analytic in $x$, the Prandtl solution obeys the bounds in Lemma 6.1 for suitable $(T_2, \lambda_2, \mu_2)$. Define the parameters $\mu_*, T_*, \lambda_* \in (0,1]$ as in (7.3), and let $\gamma_* \geq 2$ be a free parameter. With these fixed parameters, define the norm $|\!|\!|\cdot|\!|\!|_t$ by (4.10).

By combining Proposition 7.1 and Proposition 8.4, and using that by (7.3) we have $T_* \lesssim \gamma_*^{-1}$, we obtain the following a priori estimate for the cumulative error vorticity:

$$\begin{aligned}
|\!|\!|\omega_{\mathrm{e}}(t)|\!|\!|_t \leq C_0 &+ \frac{C_0}{\gamma_*}\left(\sup_{0\leq s\leq t} |\!|\!|\omega_{\mathrm{e}}|\!|\!|_s + \epsilon \sup_{0\leq s\leq t} |\!|\!|\omega_{\mathrm{e}}|\!|\!|_s^2\right) \\
&+ C_0\left(1 + \frac{C_0}{\gamma_*}(1 + \sup_{0\leq s\leq t} |\!|\!|\omega_{\mathrm{e}}|\!|\!|_s)^3\right) \exp\left(\frac{C_0}{\gamma_*}(1 + \sup_{0\leq s\leq t} |\!|\!|\omega_{\mathrm{e}}|\!|\!|_s)\right),
\end{aligned} \tag{9.1}$$

for a sufficiently large constant $C_0$ which is independent of $\gamma_*$ and $\epsilon$, and for all $t \in [0, T_*]$. Moreover, due to (7.1), the definitions (4.10) and (7.3), we also have

$$|\!|\!|\omega_{\mathrm{e}0}|\!|\!|_0 \leq C_0\,,$$

by possibly enlarging the value of $C_0$. Since $\epsilon \leq 1$, we deduce that upon choosing $\gamma_* \geq 2$ to be sufficiently large, solely in terms of $C_0$, we have

$$\sup_{t\in[0,T_*]} |\!|\!|\omega_{\mathrm{e}}(t)|\!|\!|_t \leq 2C_0\,, \tag{9.2}$$

which completes the proof upon letting $C_* = 2C_0$. □

## 10. Proof of Corollary 3.2

We conclude the paper by deducing the main corollary.

PROOF OF COROLLARY 3.2. We start by proving the inequality (3.2), which in light of (3.4)–(3.5) amounts to showing that $u_{\rm e}$ and $v_{\rm e}$ are $\mathcal{O}(1)$, uniformly in $\epsilon$ with respect to the $L^\infty(\mathbb{H})$ norm.

First, by (3.1), we have

$$\|\!|\omega_{\rm e}(\cdot,t)|\!\|_t \lesssim 1, \qquad t \in [0,T_*]\,. \tag{10.1}$$

Using (7.40) with $i=j=0$, we get

$$\|u_{\rm e}\|_{Y_{\lambda_*,\mu_*,\infty}} \lesssim \|\omega_{\rm e}\|_{Y_{\lambda_*,\mu_*}\cap S_\mu} + 1 \lesssim \|\!|\omega_{\rm e}(\cdot,t)|\!\|_t + 1 \lesssim 1\,, \tag{10.2}$$

where $\lambda_*$ and $\mu_*$ are as in the beginning of Section 4.3. Similarly, the bound (7.42) with $i=0$ analogously implies

$$\|v_{\rm e}\|_{Y_{\lambda_*,\mu_*,\infty}} \lesssim 1\,.$$

Next, using (8.4) with $i=j=0$ and (10.1), we get

$$\|u_{\rm e}(t)\|_{L^\infty_{x,y}(y\geq 1/2)} + \|v_{\rm e}(t)\|_{L^\infty_{x,y}(y\geq 1/2)} \lesssim 1. \tag{10.3}$$

Combining (10.2)–(10.3), and recalling the definition (4.5), we get

$$\|u_{\rm e}(t)\|_{L^\infty_{x,y}} + \|v_{\rm e}(t)\|_{L^\infty_{x,y}} \lesssim 1,$$

and (3.2) follows.

Next, we turn to the second assertion, (3.3). Let $K\subset H$ be such that $\operatorname{dist}(K,\partial\mathbb{H}) =: d_K > 0$. The inequality (3.3) then follows from (3.2) and (3.4)–(3.5) by observing that

$$\|\widetilde{u}^{\rm P}\|_{L^\infty(Y\geq d_K/\epsilon)} + \|\bar{v}^{\rm P}\|_{L^\infty(Y\geq d_K/\epsilon)} \lesssim \epsilon\,,$$

which follows from the bounds (6.9)–(6.10), due to the fact that $\gamma$ was chosen sufficiently large. Note that the bound (3.3) is not uniform as $d_K\to 0$. □

REMARK 10.1. The conclusion of Theorem 3.1 is stronger than the fact that the vanishing viscosity limit holds with respect to the energy norm. Namely, if in addition to the assumptions of Theorem 3.1 (or Remark 3.1), we assume that the Navier-Stokes data belongs to $L^2(\mathbb{H})$, and suppose that $\lim_{\epsilon\to 0}\left\|(u_0^{\rm NS}-u_0^{\rm E}, v_0^{\rm NS}-v_0^{\rm E})\right\|_{L^2(\mathbb{H})} = 0$, then the vanishing viscosity limit holds in the energy norm:

$$\lim_{\epsilon\to 0}\sup_{t\in[0,T_0]}\|(u^{\rm NS}-u^{\rm E}, v^{\rm NS}-v^{\rm E})(\cdot,t)\|_{L^2(\mathbb{H})} = 0\,. \tag{10.4}$$

In order to verify (10.4), denote the strip $S=\{(x,y)\in\mathbb{H}: 0\leq y\leq 1\}$. By (7.40), (7.42), and (10.1), we have that

$$\|\partial_x u_{\rm e}\|_{L^2(S)} + \|\partial_x v_{\rm e}\|_{L^2(S)} \lesssim 1\,. \tag{10.5}$$

Similarly, using (5.37) and (5.38) we get

$$\|\partial_x u^{\rm E}\|_{L^2(S)} + \|\partial_x v^{\rm E}\|_{L^2(S)} \lesssim 1\,,$$

and finally, (6.9) and (6.10) give

$$\|\partial_x \widetilde{u}^{\rm P}\|_{L^2(S)} + \|\partial_x \bar{v}^{\rm P}\|_{L^2(S)} \lesssim 1\,. \tag{10.6}$$

From the inequalities (10.5)–(10.6), together with the ansatz (3.4)–(3.5), we obtain

$$\|\partial_x u^{\rm NS}\|_{L^2(S)} + \|\partial_x v^{\rm NS}\|_{L^2(S)} \lesssim 1\,,$$

uniformly in $\epsilon\in(0,1]$. Applying the criterion (2.13) in **[58]** with $\alpha = 3/4$, we conclude that (10.4) holds.

### Acknowledgments

IK was supported in part by the NSF grant DMS-1907992, while VV was supported in part by the NSF CAREER Grant DMS-1911413.

DEPARTMENT OF MATHEMATICS, UNIVERSITY OF SOUTHERN CALIFORNIA, LOS ANGELES, CA 90089
*Email address*: kukavica@usc.edu

DEPARTMENT OF MATHEMATICS, UNIVERSITY OF SOUTHERN CALIFORNIA, LOS ANGELES, CA 90089
*Email address*: tnguyen5@usc.edu

COURANT INSTITUTE OF MATHEMATICAL SCIENCES, NEW YORK UNIVERSITY, NEW YORK, NY 10012
*Email address*: vicol@cims.nyu.edu

DEPARTMENT OF MATHEMATICS, SHANGHAI JIAO TONG UNIVERSITY, SHANGHAI, 200240
*Email address*: fwang256@sjtu.edu.cn